\documentclass[twoside,11pt]{article}

\usepackage{jmlr2e}
\usepackage{amssymb,amsmath,multirow,color,tabularx,array,makecell,booktabs}

\counterwithin*{equation}{section}

\newtheorem{thm}{Theorem}[section]
\newtheorem{lem}{Lemma}[section]
\newtheorem{cor}{Corollary}[section]
\newtheorem{prop}{Proposition}[section]           
\newtheorem{defi}{Definition}[section]
\newcommand{\beq}{\begin{equation}}
\newcommand{\eeq}{\end{equation}}
\newcommand{\bed}{\begin{defi}}
\newcommand{\eed}{\end{defi}}
\newcommand{\FP}{\mathrm{FP}}
\newcommand{\FN}{\mathrm{FN}}
\newcommand{\TP}{\mathrm{TP}}

\newcommand{\Hamm}{\mathrm{Hamm}}
\newcommand{\FDR}{\mathrm{FDR}}
\newcommand{\TPR}{\mathrm{TPR}}
\newcommand{\sgn}{\mathrm{sgn}}
\newcommand{\matB}{\begin{bmatrix} 1 \;\;\; & \rho\\\rho\;\;\; & 1\end{bmatrix}}
\newcommand{\matBabs}{\begin{bmatrix} 1 \;\;\; & |\rho|\\|\rho|\;\;\; & 1\end{bmatrix}}
\newcommand{\matBa}{\begin{bmatrix} 1 \;\;\; & a\\a\;\;\; & 1\end{bmatrix}}

\usepackage{lastpage}
\jmlrheading{25}{2024}{1-\pageref{LastPage}}{9/21; Revised
11/23}{1/24}{21-1137}{Ke, Liu and Ma}
\ShortHeadings{Power Analysis of Knockoff}{Zheng Tracy Ke, Jun S. Liu and Yucong Ma}

\firstpageno{1}

\begin{document}

\title{Power of Knockoff: The Impact of Ranking Algorithm, Augmented Design, and Symmetric Statistic}

\author{\name Zheng Tracy Ke \email zke@fas.harvard.edu \\
       \addr Department of Statistics\\
       Harvard University\\
       Cambridge, MA 02138, USA
       \AND
       \name Jun S. Liu \email jliu@stat.harvard.edu \\
       \addr Department of Statistics\\
       Harvard University\\
       Cambridge, MA 02138, USA 
       \AND
       \name Yucong Ma \email yucongma@g.harvard.edu\\
       \addr Department of Statistics\\
       Harvard University\\
       Cambridge, MA 02138, USA 
       }

\editor{Ji Zhu}

\maketitle

\begin{abstract}
The knockoff filter 
is a recent false discovery rate (FDR) control method for high-dimensional linear models. We point out that knockoff has three key components: ranking algorithm, augmented design, and symmetric statistic, and each component admits multiple choices. By considering various combinations of the three components, we obtain a collection of variants of knockoff. All these variants guarantee finite-sample FDR control, and our goal is to compare their power. 
We assume a Rare and Weak signal model 
on regression coefficients and compare the power of different variants of knockoff by deriving explicit formulas of false positive rate and false negative rate.  
Our results provide new insights on how to improve power when controlling FDR at a targeted level. 
We also compare the power of knockoff with its {\it propotype} - a method that uses the same ranking algorithm but has access to an ideal threshold. The comparison reveals the additional price one pays by finding a data-driven threshold to control FDR.  
\end{abstract}


\section{Introduction} \label{sec:Intro}
We consider a linear regression model, where $y\in\mathbb{R}^n$ is the vector of responses and $X\in\mathbb{R}^{n\times p}$ is the design matrix. We assume
\beq \label{linearMod}
y=X\beta+e, \qquad X=[X_1,X_2,\ldots,X_n]'\in\mathbb{R}^{n\times p}; \ \ \beta\in \mathbb{R}^p, \ \ e \sim N(0, \sigma^2 I_n). 
\eeq
Driven by the interests of high-dimensional data analysis, we assume $p$ is large and $\beta$ is a sparse vector (i.e., many coordinates of $\beta$ are zero). Variable selection is the problem of estimating the true support of $\beta$.  
Let $\hat{S}\subset\{1,2,\ldots,p\}$ denote the set of selected variables. The false discovery rate (FDR) is defined to be 
\[
\mathbb{E}\biggl[ \frac{\#\{j: \beta_j= 0, j\in \hat{S}\}}{\#\{j: j\in\hat{S}\}\vee 1} \biggr]. 
\]
Controlling FDR is a problem of great interest. When the design is orthogonal (i.e., $X'X$ is a diagonal matrix), 
the BH-procedure \citep{benjamini1995controlling} can be employed to control FDR at a targeted level. When the design is non-orthogonal, the BH-procedure faces challenges, 
and several recent FDR control methods were proposed, such as the knockoff filter \citep{barber2015controlling}, model-X knockoff \citep{candes2018panning}, Gaussian mirror \citep{xing2019controlling}, and multiple data splits \citep{dai2020false}. All these methods are shown to control FDR at a targeted level, but their power is less studied. This paper aims to provide a theoretical understanding to the power of FDR control methods. 

We introduce a unified framework that captures the key ideas behind recent FDR control methods. Starting from the seminal work of \cite{barber2015controlling}, this framework has been implicitly used in the literature, but it is the first time that we abstract it out: 
\begin{itemize}
\item[(a)] There is a {\it ranking algorithm}, which assigns an importance metric to each variable. 
\item[(b)] An FDR control method creates an {\it augmented design matrix} by adding fake variables. 
\item[(c)] The augmented design and the response vector $y$ are supplied to the ranking algorithm as input, and the output is converted to a (signed) importance metric for each original variable through a {\it symmetric statistic}.  
\end{itemize}
The three components, ranking algorithm, augmented design, and symmetric statistic, should coordinate so that the resulting importance metrics for null variables ($\beta_j=0$) have symmetric distributions and that the importance metrics for non-null variables ($\beta_j\neq 0$) are positive with high probability. When these requirements are satisfied, one can mimic the BH procedure \citep{benjamini1995controlling} to control FDR at a targeted level. 

The choices of the three components are not unique.  For example, we may use any linear regression method $\hat{\beta}$ as the ranking algorithm, where we assign $|\hat{\beta}_j|$ as the importance metric for variable $j$. Similarly, the other two components also admit multiple choices. This leads to many different combinations of the three components. 
The literature has revealed insights on how to choose these components to get a valid FDR control, but there is little understanding on how to design them to boost power.  
The main contribution of this paper is dissecting and detailing the impact of each component on the power. 

\subsection{Main results and discoveries} \label{subsec:discoveries}

We start from the orthodox knockoff  in \cite{barber2015controlling}, which uses Lasso as the ranking algorithm, a semi-definite programming (SDP) procedure to construct the augmented design, and the signed maximum function as the symmetric statistic. We then replace each component of the orthodox knockoff by a popular alternative choice in the literature. We compare the power of the resulting variant of knockoff with the power of the orthodox one. This serves to reveal the impact of each component on power. 

Our results lead to some noteworthy discoveries: (i) For the choice of symmetric statistic, the signed maximum is better than a popular alternative - the difference statistic; (ii) For the choice of the augmented design, the SDP approach in orthodox knockoff is less favored than a recent alternative - the conditional independence approach \citep{liu2019power}; (iii) For the choice of ranking algorithm, we compare Lasso and least-squares and find that Lasso has an advantage when the signals are extremely sparse and least-squares has an advantage when the signals are moderately sparse. 

For each variant of knockoff, we also consider its {\it prototype}, which applies the ranking algorithm to the original design $X$ and selects variables by applying an ideal threshold on the importance metrics output by the ranking algorithm. 
We note that the core idea of knockoff is hinged on the other two components - augmented design and symmetric statistic, as these two components serve to find a data-driven threshold on importance metrics. Therefore, the comparison of knockoff and its prototype reveals the key difference between FDR control and variable selection - we need to pay an extra price to find a data-driven threshold. If an FDR control method is designed effectively, it should have a negligible power loss compared with its prototype. In the knockoff framework, when the design is orthogonal or blockwise diagonal and when the ranking algorithm is Lasso, we can show that knockoff (with proper choices of augmented design and symmetric statistic) indeed yields a negligible power loss compared with its own prototype. On the other hand, this is not true for a general design or when the ranking algorithm is not Lasso.

\subsection{The theoretical framework and criteria of power comparison}
Let $G=X'X\in\mathbb{R}^{p\times p}$ be the Gram matrix. 
Without loss of generality, we assume that each column of $X$ has been normalized so that the diagonal entries of $G$ are all equal to $1$.\footnote{We use a conventional normalization of $X$ in the study of Rare/Weak signal models. It is different from the standard normalization where the diagonal entries of $G$ are assumed to be $n$. We note that our $\beta$ is actually $\sqrt{n} \beta$ in the standard normalization. This is why $n$ disappears in the order of nonzero $\beta_j$.}
We study a challenging regime of ``Rare and Weak signals'' \citep{donoho2015special,jin2016rare}, where 
for some constants $\vartheta\in (0,1)$ and $r>0$, we consider settings where
\[
\#\{j: \beta_j\neq 0\}  \; \sim \; p^{1-\vartheta},\qquad
|\beta_j| \; \sim \; \sqrt{2r\log(p)} \text{ if } \beta_j\neq 0.
\]
The two parameters, $\vartheta$ and $r$, characterize the signal rarity and signal weakness, respectively. Here, $\sqrt{\log(p)}$ is the minimax order for a successful inference of the support of $\beta$ \citep{genovese2012comparison}, and the constant factor $r$ drives subtle phase transitions. 
This model is widely used in multiple testing \citep{donoho2004higher,
arias2011global,barnett2017generalized} and variable selection \citep{ji2012ups, jin2014optimality, ke2014covariance}.

The power of an FDR control method depends on the target FDR level $q$. Instead of fixing $q$, we derive a trade-off diagram between FDR and the true positive rate (TPR) as $q$ varies. This trade-off diagram provides a full characterization of power,  for any given model parameters $(\vartheta,r)$. 
We also derive a phase diagram \citep{jin2016rare} for each FDR control method. The phase diagram is a partition of the two-dimensional space $(\vartheta,r)$ into different regions, 
according to the asymptotic behavior of the Hamming error (i.e., the expected sum of false positives and false negatives). 
The phase diagram provides a visualization of power for all $(\vartheta, r)$ together.  
Both the FDR-TPR trade-off diagram and phase diagram can be used as criteria of power comparison. 
We prefer the phase diagram, because a single phase diagram covers the whole parameter range (in contrast, the FDR-TPR trade-off diagram is tied to a specified $(\vartheta, r)$). 
Throughout the paper, we use phase diagram to compare different variants of knockoff. At the same time, we also give explicit forms of false positive rate and false negative rate, from which the FDR-TPR trade-off diagram can be deduced easily.

\subsection{Related literature} \label{subsec:literature}

Power analysis of FDR control methods is a small body of literature. \cite{su2017false} set up a framework for studying the trade-off between false positive rate and true positive rate across the lasso solution path.  \cite{weinstein2017power} and \cite{weinstein2020power} extended this framework to find a trade-off for the knockoff filter, when the ranking algorithm is the Lasso and thresholded Lasso, respectively.  These trade-off diagrams are for linear sparsity (number of nonzero coefficients of $\beta$ is a constant fraction of $p$) and independent Gaussian designs ($X(i,j)$ are iid $N(0,n^{-1/2})$ variables). 
However, their analysis and results do not apply to our setting: 
In our setting, $\beta$ is much sparser, and the overall signal strength as characterized by $\|\beta\|$ is much smaller. Furthermore, we are primarily interested in correlated designs, but their study 
is mostly focused on iid Gaussian designs.

For correlated designs, \cite{liu2019power} gave sufficient and necessary conditions on $X$ such that knockoff has a full power, but they did not provide an explicit trade-off diagram. Moreover, their analysis does not apply to the orthodox knockoff but only to a variant of knockoff that uses de-biased Lasso as the ranking algorithm.   
Beyond linear sparsity, \cite{fan2019rank} studied the power of model-X knockoff for arbitrary sparsity, 
under a stronger signal strength: we assume $|\beta_j|\asymp \sqrt{\log(p)}$, while they assumed $|\beta_j|\gg \sqrt{\log(p)}$. In a similar setting, \cite{javanmard2019false} studied the power of using de-biased Lasso for FDR control. 
Our paper differs from these works because we study the regime of weaker signals and also derive the explicit FDR-TPR trade-off diagrams and phase diagrams.

\cite{wang2020power} and \cite{spector2020powerful} studied the power of model-X knockoff and conditional randomization tests. They considered linear sparsity and iid Gaussian designs, and found a disadvantage of power by constructing augmented design as in the orthodox knockoff (with least-squares as the ranking algorithm). This qualitatively agrees with some of our conclusions in Section~\ref{sec:tamperDesign}, but it is for a different setting with uncorrelated variables and linear sparsity.
We also study more variants of knockoff than those considered in aforementioned works. Recently, \cite{li2021whiteout} recast the fix-X knockoff as a conditional post-selection inference method and studied its power. 

In a sequel of papers \citep{ji2012ups, jin2014optimality, ke2014covariance}, the Rare/Weak signal model was used to study variable selection. They focused on the class of Screen-and-Clean methods for variable selection and proved its optimality under various design classes. We borrowed the notion of phase diagram from these works. 
However, they did not consider any FDR control method and their methods do not apply to knockoff.
Different from the proof techniques employed by the aforementioned work,  
our proof is based on a geometric approach, where the key is studying the geometric properties of the ``rejection region'' induced by knockoff (see Section~\ref{sec:insight} for details).

\subsection{Organization}
The remainder of this paper is organized as follows. Section~\ref{sec:FDRmethods} reviews the idea of knockoff. Section~\ref{sec:background} introduces the Rare/Weak signal model and explains how to use it as a theoretical platform to study and compare the power of FDR control methods. Sections~\ref{sec:orthogonal}-\ref{sec:RankingAlg} contain the main results, where we study the impact of symmetric statistic, augmented design, and ranking algorithm, respectively.  Section~\ref{sec:insight} sketches the proof and explains the geometrical insight behind the proof. Section~\ref{sec:simu} contains simulation results.  Section~\ref{sec:discuss} concludes with a short discussion.  Detailed proofs are relegated to the Appendix.  

\section{The knockoff filter, its variants and prototypes} \label{sec:FDRmethods}
Let us first review the orthodox knockoff filter \citep{barber2015controlling}. 
Write $G=X'X$ and let $\mathrm{diag}(s)$, with $s\in \mathbb{R}^p$, be a nonnegative diagonal matrix (to be chosen by the user) such that $\mathrm{diag}(s)\preceq 2G$. The knockoff first creates a design matrix $\tilde{X}\in\mathbb{R}^{n\times p}$ such that
\beq \label{knockoff1}
\tilde{X}'\tilde{X}=G, \qquad X'\tilde{X}=G-\mathrm{diag}(s). 
\eeq
Let $x_j$ and $\tilde{x}_j$ be the $j$th column of $X$ and $\tilde{X}$, respectively, $1\leq j\leq p$. Here, $\tilde{x}_j$ is called a {\it knockoff} of variable $j$. For any $\lambda>0$, let $\hat{\beta}(\lambda)\in\mathbb{R}^{2p}$ be the solution of Lasso \citep{tibshirani1996regression} on the expanded design matrix $[X,\tilde{X}]$ with a tuning parameter $\lambda$:
\beq \label{knockoff2}
    \hat{\beta}(\lambda)=\mathrm{argmin}_b\big\{ \|y-[X,\tilde{X}]b\|^2/2+\lambda \|b\|_1\big\}. 
\eeq
For each $1\leq j\leq p$, let $Z_j=\sup\{\lambda>0: \hat{\beta}_j(\lambda)\neq 0\}$ and $\tilde{Z}_j=\sup\{\lambda>0: \hat{\beta}_{p+j}(\lambda)\neq 0\}$. The importance of variable $j$ is measured by a {\it symmetric statistic} 
\beq \label{Wj-knockoff}
W_j=f(Z_j, \tilde{Z}_j),
\eeq
where $f(\cdot,\cdot)$ is a bivariate function satisfying $f(v,u)=-f(u,v)$. Here $\{W_j\}_{j=1}^p$ are (signed) importance metrics for variables. Under some regularity conditions, it can be shown that $W_j$ has a symmetric distribution when $\beta_j=0$ and that $W_j$ is positive with high probability when $\beta_j\neq 0$. 
Given a threshold $t>0$, the number of false discoveries is equal to 
\[
\#\{j: \beta_j=0, W_j>t\}\;\; \approx\;\; \#\{j: \beta_j=0, W_j< - t\}\;\; \approx \;\; \#\{j: W_j< - t\}, 
\]
where the first approximation is based on the symmetry of the distribution of $W_j$ for null variables and the second approximation comes from the sparsity of $\beta$. 
The right hand side gives an estimate of the number of false discoveries. Hence, a data-driven threshold to control FDR at $q$ is 
\beq \label{knockoff3}
\hat{T}(q)=\min\biggl\{ t>0: \frac{\#\{j: W_j<-t\}}{\#\{j: W_j>t\}\vee 1}\leq q \biggr\}.
\eeq
The set of selected variables is $\hat{S}(q)=\{j: W_j>\hat{T}(q)\}$. As long as $\tilde{X}$ in \eqref{knockoff1} exists, it can be shown that the FDR associated with $\hat{S}(q)$ is guaranteed to be $\leq q$. 

\subsection{Variants of knockoff}
The idea of knockoff provides a general framework for FDR control. It consists of three key components, as summarized below:  
\begin{itemize}
\item[(a)] A ranking algorithm, which takes $y$ and an arbitrary design and assigns an importance metric $Z_j$ to each variable in the design. In \eqref{knockoff2}, it uses a particular ranking algorithm based on the solution path of Lasso. 
\item[(b)] An augmented design, which is the  $n\times (2p)$ matrix $[X, \tilde{X}]$, where a knockoff $\tilde{x}_j$ is created for each original $x_j$. We supply the augmented design to the ranking algorithm to get importance metrics $Z_j$ and $\tilde{Z}_j$ for each variable $j$ and its knockoff. 
\item[(c)] A symmetric statistic $f(\cdot,\cdot)$, which combines the two importance metrics $Z_j$ and $\tilde{Z}_j$ to an ultimate importance metric $W_j$ for variable $j$. 
\end{itemize}
The choice of each component is non-unique. For (c), $f$ can be any anti-symmetric function. Two popular choices are the {\it signed maximum} statistic and the {\it difference} statistic:  
\beq \label{knockoff4}
f^{\mathrm{sgm}}(u,v)=\sgn(u-v)\cdot\max\{u,v\}, \qquad \mbox{and}\qquad f^{\mathrm{dif}}(u,v)=u-v. 
\eeq
For (b), the freedom comes from choosing $\mathrm{diag}(s)$ and constructing $\tilde{X}$.  In fact, once $\mathrm{diag}(s)$ is given, it can be shown that any $\tilde{X}$ satisfying \eqref{knockoff1} yields the same 
asymptotic performance for knockoff. 
Hence, the choice of the augmented design boils down to choosing $\mathrm{diag}(s)$. A popular option is the {\it SDP-knockoff}, which solves $\mathrm{diag}(s)$ from a semi-definite programming:
\beq\label{knockoff5}
\min \sum_j(1-s_j), \qquad\mbox{subject to}\quad 0\leq s_j\leq 1, \; \mathrm{diag}(s)\preceq 2G. 
\eeq
Another option is the {\it CI-knockoff} \citep{liu2019power}: 
\beq \label{CI-knockoff}
\mathrm{diag}(s) =c\cdot [\mathrm{diag}(G^{-1})]^{-1}, \qquad\mbox{where }c=\sup\bigl\{0<\tilde{c}\leq1: c [\mathrm{diag}(G^{-1})]^{-1}\preceq 2G \bigr\}.
\eeq
For (a), Lasso is currently used as the ranking algorithm (see \eqref{knockoff2}), but it can be replaced by other linear regression methods. Take the least-squares $\hat{\beta}^{ols}=(X'X)^{-1}X'y$ for example. We can define a ranking algorithm that outputs $|\hat{\beta}_j^{ols}|$ as the importance metric. 
If we supply the augmented design $[X, \tilde{X}]$ to this ranking algorithm, then we have
\beq\label{knockoff6}
\hat{\beta}=([X, \tilde{X}]'[X, \tilde{X}])^{-1}[X, \tilde{X}]'y.
\eeq 
We can set $Z_j=|\hat{\beta}_j|$ and $\tilde{Z}_j=|\hat{\beta}_{j+p}|$ and plug them into \eqref{Wj-knockoff}. 

\begin{figure}[tb]
\centering
\includegraphics[width=.72\textwidth]{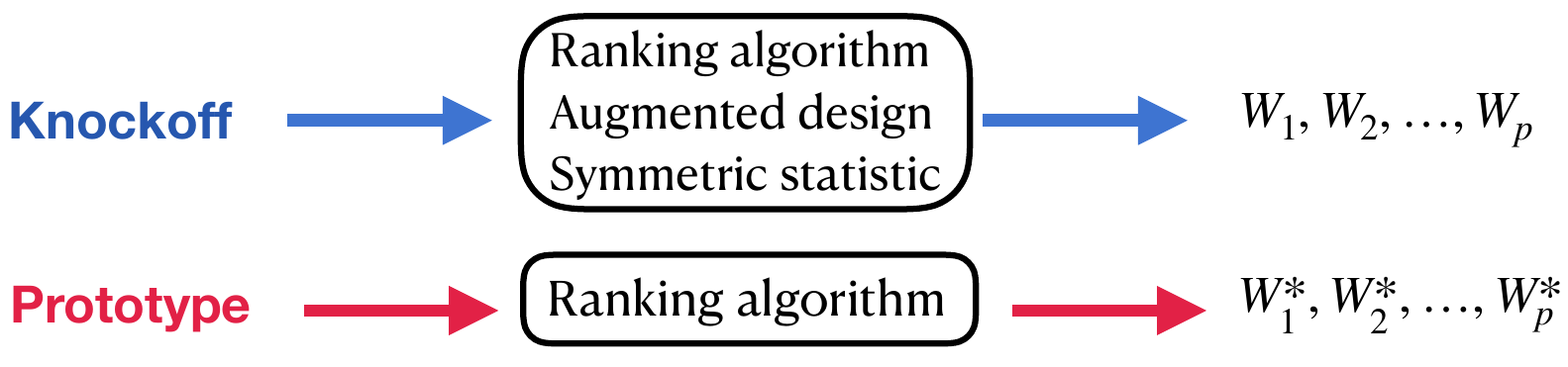}
\caption{An illustration of the knockoff and its prototype.} \label{fig:Illustration}
\end{figure}

In summary, the flexibility of the three components gives rise to many different variants of knockoff. For example, we can use \eqref{knockoff2} or \eqref{knockoff6} as the ranking algorithm, \eqref{knockoff5} or \eqref{CI-knockoff} as the augmented design, and either one in \eqref{knockoff4} as the symmetric statistic; this already gives $2\times 2\times 2=8$ variants of knockoff. 
By the theory in  \cite{barber2015controlling}, each variant guarantees the finite-sample FDR control. When the FDR is under control, the user always wants to select as many true signals as possible, i.e., to maximize the power. 
In this paper, one of our goals is to understand and compare the power of different variants of knockoff. 

To this end, we start from the default choices of the three components in the orthodox knockoff, where the ranking algorithm is Lasso as in \eqref{knockoff2}, the augmented design is the SDP-knockoff as in \eqref{knockoff5}, and the symmetric statistic is the signed maximum in \eqref{knockoff4}. In Sections~\ref{sec:orthogonal}-\ref{sec:RankingAlg}, we successively alter each component and study its impact on the power.

\subsection{Prototypes of knockoff} \label{subsec:prototype}

Given a variant of knockoff (where a specific choice of the ranking algorithm is applied to the augmented design $[X,\tilde{X}]$), we define the corresponding {\it prototype} method as follows: It runs the same ranking algorithm on the original design matrix $X$, and outputs $W_1^*, W_2^*,\ldots, W_p^*$ as importance metrics. 
The method selects variables by thresholding $W_j^*$ at $T_q^*$, where 
\beq \label{knockoff3(2)}
T_q^* =\min\biggl\{ t>0: \mathbb{E}\biggl[\frac{\#\{j: \beta_j=0, W^*_j>t\}}{\#\{j: W^*_j>t\}\vee 1}\biggr]\leq q \biggr\}.
\eeq
Compared with knockoff, the prototype ranks variables by $W_j^*$, whose induced ranking may be different from the one by $W_j$. Additionally, the prototype has access to an ideal threshold $T_q^*$
that guarantees an exact FDR control but is practically infeasible.  
In contrast,  knockoff has to find a data-driven threshold from \eqref{knockoff3}. See Figure~\ref{fig:Illustration}.

We look at two examples. Consider the orthodox knockoff, where the ranking algorithm is Lasso (see \eqref{knockoff2}). Its prototype runs Lasso on $X$ to get $\hat{\beta}^{\mathrm{lasso}}(\lambda)=\mathrm{argmin}_b \{\|y-Xb\|^2/2+\lambda\|b\|_1\}$ and assigns an importance metric to variable $j$ as
\beq \label{lasso}
W_j^*= \sup\bigl\{\lambda>0: \hat{\beta}_j^{\mathrm{lasso}}(\lambda) \neq 0\bigr\}. 
\eeq
It then selects variables by thresholding $W_j^*$ using the ideal threshold in \eqref{knockoff3(2)}. 
We call this method the {\it Lasso-path}. It is the prototype of all the variants of knockoff that use Lasso as ranking algorithm. 
For all variants of knockoff that use least-squares (see \eqref{knockoff6}) as the ranking algorithm, they share the same prototype, which
computes $\hat{\beta}^{\mathrm{ols}}=(X'X)^{-1}X'y$ and assigns an importance metric to variable $j$ as
\beq \label{ols}
W_j^* = |\hat{\beta}_j^{\mathrm{ols}}|=|e_j'G^{-1}X'y|. 
\eeq
It then selects variables by thresholding $W_j^*$ using the ideal threshold in \eqref{knockoff3(2)}. 
We call this method the {\it least-squares}.

In this paper, besides comparing different variants of knockoff, we also aim to compare each variant with its prototype. 
Here is the motivation: 
 FDR control splits into two tasks: (1) ranking the importance of variables and (2) finding a data-driven threshold. 
 When the FDR is under control, the power depends on how well Task 1 is performed. 
Importantly, in knockoff, although the augmented design and the ranking algorithm are meant to carry Task 2 only, they do affect the final ranking of variables 
 because the ranking by $W_j$'s is usually different from the ranking by $W_j^*$'s.  
This yields a potential power loss compared with its prototype -- a price we pay for finding a data-driven threshold.  Hence, a power comparison between knockoff and its prototype helps us understand how large this price is. 

\medskip

{\bf Remark 1}. In this paper, we focus on two ranking algorithms, Lasso-path and least-squares. They are both tuning-free. When the ranking algorithm has tuning parameters, we should not set the tuning parameters in the prototype the same as in the original knockoff. 
For example, we can use $|\hat{\beta}_j^{\mathrm{lasso}}(\lambda)|$ to rank variables, treating $\lambda$ as a tuning parameter. The prototype runs lasso on the original design with $p$ variables, while knockoff runs lasso on the augmented design with $2p$ variables. 
The optimal $\lambda$ that minimizes the expected Hamming error is different in the two scenarios. 
A reasonable approach of comparing knockoff and its prototype is to use their respective optimal $\lambda$. This will require computing the expected Hamming error of lasso for an arbitrary $\lambda$ \citep{ji2012ups, ke2021comparison}.

\section{Rare/Weak signal model and criteria of power comparison} \label{sec:background}
We introduce our theoretical framework of power comparison. 
Recall that we consider a linear model $y=X\beta + e$, where $y\in \mathbb{R}^n$, $X=[X_1,X_2,\ldots,X_n]' \in\mathbb{R}^{n\times p}$, and $e\sim N(0, \sigma^2 I_n)$. Without loss of generality, fix $\sigma=1$. 
Given $p$, we allow $n$ to be any integer such that $n\geq 2p$. This is from the requirement of knockoff (it needs $n\geq 2p$ to guarantee the existence of $\tilde{X}$ in \eqref{knockoff1}) and should not be viewed as a limitation of our theory. Our results are extendable to $n<2p$, provided that knockoff is replaced by its extension in this case (see Section~\ref{sec:discuss} for a discussion). 
The Gram matrix is  
\beq \label{Gram}
G: = X'X \in \mathbb{R}^{p\times p}, \qquad \mbox{where we assume } G_{jj}=1, \mbox{ for all } 1\leq j\leq p. 
\eeq
We adopt the Rare/Weak signal model \citep{donoho2004higher} to assume that $\beta$ satisfies: 
\beq\label{RWmodel1} 
\beta_j\; \overset{iid}{\sim}\; (1-\epsilon_p)\nu_0 + \epsilon_p\nu_{\tau_p}, \qquad 1\leq j\leq p,  
\eeq
where $\nu_a$ denotes a point mass at $a$. Here, $\epsilon_p\in (0,1)$ is the expected fraction of signals, and $\tau_p>0$ is the signal strength. We let $p$ be the driving asymptotic parameter and tie $(\epsilon_p, \tau_p)$ with $p$ through fixed constants $\vartheta\in (0,1)$ and $r>0$:
\beq \label{RWmodel2}
\epsilon_p = p^{-\vartheta}, \qquad\tau_p = \sqrt{2r\log(p)}. 
\eeq 
The parameters, $\vartheta$ and $r$, characterize the signal rarity and the signal strength, respectively. 
Here, $n$ does not appear in the order of nonzero $\beta_j$, because we have already re-parameterized $(X, \beta)$ such that the diagonals of $G$ are $1$ (see Footnote 1).

Under the Rare/Weak signal model \eqref{RWmodel1}-\eqref{RWmodel2}, we define two diagrams for characterizing the power of knockoff. Let $W_j$ be the ultimate importance metric \eqref{Wj-knockoff} assigned to variable $j$, and consider the set of selected variables at a threshold $\sqrt{2u\log(p)}$:
\[
\hat{S}(u)=\bigl\{1\leq j\leq p: W_j>\sqrt{2u\log(p)}\bigr\}.
\]
Let $S=\{1\leq j\leq p:\beta_j\neq 0\}$. 
Define $
\FP_p(u)= \mathbb{E}(|\hat{S}(u)\backslash S|)$, $\FN_p(u)=\mathbb{E}(|S\backslash \hat{S}(u)|)$, and $\TP_p(u)=\mathbb{E}(|S\cap \hat{S}(u)|)$, 
where the expectation is taken with respect to the randomness of both $\beta$ and $y$. Write $s_p=p\epsilon_p$, and define
\[
\Hamm_p(u) = \FP_p(u) + \FN_p(u), \quad \FDR_p(u) = \frac{\FP_p(u)}{\FP_p(u)+\TP_p(u)}, \quad \TPR_p(u)=\frac{\TP_p(u)}{s_p}.  
\]
The first quantity is the expected Hamming error. The last two quantities are proxies of the false discovery rate and true positive rate, 
respectively.\footnote{The $\FDR_p(u)$ we consider here is the ratio of expectations of false positives and total discoveries (called mFDR in some literature), not the original definition of FDR, which is the expectation of ratios of false positives and total discoveries. According to \cite{javanmard2018online}, mFDR and FDR can be different in situations with high variability, but it is not the case here. In our setting, the expectation of total discoveries grows to infinity as a power of $p$, hence, mFDR and FDR have a negligible difference.}

The following definition is conventional in the study of Rare/Weak signal models \citep{genovese2012comparison,ji2012ups} and will be used frequently in our theoretical results: 
\bed[Multi-$\log(p)$ term] \label{Multi-Log(p)}
Consider a sequence $\{a_p\}_{p=1}^\infty$. If for any fixed $\delta>0$, $a_pp^\delta\to \infty$ and $a_p p^{-\delta}\to 0$, we call $a_p$ a multi-$\log(p)$ term and write $a_p=L_p$ ($L_p$ is a generic notion for all multi-$\log(p)$ terms).  
If there is a constant $b_0\in\mathbb{R}$ such that $a_pp^{-b_0}=L_p$, we write $a_p=L_pp^{b_0}$, which means for any fixed $\delta>0$, $a_pp^{-b_0+\delta}\to\infty$ and $a_pp^{-b_0-\delta}\to 0$. 
\eed
 
In the Rare/Weak signal model, for many classes of designs of interest, $\FDR_p(u)$ and $\TPR_p(u)$ satisfy a property: There exist two fixed functions  $g_{\FDR}(u ;\vartheta,r)$ and $g_{\TPR}(u ;\vartheta,r)$ such that, for any $(\vartheta,r, u)$, as $p\to\infty$, 
\beq \label{FDR-exponent}
\FDR_p(u)=L_p p^{-g_{\FDR}(u;\vartheta,r)}, \qquad 1-\TPR_p(u)=L_pp^{-g_{\TPR}(u;\vartheta,r)}. 
\eeq
The two functions $g_{\FDR}$ and $g_{\TPR}$ depend on the choice of the three components in knockoff and the design class. We propose the FDR-TPR trade-off diagram as follows:
\bed[FDR-TPR trade-off diagram] \label{def:trade-off-diagram}
Given a variant of knockoff and a sequence of designs indexed by $p$, if $\FDR_p(u)$ and $\TPR_p(u)$ satisfy \eqref{FDR-exponent}, the FDR-TPR trade-off diagram associated with $(\vartheta, r)$ is the plot with $g_{\FDR}(u ;\vartheta,r)$ in the y-axis and $g_{\TPR}(u ;\vartheta,r)$ in the x-axis, as $u$ varies.
\eed

\begin{figure}[tb] 
\centering
\includegraphics[height=.26\textwidth]{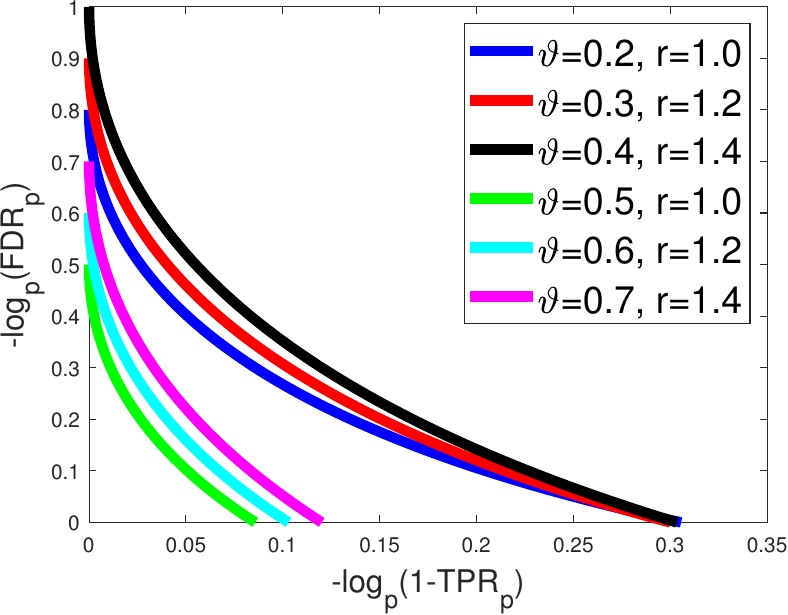}
$\qquad\quad$
\includegraphics[height=.26\textwidth]{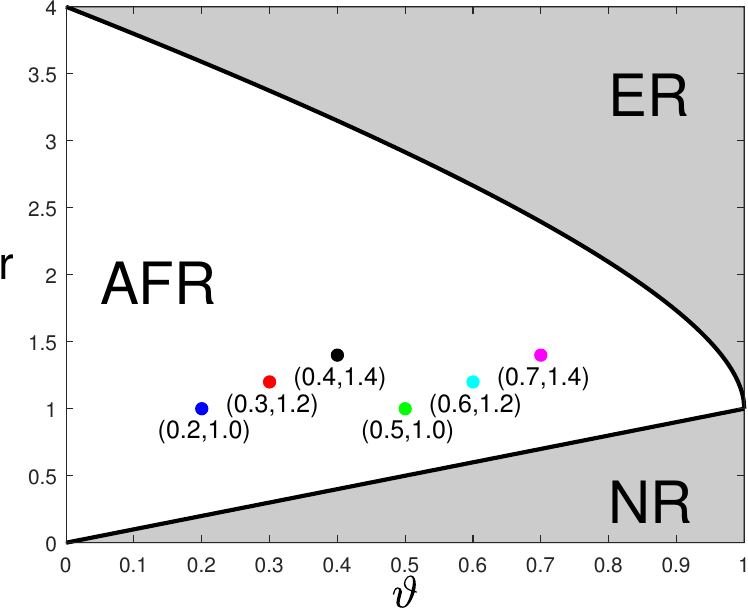} 
\caption{Left: the FDR-TPR trade-off diagram for a few values of $(\vartheta,r)$. Right: the phase diagram. The design is orthogonal, and the importance metric is as in \eqref{BH}. Each FDR-TPR trade-off diagram corresponds to one point in the phase diagram.}\label{fig:HT}
\end{figure}

An FDR-TPR trade-off diagram depends on $(\vartheta,r)$. 
To compare the performance of two variants of knockoff, we have to draw many curves for different values of $(\vartheta,r)$. 
Hence, we introduce another diagram, which characterizes the power simultaneously at all $(\vartheta,r)$. Define $\Hamm_p^*\equiv \min_u\{\FP_p(u)+\FN_p(u)\}$. This is the minimum expected Hamming selection error when the threshold $u$ is chosen optimally. For each variant of knockoff and each class of designs of interest, there exists a bivariate function $f^*_{\Hamm}(\vartheta,r)\in [0,1]$ such that
\beq \label{Hamm-exponent}
\Hamm_p^*=L_p p^{f^*_{\Hamm}(\vartheta,r)}.
\eeq
The phase diagram is defined as follows: 
\bed[Phase diagram] \label{def:PD}
When $\Hamm^*_p$ satisfies \eqref{Hamm-exponent}, the phase diagram is defined to the partition of the two-dimensional space $(\vartheta,r)$ into three regions: 
\begin{itemize} \itemsep 0pt
\item Region of Exact Recovery (ER): $\{(\vartheta,r): f^*_{\Hamm}(\vartheta,r)<0\}$.
\item Region of Almost Full Recovery (AFR): $\{(\vartheta,r): 0<f^*_{\Hamm}(\vartheta,r)<1-\vartheta\}$. 
\item Region of No Recovery (NR): $\{(\vartheta,r): f_{\Hamm}^*(\vartheta,r)\geq 1-\vartheta\}$. 
\end{itemize}
The curves separating different regions are called phase curves. We use $h_{\mathrm{AFR}}(\vartheta)$ to denote the curve between NR and AFR, and  $h_{\mathrm{ER}}(\vartheta)$ the curve between AFR and ER.  
\eed

In the ER region, the expected Hamming error, $\Hamm_p^*$, tends to zero. Therefore, with high probability, the support of $\beta$ is exactly recovered. In the AFR region, 
$\Hamm_p^*$ does not tend to zero but is much smaller than $p\epsilon_p$ (which is the expected number of signals). As a result, with high probability, the majority of signals are correctly recovered. In the region of NR, $\Hamm_p^*$ is  comparable with the number of signals, and variable selection fails. 
The phase diagram was introduced in the literature \citep{genovese2012comparison,ji2012ups} but has never been used to study FDR control methods.

We illustrate these definitions with an example. Both the FDR-TPR trade-off diagram and phase diagram only depend on the importance metrics assigned to variables. Therefore, they are also well-defined for the prototypes in Section~\ref{subsec:prototype}. We consider a special class of designs, where $X'X=I_p$, and a prototype that assigns the importance metrics
\beq \label{BH}
W^*_j = |x_j'y|, \qquad 1\leq j\leq p. 
\eeq
The next proposition is adapted from literature \citep{donoho2004higher,ji2012ups} and proved in the Appendix. We use $a_+$ to denote $\max\{a,0\}$, for any $a\in\mathbb{R}$.
\begin{prop} \label{prop:Identity}
Suppose $X'X=I_p$ and consider the importance metric in \eqref{BH}.  When $r>\vartheta$, the FDR-TPR trade-off diagram is given by $g_{\FDR}(u;\vartheta,r)=(u-\vartheta)_+$ and $g_{\TPR}(u;\vartheta,r)=(\sqrt{r}-\sqrt{u})_+^2$. The phase diagram is given by $h_{\mathrm{AFR}}(\vartheta)=\vartheta$ and $h_{\mathrm{ER}}(\vartheta)=(1+\sqrt{1-\vartheta})^2$.
\end{prop}
\noindent
These diagrams are visualized in Figure~\ref{fig:HT}. 

As we have mentioned in Section~\ref{subsec:prototype}, FDR control is composed of the task of ranking variables and the task of finding a data-driven threshold. 
It is the variable ranking that determines the power of an FDR control method simultaneously at all FDR levels $q$. The two diagrams in Definitions~\ref{def:trade-off-diagram}-\ref{def:PD} only depend on the importance metrics assigned to variables, hence, they measure the {\it quality of variable ranking}, which is fundamental for power comparison of different FDR control methods when they all control FDR at the same level.

%

\section{Impact of the symmetric statistic} \label{sec:orthogonal}
We fix the choice of ranking algorithm and augmented design in the orthodox knockoff and compare using the two symmetric statistics in \eqref{knockoff4}. 
For simplicity, we consider the orthogonal design where $X'X=I_p$. In this special case, the ranking algorithm in \eqref{knockoff2} reduces to calculating marginal regression coefficients and the output $Z_j$ and $\tilde{Z}_j$ reduce to
\beq \label{Zj-ortho}
Z_j = |x_j'y|, \qquad\mbox{and}\qquad \tilde{Z}_j= |\tilde{x}_j'y|, \qquad 1\leq j\leq p. 
\eeq
In the augmented design, the choice of $\mathrm{diag}(s)$ in \eqref{knockoff5} reduces to $\mathrm{diag}(s)=I_p$. We consider a slightly more general form: 
\beq \label{diag(s)-ortho}
\mathrm{diag}(s)=(1-a)I_p, \qquad\mbox{where }-1< a< 1 \mbox{ is a fixed constant}. 
\eeq
Given $\mathrm{diag}(s)$, the construction of $\tilde{X}$ is not unique, but all constructions lead to exactly the same FDR-TPR trade-off diagram and the same phase diagram. For this reason, we only specify $\mathrm{diag}(s)$ as in \eqref{diag(s)-ortho}, but not the actual $\tilde{X}$. 
Fixing the above choices \eqref{Zj-ortho}-\eqref{diag(s)-ortho}, we consider the two symmetric statistics in \eqref{knockoff4}, which lead to the importance metrics of
\beq \label{Wj}
W_j^{\mathrm{sgm}} = (Z_j \vee \tilde{Z}_j)\cdot
\begin{cases}
+1,&\text{if } Z_j>\tilde{Z}_j\\
-1, &\text{if } Z_j\leq \tilde{Z}_j
\end{cases}, 
\qquad\mbox{and}\qquad W_j^{\mathrm{dif}} = Z_j-\tilde{Z}_j. 
\eeq
We call the two variants of knockoff {\it knockoff-sgm} and {\it knockoff-diff}, respectively. 
The next theorem gives the explicit forms of $\FP_p(u)$ and $\FN_p(u)$ associated with these two variants. Its proof can be found in the Appendix.

\begin{thm} \label{thm:knockoff}
Consider a linear regression model where \eqref{Gram}-\eqref{RWmodel2} hold. Suppose $n\geq 2p$ and $G=I_p$. We construct $\tilde{X}$ in knockoff as in \eqref{diag(s)-ortho}, for a constant $a\in (-1,1)$, and let $Z_j$ and $\tilde{Z}_j$ be as in \eqref{Zj-ortho}. 
For any constant $u>0$, let $\FP_p(u)$ and $\FN_p(u)$ be the expected numbers of false positives and false negatives, by selecting variables with $W_j>\sqrt{2u\log(p)}$.  
When $W_j$ is the signed maximum statistic in \eqref{Wj}, as $p\to\infty$,  
\[
\FP_p(u) = L_p p^{1-u}, \qquad \FN_p(u) = L_p p^{1-\vartheta-\min\bigl\{\frac{(1-|a|)r}{2},\ (\sqrt{r}-\sqrt{u})_+^2\bigr\}}.
\]
When $W_j$ is the difference statistic in \eqref{Wj}, as $p\to\infty$, 
\[
\FP_p(u) = L_p p^{1-u}, \qquad \FN_p(u) = L_p p^{1-\vartheta-\frac{(1-|a|)}{2}(\sqrt{r}-\sqrt{u})_+^2}.
\]
\end{thm}
Here, $L_p$ is the generic multi-$\log(p)$ notion in Definition~\ref{Multi-Log(p)}. 
For Theorem~\ref{thm:knockoff}, using Mills' ration, we actually know that $L_p\asymp 1/\sqrt{\log(p)}$. 
For other theorems, we do not always know the exact order of $L_p$, but we can intuitively regard it as a polynomial of $\log(p)$. 

Given Theorem~\ref{thm:knockoff} and Definitions~\ref{def:trade-off-diagram}-\ref{def:PD}, we can derive the explicit FDR-TPR tradeoff diagrams and the phase diagrams:

\begin{cor} \label{cor:orthoPD-knockoff}
In the setting of Theorem~\ref{thm:knockoff}, when $r>\vartheta$, 
the FDR-TPR trade-off diagram is given by \[
g_{\FDR}(u;\vartheta,r)=(u-\vartheta)_+, \quad g_{\TPR}(u)=
\begin{cases}
\min\bigl\{\frac{(1-|a|)r}{2},\; (\sqrt{r}-\sqrt{u})_+^2\bigr\}, & \mbox{if}\;\; W_j=W_j^{\mathrm{sgn}},\\
\frac{(1-|a|)}{2}(\sqrt{r}-\sqrt{u})_+^2& \mbox{if}\;\; W_j=W_j^{\mathrm{dif}}.
\end{cases}
\]
The phase diagram is given by 
\[
h_{\mathrm{AFR}}(\vartheta)=\vartheta, \qquad h_{\mathrm{ER}}(\vartheta) = \begin{cases} \max\bigl\{\frac{2-2\vartheta}{1-|a|},\; (1+\sqrt{1-\vartheta})^2\bigr\}, & \mbox{if}\;\; W_j=W_j^{\mathrm{sgn}}, \\
\Bigl(1+\sqrt{\frac{2-2\vartheta}{1-|a|}}\Bigr)^2, & \mbox{if}\;\; W_j=W_j^{\mathrm{dif}}. \\
\end{cases}
\]
\end{cor}

\begin{figure}[tb]
\centering
\includegraphics[height=.258\textwidth]{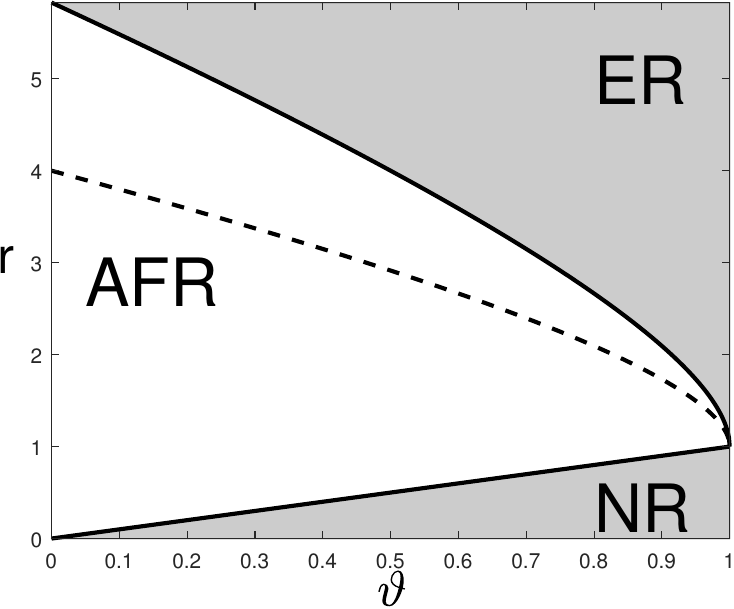}
\includegraphics[height=.26\textwidth]{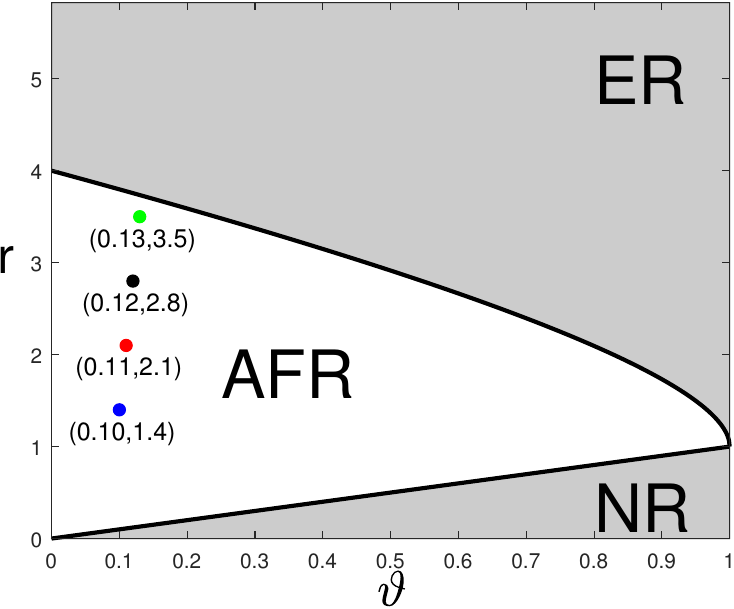}
\includegraphics[height=.26\textwidth]{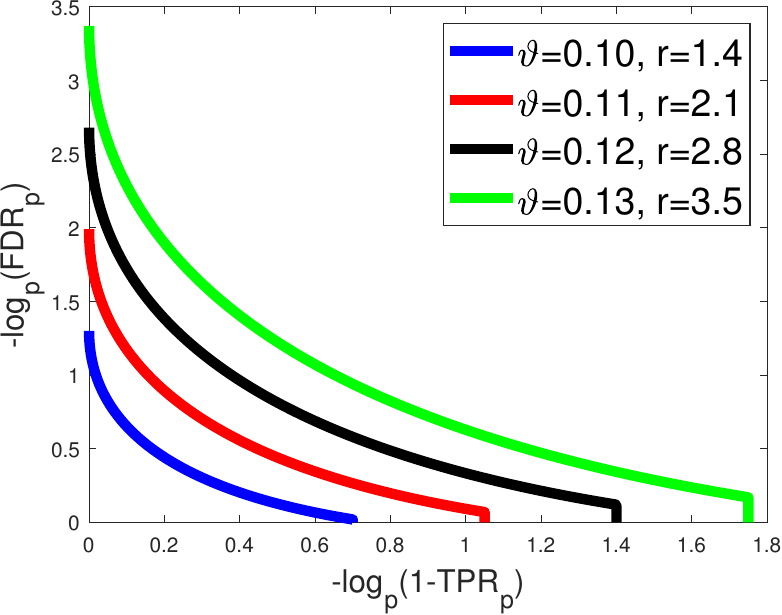}
\caption{Power comparison of knockoff with different symmetric statistics (orthogonal design; ranking algorithm is Lasso, and augmented design is such that $\mathrm{diag}(s)=I_p$). The left two panels are the phase diagrams of knockoff-diff (left) and knockoff-sgm (middle), where the dashed lines are the phase curves of their common prototype. 
The right panel is the FDR-TPR trade-off diagram of knockoff-sgm, where 
each trade-off curve corresponds to one point in the phase diagram.} \label{fig:FDR-orthogonal}
\end{figure}

For both knockoff-sgm and knockoff-diff, by Corollary~\ref{cor:orthoPD-knockoff}, the best choice of $a$ is $a=0$.
In the remaining of this section, we fix $a=0$. Figure~\ref{fig:FDR-orthogonal} gives visualizations of the FDR-TPR trade-off diagrams and phase diagrams.
The prototype is Lasso-path (see \eqref{lasso}). In the orthogonal design $X'X=I_p$, Lasso-path reduces to the  prototype in \eqref{BH}, whose FDR-TPR trade-off diagram and phase diagram are given in Figure~\ref{fig:HT}. 

\paragraph{Comparison of two symmetric statistics:}
First, we compare the phase diagrams in Figures~\ref{fig:HT}-\ref{fig:FDR-orthogonal} and find that (i) knockoff-sgm has a strictly better phase diagram than knockoff-diff, and (ii) knockoff-sgm has the same phase diagram as the prototype.
It suggests that signed maximum is a better choice of symmetric statistic. It also suggests that knockoff-sgm yields a negligible power loss relative to its prototype. 

We also point out that knockoff-sgm is already ``optimal" among all symmetric statistics, in this orthogonal design. The reason is that, when $X'X=I_p$, the Hamming error ($\Hamm_p^*$) has an information-theoretical lower bound \citep{genovese2012comparison,ji2012ups}, whose induced phase diagram coincides with the phase diagram of the prototype, which is also the phase diagram of knockoff-sgm. This is the optimal phase diagram any method can achieve (including all variants of knockoff with other symmetric statistics).

Next, we compare the FDR-TPR trade-off diagrams of knockoff and the prototype.  We focus on knockoff-sgm, whose trade-off diagram is in Figure~\ref{fig:FDR-orthogonal} (right panel). The trade-off diagram of the prototype is in Figure~\ref{fig:HT} (right panel). We find that the trade-off diagram of knockoff-sgm is slightly different from the one of the prototype. By Theorem~\ref{thm:knockoff}, $(1-\TPR_p)=\FN_p/s_p\geq L_pp^{-r/2}$; hence, the FDR-TRP trade-off curve is truncated at $r/2$ in the x-axis. For large $\vartheta$, the curve hits zero before the x-axis reaches $r/2$, and the truncation has no impact. However, for small $\vartheta$, the curve has changed due to the truncation  (Figure~\ref{fig:FDR-orthogonal}, right panel, all but the blue curve).


\paragraph{Some geometric insights, especially why signed maximum is ``optimal".}

\begin{figure}[tb]
\centering
\includegraphics[height=.26\textwidth]{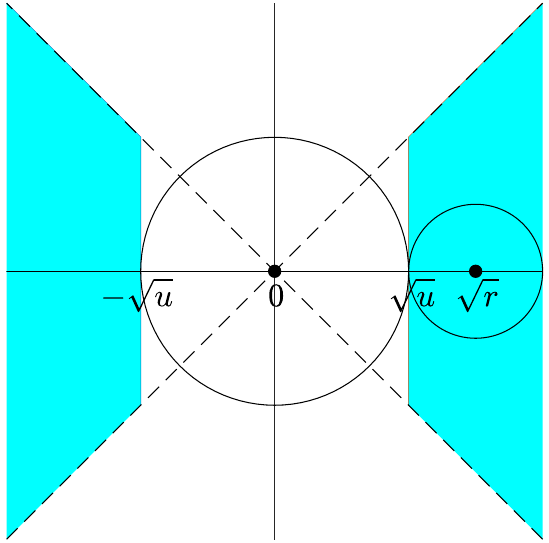} $\qquad$
\includegraphics[height=.26\textwidth]{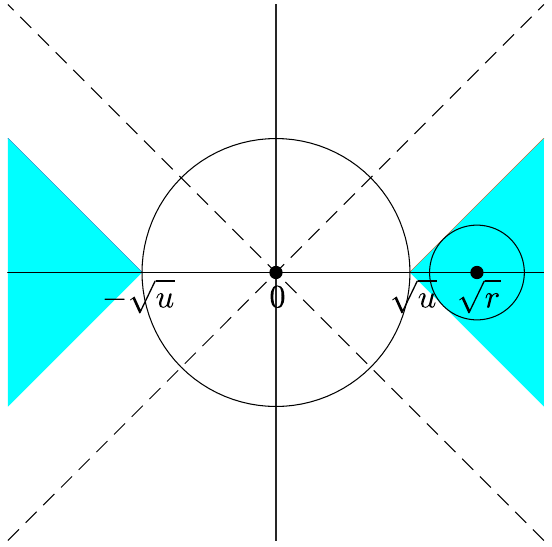} $\qquad$
\includegraphics[height=.26\textwidth]{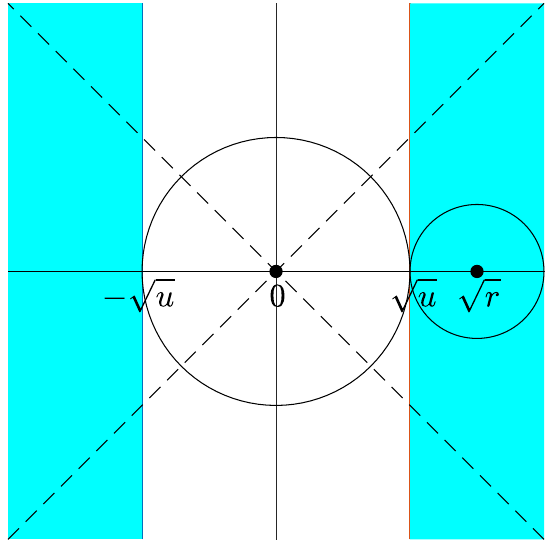}
\caption{The rejection region of the symmetric statistics (orthogonal design, $a=0$ in the construction of knockoff variables). Left: the signed maximum statistic. Middle: the difference statistic. Right: the prototype.} 
\label{fig:rejection-region}
\end{figure}

By \eqref{Zj-ortho} and \eqref{Wj}, 
the importance metrics produced by knockoff can be written as $W_j=I(x_j'y, \tilde{x}_j'y)$, where $x_j$ and $\tilde{x}_j$ are the $j$th variable and its knockoff, and $I(\cdot,\cdot)$ is a bivariate function depending on the choice of symmetric statistic. Define the ``rejection region'' as
\[
{\cal R}=\left\{(h_1, h_2)\in\mathbb{R}^2:\; I\Bigl( h_1\sqrt{2\log(p)},\, h_2\sqrt{2\log(p)}\Bigr)>\sqrt{2u\log(p)}\right\}. 
\]
Figure~\ref{fig:rejection-region} shows the rejection region induced by knockoff-sgm, knockoff-diff, and their prototype. Write $\hat{h}_1=x_j'y/\sqrt{2\log(p)}$ and $\hat{h}_2=\tilde{x}_j'y/\sqrt{2\log(p)}$. 
The random vector $(\hat{h}_1,\hat{h}_2)'$ follows the bivariate normal distribution with covariance matrix $\frac{1}{\log(p)}I_2$. Its mean vector is $(0,0)'$ when $\beta_j=0$ and $(\sqrt{r},0)'$ when $\beta_j=\tau_p$.
By Lemma~\ref{lem:tool} (to be introduced in Section~\ref{sec:insight}), 
the exponent in $\FP_p$ is determined by the Euclidean distance from $(0,0)'$ to ${\cal R}$ and the exponent in $\FN_p$ is determined by the Euclidean distance from $(\sqrt{r}, 0)'$ to ${\cal R}^c$. 
From Figure~\ref{fig:rejection-region}, it is clear that the difference statistic is inferior to the signed maximum statistic because the distance from $(\sqrt{r}, 0)'$ to ${\cal R}^c$ is strictly smaller in the former.

As we have mentioned earlier, signed maximum is ``optimal" among all symmetric statistics, because its phase diagram already matches with the information-theoretic lower bound. Now, we use Figure~\ref{fig:rejection-region} to provide a geometric interpretation of why signed maximum is ``optimal". 
We call a subset ${\cal R}$ an eligible rejection region if there exists a symmetric statistic $f(\cdot,\cdot)$ whose induced rejection region is ${\cal R}$. 
It is not hard to see that any eligible ${\cal R}$ should be symmetric with respect to both x-axis and y-axis. In addition, by the anti-symmetry requirement $f(v,u)=-f(u,v)$, an eligible rejection region also needs to satisfy the following necessary condition: {\it ${\cal R}\cap {\cal R}_{\pm}=\emptyset$, where ${\cal R}_{\pm}$ is the reflection of ${\cal R}$ with respect to the line $y=\pm x$}. 
The prototype has the optimal phase diagram, but its rejection region ${\cal R}_0$ (Figure~\ref{fig:rejection-region}, right panel) does not satisfy this condition. We can see that the rejection region of knockoff-sgm (left panel) is a {\it minimal} modification of ${\cal R}_0$ to tailor to this condition.
This partially explains why signed maximum is already the best choice in this setting. 

\smallskip

{\bf Remark 2}. We consider a non-stochastic threshold $\sqrt{2u\log(p)}$ in Theorem~\ref{thm:knockoff}. For a data-driven threshold $\sqrt{2\hat{u}\log(p)}$, if there is $\epsilon_p=o(1)$ such that $|\hat{u}-u|\leq \epsilon_p$ with probability $1-o(p^{-2})$, then Theorem~\ref{thm:knockoff} continues to hold. 

{\bf Remark 3}. In this section, we conduct power comparison only on the orthogonal design. Our rationale is as follows: A good method has to at least perform well in the simplest case. If a method is inferior to others in the orthogonal design, then we do not expect it to have a good potential in real applications where the designs can be much more complicated, i.e., the results on orthogonal designs help us filter out those methods that have little potential in practice. Such insights are valuable to users.

{\bf Remark 4}. Recently, \cite{weinstein2020power} showed a remarkable result: Under linear sparsity and random Gaussian design, they found that knockoff with the difference symmetric statistic has great power. We note that the prototype of their knockoff is {\it thresholded Lasso}, not Lasso, and so the gain of power is primarily from prototype instead of symmetric statistic. Also, see  \cite{ke2021comparison} for a comparison of thresholded-Lasso v.s. Lasso.

\section{Impact of the augmented design} \label{sec:tamperDesign}
We fix the choice of ranking algorithm and symmetric statistic as in the orthodox knockoff and compare using two augmented designs, the SDP-knockoff in \eqref{knockoff5} and the CI-knockoff in \eqref{CI-knockoff}. We also call the two respective variants of knockoff the {\it SDP-knockoff} and {\it CI-knockoff}, so that  ``SDP-knockoff" (say) has two meanings, an augmented design or a variant of knockoff, depending on the context.

When $X'X=I_p$, SDP-knockoff and CI-knockoff are the same and reduce to $\mathrm{diag}(s)=I_p$, so it is impossible to tell their difference in power.  
We must consider non-orthogonal designs. 
However, since there is no explicit form of the Lasso solution path, the results for a general $X$ are difficult to obtain;
despite technical challenges, the phase diagrams may be too messy to provide any useful insight. We hope to find a class of non-orthogonal designs such that (i) it is mathematically tractable, (ii) it is considerably different from the orthogonal design and allows $G$ to have some ``large" off-diagonal entries, 
and (iii) it captures some key features of real applications. 
We start from a class of row-wise sparse designs, which approximate the designs in many real applications (e.g., in bioinformatics and in compressed sensing). The next proposition is adapted from Lemma 1 of \cite{jin2014optimality}, whose proof is omitted.

\begin{prop} \label{prop:graph}
Consider a linear model where \eqref{Gram}-\eqref{RWmodel2} hold. Suppose each row of $G$ has at most $L_p$ nonzero entries, where $L_p$ is a multi-$\log(p)$ term as in Definition~\ref{Multi-Log(p)}. Let $S$ be the support of $\beta$. There exists a constant integer $m_0=m_0(\vartheta)$ such that with probability $1-o(1)$, $G_{SS}$ is a blockwise diagonal matrix afte a permutation of indices, and the maximum block size is bounded by $m_0$.
\end{prop}

\noindent
Proposition~\ref{prop:graph} is a consequence of the interplay between design sparsity and signal sparsity: Under the Rare/Weak signal model \eqref{RWmodel1} and a sparse design, the true signals in $S$ appear in groups, where each group contains only a small number of variables and distinct groups are mutually uncorrelated.  
This motivates us to consider a simpler setting where $G$ is blockwise diagonal {\it by itself}.
While these two settings look so different from each other, 
the asymptotic behavior of Hamming error is closely related. For example, when $G$ is a tridiagonal matrix with equal values in the sub-diagonal, the optimal phase diagram is the same as in the case where $G$ is blockwise diagonal with $2\times 2$ blocks \citep{ji2012ups,jin2014optimality}.  Inspired by these observations, we study a class of blockwise diagonal designs \citep{jin2016rare}:
For some $\rho\in (-1,1)$ and a $p\times p$ permutation matrix ${\cal T}$,
\begin{equation}\label{block}
G={\cal T}\mathrm{diag}(B, B, \ldots, B, B_1){\cal T},\quad\mbox{where}\;   B= \begin{bmatrix} 1 \;\;\; & \rho\\\rho\;\;\; & 1\end{bmatrix},  \;\; B_1= \begin{cases} B, & \mbox{if $p$ is even}, \\ 1, &\mbox{if $p$ is odd}. \end{cases} 
\end{equation}
This design serves the aforementioned purposes (i)-(iii): It has only one parameter $\rho$, so is mathematically tractable. 
The nonzero off-diagonal entries of $G$ are at the constant order, so this design is sufficiently different from the orthogonal design (in contrast, many literature consider the independent random Gaussian design, for which the maximum absolute off-diagonal entry of $G$ is only $o(1)$). Also, as we have argued above, studying this design helps us draw useful insights that will likely continue to hold for general sparse designs.

\subsection{The prototype, Lasso-path} \label{subsec:Lasso-path}
Before studying SDP-knockoff and CI-knockoff, we first study their prototype, Lasso-path. The next theorem characterizes $\FP_p(u)$ and $\FN_p(u)$ for Lasso-path and is proved in the Appendix.

\begin{thm} \label{thm:lasso}
Consider a linear regression model where \eqref{Gram}-\eqref{RWmodel2} hold. Suppose $n\geq 2p$ and $G$ is as in \eqref{block} with a correlation parameter $\rho\in (-1,1)$.  Let $W_j^*$ be as in \eqref{lasso}. 
For any constant $u>0$, let $\FP_p(u)$ and $\FN_p(u)$ be the expected numbers of false positives and false negatives, by selecting variables with $W^*_j>\sqrt{2u\log(p)}$. As $p\to\infty$, 
\[
\FP_p(u) = L_p p^{1-\min\left\{u,\;\, \vartheta+(\sqrt{u}-|\rho|\sqrt{r})^2+(\xi_\rho\sqrt{r}-\eta_\rho \sqrt{u})_+^2-(\sqrt{r}-\sqrt{u})_+^2\right\}},
\]
and
\[
\FN_p(u) = \begin{cases}
L_p p^{1- \vartheta - \{(\sqrt{r}-\sqrt{u})_+ - [(1-\xi_\rho)\sqrt{r}-(1-\eta_\rho)\sqrt{u}]_+\}^2},& \rho\geq 0,\\
L_p p^{1- \min\left\{ \vartheta+\{(\sqrt{r}-\sqrt{u})_+ - [(1-\xi_\rho)\sqrt{r}-(1-\eta_\rho)\sqrt{u}]_+\}^2,\;\, 2\vartheta + (\xi_\rho\sqrt{r}-\eta_\rho^{-1} \sqrt{u})_+^2\right\}},  &\rho<0,
\end{cases}
\]
where $\xi_\rho=\sqrt{1-\rho^2}$ and $\eta_\rho=\sqrt{(1-|\rho|)/(1+|\rho|)}$. 
\end{thm}

\begin{figure}[tb]
\centering
\includegraphics[height=.25\textwidth]{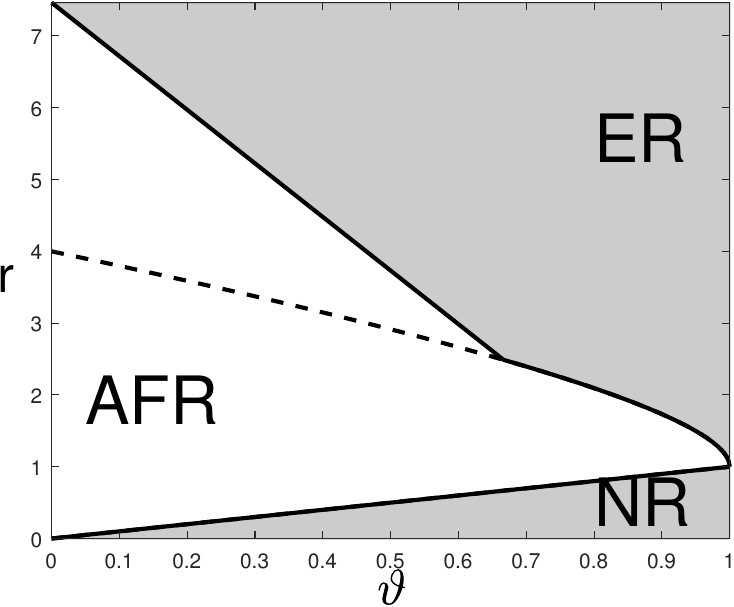}
\includegraphics[height=.25\textwidth]{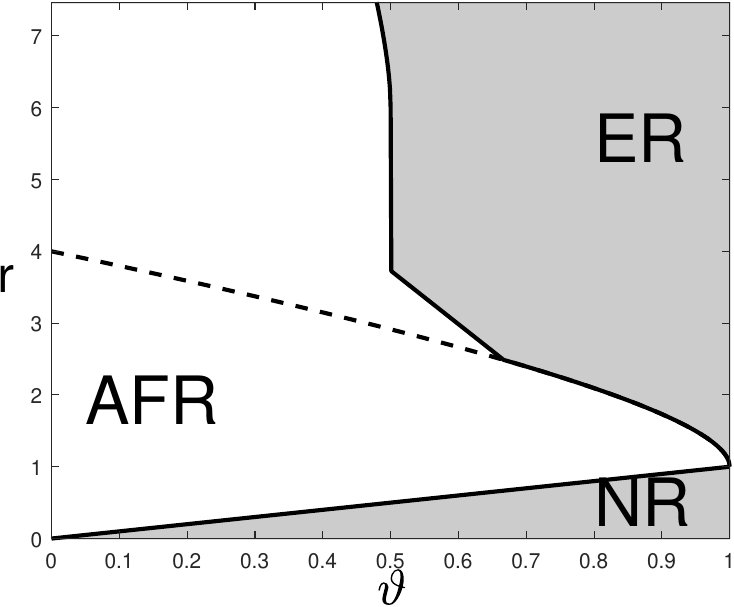}
\includegraphics[height=.25\textwidth]{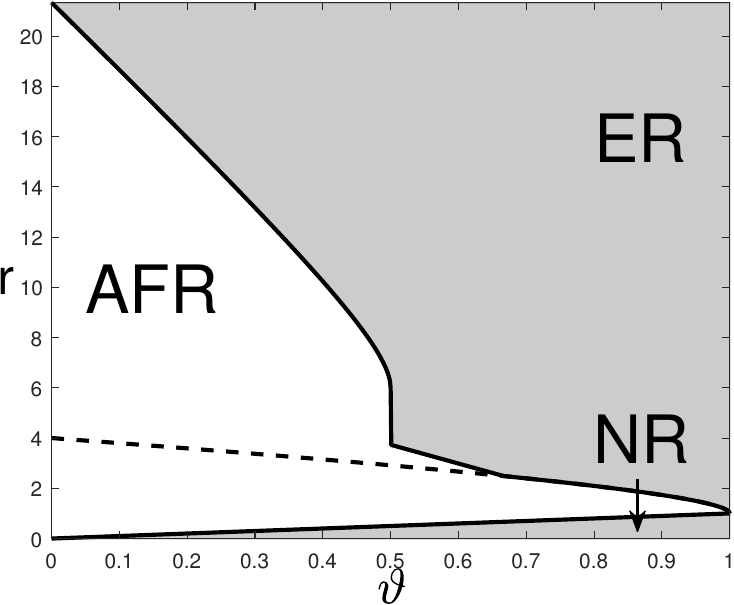}
\caption{The phase diagrams of Lasso-path (block-wise diagonal designs). Left: $\rho=0.5$. Middle: $\rho=-0.5$. Right: zoom-out of the middle panel. 
In all three panels, the dashed lines are the phase curves for orthogonal designs ($\rho=0$), as a reference.
}\label{fig:FDR-blockwise}
\end{figure}

Using Theorem~\ref{thm:lasso}, we can deduce the FDR-TPR tradeoff diagram and phase diagram. To save space, we only present the phase diagram:

\begin{cor}[Phase diagram of Lasso-path]  \label{cor:blockPD-Lasso}
In the setting of Theorem~\ref{thm:lasso}, the phase diagram of Lasso-path is given by 
\[
h_{\mathrm{AFR}}(\vartheta)=\vartheta,\qquad h_{\mathrm{ER}}(\vartheta)=\begin{cases}\max\{h_1(\vartheta), h_2(\vartheta)\},&\mbox{when}\;\;\rho\geq 0,\cr
\max\{h_1(\vartheta), h_2(\vartheta), h_3(\vartheta), h_4(\vartheta)\}, &\mbox{when}\;\; \rho<0, 
\end{cases}
\]
where $h_1(\vartheta)=(1+\sqrt{1-\vartheta})^2$, $h_2(\vartheta) = \bigl(1+\sqrt{\frac{1+|\rho|}{1-|\rho|}}\bigr)^2(1-\vartheta)$, 
$h_3(\vartheta) =\frac{1}{(1+\rho)^2}\bigl(\sqrt{\frac{1+\rho}{1-\rho}}\sqrt{1-2\vartheta}+\sqrt{\frac{1-\rho}{1+\rho}}\sqrt{1-\vartheta}\bigr)^2\cdot 1\{\vartheta<1/2\}$, and $h_4(\vartheta)=\frac{1}{(1+\rho)^2}\bigl(1+\sqrt{\frac{1+\rho}{1-\rho}}\cdot \sqrt{1-2\vartheta}\bigr)^2 \cdot 1\{\vartheta<1/2\}$.
 \end{cor}
\noindent
A visualization of the phase diagram for $\rho=\pm 0.5$ is in Figure~\ref{fig:FDR-blockwise}.

When $\rho=0$, the blockwise diagonal design reduces to the orthogonal design, Lasso-path reduces to \eqref{BH}, and the phase diagram reduces to the one in Proposition~\ref{prop:Identity}. Comparing Figure~\ref{fig:FDR-blockwise} with Figure~\ref{fig:HT}, the phase diagram of Lasso-path is inferior to the one for orthogonal designs. This suggests that the strength of design correlations can have a significant impact on the performance of variable selection. 

Another observation from Figure~\ref{fig:FDR-blockwise} is that the sign of $\rho$ plays a crucial role. This is related to the ``signal cancellation" phenomenon \citep{ke2014covariance,ke2017covariate}. Suppose $\{j, j+1\}$ is a block and both $\beta_j$ and $\beta_{j+1}$ are signals. It is seen that $\mathbb{E}[x_j'y|\beta]=(1+\rho)\tau_p$, whose absolute value is strictly smaller than $\tau_p$ for a negative $\rho$. Hence, when $\rho$ is negative, the signal at $j+1$ creates a ``cancellation effect" and makes $x_j$ marginally less correlated with $y$.  
Lasso is known to be quite vulnerable  to ``signal cancellation" \citep{zhao2006model}. This is why the phase diagram becomes worse when the sign of $\rho$ is flipped. We will provide a more rigorous explanation in Section~\ref{sec:insight} using geometric properties of the Lasso solution.

\subsection{SDP-knockoff} \label{subsec:FDR-lasso}
We now study the SDP-knockoff, where $\mathrm{diag}(s)$ is as in \eqref{knockoff5}. For the block-wise diagonal design parameterized by $\rho$, we have an explicit form of $\mathrm{diag}(s)$:
\beq \label{diag-s}
\mathrm{diag}(s)=(1-a)I_p, \qquad \mbox{where}\quad a = 
\begin{cases}
2|\rho|-1, & |\rho|\geq 1/2,\\
0, & |\rho|<1/2. 
\end{cases}
\eeq 
The value of $a$ controls the correlation between $x_j$ and $\tilde{x}_j$. 
SDP-knockoff aims to minimize this correlation, subject to the eligibility constraints. We first study the case $|\rho|\geq 1/2$.  



\begin{thm}[The case of $|\rho|\geq 1/2$] \label{thm:knockoff-block}
Consider a linear model where \eqref{Gram}-\eqref{RWmodel2} hold. Suppose $n\geq 2p$ and $G$ is as in \eqref{block}, where $|\rho|\geq 1/2$. We construct $\tilde{X}$ in knockoff with $\mathrm{diag}(s)$ as in \eqref{diag-s}. Let $Z_j$, $\tilde{Z}_j$ and $W_j$ be as in \eqref{knockoff2}-\eqref{Wj-knockoff}, where $f$ is the signed maximum in \eqref{knockoff4}. 
For any constant $u>0$, let $\FP_p(u)$ and $\FN_p(u)$ be the expected numbers of false positives and false negatives, by selecting variables with $W_j>\sqrt{2u\log(p)}$. As $p\to\infty$,
\[
\FP_p(u) = L_p p^{1-\min\left\{u,\;\, \vartheta+(\sqrt{u}-|\rho|\sqrt{r})^2+(\xi_\rho\sqrt{r}-\eta_\rho \sqrt{u})_+^2-(\sqrt{r}-\sqrt{u})_+^2\right\}},
\]
and for $\rho\geq 1/2$, 
\[
\FN_p(u) = L_p p^{1- \vartheta-\left\{ (\sqrt{r}-\sqrt{u})_+ - [(1-\xi_\rho)\sqrt{r}-(1-\eta_\rho)\sqrt{u}]_+ - (\lambda_\rho\sqrt{r}-\eta_\rho \sqrt{u})_+\right\}^2 }, 
\]
and for $\rho\leq -1/2$, 
\[
\FN_p(u) = L_p p^{1-\min\bigl\{\vartheta+\left\{ (\sqrt{r}-\sqrt{u})_+ - [(1-\xi_\rho)\sqrt{r}-(1-\eta_\rho)\sqrt{u}]_+ - (\lambda_\rho\sqrt{r}-\eta_\rho \sqrt{u})_+\right\}^2,\;\;\; 2\vartheta\bigr\}}, 
\]
where $\xi_\rho=\sqrt{1-\rho^2}$, $\eta_\rho=\sqrt{(1-|\rho|)/(1+|\rho|)}$, and $\lambda_\rho=\sqrt{1-\rho^2}-\sqrt{1-|\rho|}$. 
\end{thm}

When $|\rho|<1/2$, listing the separate forms of $\FP_p(u)$ and $\FN_p(u)$  is tedious. We instead present the form of $\FP_p(u)+\FN_p(u)$, which is sufficient for deriving the phase diagram.

\begin{thm} [The case of $|\rho| < 1/2$] \label{thm:knockoff-block2}
Consider a linear model where \eqref{Gram}-\eqref{RWmodel2} hold. Suppose $n\geq 2p$ and $G$ is as in \eqref{block}, where $|\rho|< 1/2$. We construct $\tilde{X}$ in knockoff with $\mathrm{diag}(s)$ as in \eqref{diag-s}. Let $Z_j$, $\tilde{Z}_j$ and $W_j$ be as in \eqref{knockoff2}-\eqref{Wj-knockoff}, where $f$ is the signed maximum in \eqref{knockoff4}. 
For any constant $u>0$, let $\FP_p(u)$ and $\FN_p(u)$ be the expected numbers of false positives and false negatives, by selecting variables with $W_j>\sqrt{2u\log(p)}$. As $p\to\infty$,
\begin{align*}
 \FP_p& (u) +  \FN_p(u) = \cr
& \begin{cases}
L_p p^{1-f^+_{\text{Hamm}}(u,r,\vartheta)}, & 0\leq \rho<1/2,\\
L_p p^{1-\min\bigl\{f^+_{\text{Hamm}}(u,r,\vartheta),\;\; 2\vartheta+(\xi_\rho \sqrt{r}-\eta_\rho^{-1} \sqrt{u})_+^2,\;\; 2\vartheta+\frac{(1+2\rho)^2(1-\rho)}{2(1+\rho)}r\bigr\}},  &-1/2<\rho<0, 
\end{cases}
\end{align*}
where 
\begin{align*}
    f^+_{\text{Hamm}}(u,r,\vartheta)& =\min\bigl\{u,\;\; \vartheta+(\sqrt{u}-|\rho|\sqrt{r})^2+((\xi_\rho \sqrt{r}-\eta_\rho \sqrt{u})_+)^2-((\sqrt{r}-\sqrt{u})_+)^2,\\
    &\qquad \vartheta+[(\sqrt{r}-\sqrt{u})_+-((1-\xi_\rho)\sqrt{r}-(1-\eta_\rho)\sqrt{u})_+]^2\bigr\},
\end{align*}
and $\xi_\rho$, $\eta_\rho$ are the same as those in Theorem~\ref{thm:knockoff-block}.
\end{thm}

\begin{figure}[tb]
\centering
\includegraphics[height=.25\textwidth]{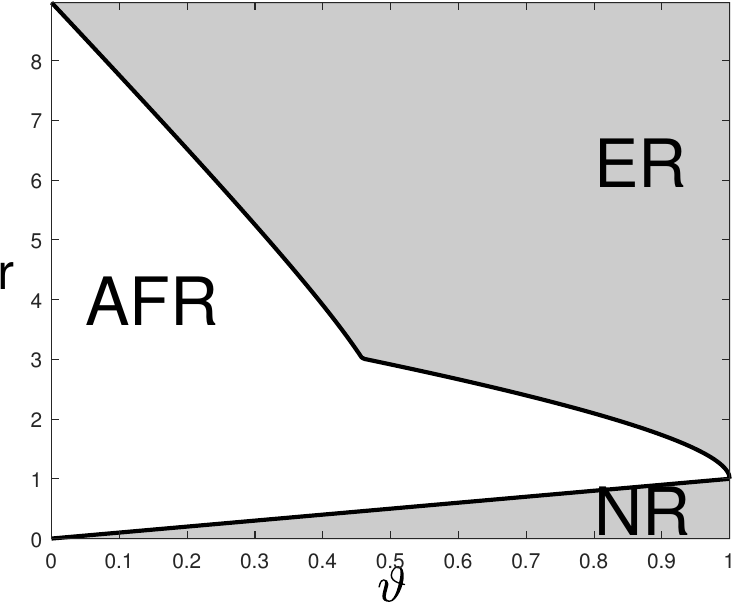}
\includegraphics[height=.25\textwidth]{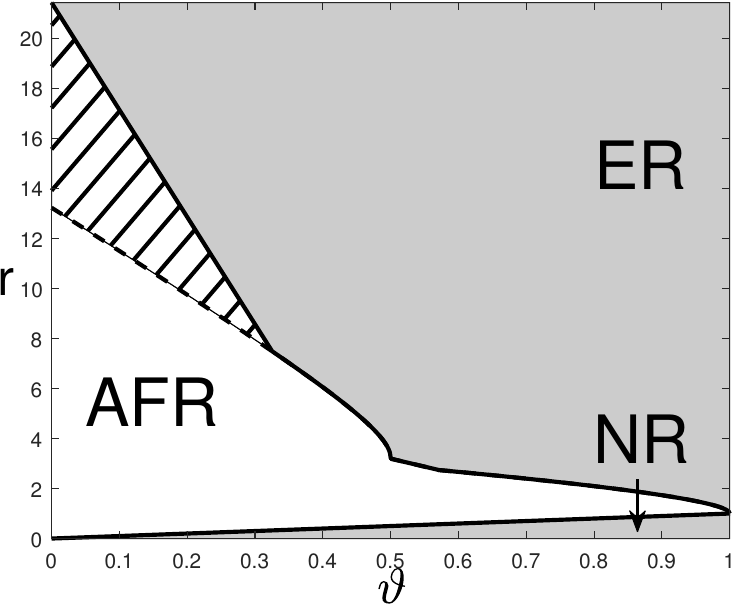}
\includegraphics[height=.25\textwidth]{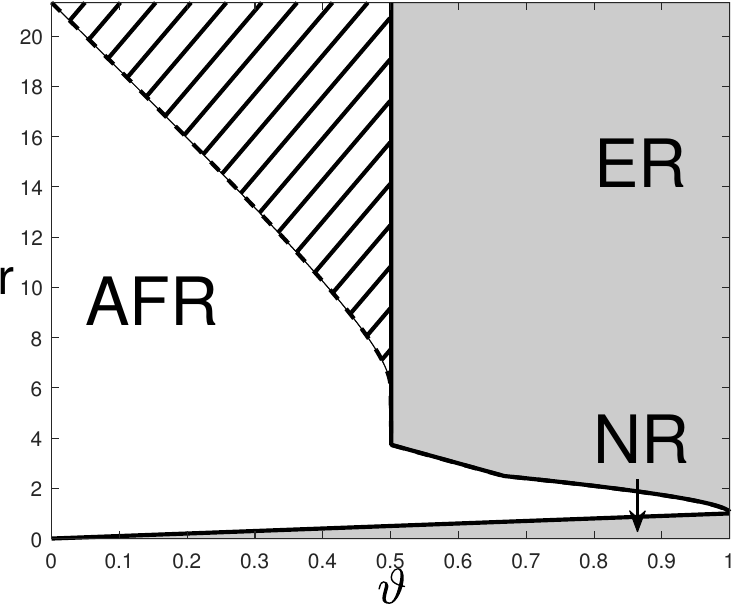}
\caption{The phase diagrams of SDP-knockoff (blockwise diagonal designs; ranking algorithm is Lasso, and symmetric statistic is signed maximum). From left to right, the correlation parameter in the design is $\rho=-0.3$, $\rho=-0.4$, and $\rho=-0.5$, respectively. They correspond to the three cases in Corollary~\ref{cor:blockPD-knockoff}. The shadowed area is the Almost Full Recovery region for SDP-knockoff but Exact Recovery region for the prototype Lasso-path. 
If SDP-knockoff is replaced by CI-knockoff, then in each of three cases the phase diagram is the same as that of Lasso-path. 
}\label{fig:knockoff-phase}
\end{figure}

\begin{cor}[Phase diagram of SDP-knockoff]  \label{cor:blockPD-knockoff}
Suppose conditions of Theorems~\ref{thm:knockoff-block}-\ref{thm:knockoff-block2} hold, where $\rho$ can be any value in $(-1,1)$. Let $h^{\mathrm{lasso}}_{AFR}(\vartheta)$ and $h^{\mathrm{lasso}}_{ER}(\vartheta)$ be the phase curves in Corollary~\ref{cor:blockPD-Lasso}. Define
\[
\rho_0 = \sqrt{2}-1-\sqrt{2-\sqrt{2}} \qquad \mbox{(note: $\rho_0\approx -0.35$)}.
\] 
For SDP-knockoff,  $h_{AFR}(\vartheta)=h^{\mathrm{lasso}}_{AFR}(\vartheta)$, and $h_{ER}(\vartheta)$ has three cases: 
\begin{itemize} \itemsep 0pt
\item When $\rho\in [\rho_0, 1)$, $
h_{ER}(\vartheta)=h^{\mathrm{lasso}}_{ER}(\vartheta)$. 
\item When $\rho\in (-0.5, \rho_0)$, $h_{ER}(\vartheta)=\max\{ h^{\mathrm{lasso}}_{ER}(\vartheta), \, h_5(\vartheta)\}$, where $h_5(\vartheta) =\frac{2(1-2\vartheta)(1+\rho)}{(1+2\rho)^2(1-\rho)}$. 
\item When $\rho\in (-1, -0.5]$, $h_{ER}(\vartheta)=h^{\mathrm{lasso}}_{ER}(\vartheta)$ if $\vartheta>1/2$, and $h_{ER}(\vartheta)=\infty$ otherwise. 
\end{itemize}
\end{cor}
A visualization of the phase diagram for three values of $\rho$ is in Figure~\ref{fig:knockoff-phase}.

Comparing Corollary~\ref{cor:blockPD-knockoff} and Corollary~\ref{cor:blockPD-Lasso}, We have the following observations:
\begin{itemize}
\item When $\rho\in [\rho_0,1)$, SDP-knockoff shares the same phase diagram as Lasso-path, i.e., SDP-knockoff yields a negligible power loss compared with its prototype. 
\item When $\rho\in (-1,\rho_0)$, the phase diagrams of SDP-knockoff is inferior to that of Lasso-path. Especially, when $\rho\in (-1, -0.5]$, the AFR region of SDP-knockoff is infinite: For any $\vartheta<1/2$, no matter how large $r$ is, SDP-knockoff never achieves Exact Recovery.  
\end{itemize}

We give an explanation of the discrepancy of the phase diagram between SDP-knockoff and Lasso-path for $\rho\in (-1, \rho_0)$. First, consider $\rho\in (-0.5, \rho_0)$. By \eqref{diag-s}, $a=0$ and $\mathrm{diag}(s)=I_p$. It follows that the $j$th knockoff is uncorrelated with the $j$th original variable. However, this knockoff is still highly correlated with the $(j+1)$th original variable. Suppose $j$ is a true signal variable. Then, a true signal at $(j+1)$ will increase the absolute correlation between $y$ and $\tilde{x}_j$ but decrease the absolute correlation between $y$ and $x_j$ (since $\rho<0$), making it more difficult for $x_j$ to stand out.  Next, consider $\rho\in (-1,-0.5]$. Suppose $\{j, j+1\}$ has two `nested' signals, i.e., $(\beta_j, \beta_{j+1})=(\tau_p, \tau_p)$. By \eqref{diag-s} and an elementary calculation, 
\[
\mathbb{E}[x_j'y|\beta] = (1+\rho)\tau_p, \qquad \mathbb{E}[\tilde{x}_j'y|\beta] = 
\begin{cases}
\rho \tau_p, & \mbox{when }-0.5<\rho<0, \\
-(1+\rho)\tau_p, & \mbox{when }-1<\rho\leq -0.5. 
\end{cases}
\]
When $\rho\leq -0.5$, variable $j$ and its knockoff have the same absolute correlation with $y$. Consequently, there is a non-diminishing probability that the true signal variable fails to dominate its knockoff variable, making it impossible to select $j$ consistently. In the Rare/Weak signal model, `nested' signals appear with a non-diminishing probability if $\vartheta<1/2$. This explains why $h_{ER}(\vartheta)=\infty$ when $\rho\leq -0.5$ and $\vartheta<1/2$. 

The rationale of SDP-knockoff is to minimize the correlation between a variable and its own knockoff, but this is not necessarily the best strategy for constructing knockoff variables when the original variables are highly correlated. In the next subsection, we will see that a proper increase of the correlation between a variable and its knockoff can boost power.

\subsection{CI-knockoff} \label{subsec:CI-knockoff}
We study the CI-knockoff, where $\mathrm{diag}(s)$ is as in \eqref{CI-knockoff}. 
 \cite{liu2019power} showed that when $c=1$ in \eqref{CI-knockoff}, the resulting $\tilde{X}$ satisfies $x_j'(I-P_{-j})\tilde{x}_j=0$, where $P_{-j}$ is the projection matrix to the linear span of $\{x_k: k\neq j\}$. It means $x_j$ and $\tilde{x}_j$ are {\it conditionally uncorrelated}, conditioning on the other $(p-1)$ original variables. For the block-wise diagonal design \eqref{block}, $\mathrm{diag}(s)$ has an explicit form:  
\beq \label{CI-knockoff2}
\mathrm{diag}(s) = (1-\rho^2)I_p, \qquad \mbox{for all }\rho\in (-1,1). 
\eeq 
Compared with \eqref{cor:blockPD-Lasso}, the value of $a$ has changed. We recall that $a$ controls the correlation between an original variable and its knockoff. In SDP-knockoff, $a$ is chosen as the minimum eligible value, but in CI-knockoff, $a$ is set at $\rho^2$. 

\begin{thm}\label{thm:knockoff-CI}
Consider a linear model where \eqref{Gram}-\eqref{RWmodel2} hold. Suppose $n\geq 2p$ and $G$ is as in \eqref{block}, with a correlation parameter $\rho\in (-1,1)$. We construct $\tilde{X}$ in knockoff with $\mathrm{diag}(s)$ as in \eqref{CI-knockoff2}. Let $Z_j$, $\tilde{Z}_j$ and $W_j$ be as in \eqref{knockoff2}-\eqref{Wj-knockoff}, where $f$ is the signed maximum in \eqref{knockoff4}. 
For any constant $u>0$, let $\FP_p(u)$ and $\FN_p(u)$ be the expected numbers of false positives and false negatives, by selecting variables with $W_j>\sqrt{2u\log(p)}$. As $p\to\infty$,
\begin{align*}
 \FP_p (u) +  \FN_p(u) = 
 \begin{cases}
L_p p^{1-f^+_{\text{Hamm}}(u,r,\vartheta)}, & \rho\geq 0,\\
L_p p^{1-\min\bigl\{f^+_{\text{Hamm}}(u,r,\vartheta),\;\; 2\vartheta+(\xi_\rho \sqrt{r}-\eta_\rho^{-1} \sqrt{u})_+^2\bigr\}},  &\rho<0, 
\end{cases}
\end{align*}
where $f^+_{\text{Hamm}}(u,r,\vartheta)$ is the same as that in Theorem~\ref{thm:knockoff-block2}. 
\end{thm} 

The exponent in Theorem~\ref{thm:knockoff-CI} is in fact the same as that in Theorem~\ref{thm:lasso}. We immediately conclude that CI-knockoff yields the same phase diagram as its prototype, Lasso-path.

\begin{cor}[Phase diagram of CI-knockoff] \label{cor:blockPD-knockoffCI} 
In the setting of Theorem~\ref{thm:knockoff-CI}, for any $\rho\in (-1,1)$, the phase curves of CI-knockoff are the same as those in Corollary~\ref{cor:blockPD-Lasso}. 
\end{cor}

The result of CI-knockoff is very encouraging. We now explain how CI-knockoff improves SDP-knockoff for $\rho\in (-1, \rho_0)$. Comparing \eqref{CI-knockoff2} with \eqref{diag-s}, we find that the correlation between $x_j$ and $\tilde{x}_j$ increases from $\max\{0, 2|\rho|-1\}$ to $\rho^2$. We revisit the scenario of two `nested' signals, i.e., $(\beta_j, \beta_{j+1})=(\tau_p, \tau_p)$. By direct calculations, 
\[
\mathbb{E}[x_j'y|\beta] = (1+\rho)\tau_p, \qquad \mathbb{E}[\tilde{x}_j'y|\beta] = \rho(1+\rho)\tau_p.
\] 
It always holds that $|\mathbb{E}[x_j'y|\beta]|>|\mathbb{E}[\tilde{x}_j'y|\beta]|$. As long as $r$ is sufficiently large, the original variable $x_j$ can standard out. This resolves the previous issue of SDP-knockoff.


Going beyond the block-wise design, it is an interesting question whether CI-knockoff still improves SDP-knockoff. We study it numerically in Section~\ref{sec:simu}, where we consider designs such as {\it Factor models}, {\it Exponential decay}, and {\it Normalized Wishart}; see Experiment 4. 

\smallskip

{\bf Remark 5}. Our theory is focused on the $2\times 2$ blockwise design in \eqref{block}. Using similar techniques, we can study other blockwise designs, such as $k\times k$ blocks or varying-size blocks. Take $k\times k$ blocks for example. In knockoff, solving Lasso in \eqref{knockoff2} reduces to solving many $2k$-dimensional problems separately. Let $J=\{j, j$+$1, \ldots,j$+$k$-$1\}$ be a block. 
The sufficient statistic for $W_j$ is $\tilde{y}=[X_J,\tilde{X}_J]'y\in\mathbb{R}^{2k}$. 
By Lemma~\ref{lem:tool}, the Hamming error at $j$ depends on the interplay between probability contour of $\tilde{y} $ and geometry of the $2k$-dimensional Lasso problem. We make such analysis for $k=2$ in the proofs of Theorems~\ref{thm:knockoff-block}-\ref{thm:knockoff-CI}, which can be extended to a general $k$.

\section{Impact of the ranking algorithm} \label{sec:RankingAlg}
We consider two options of the ranking algorithm, Lasso and least-squares. As the ranking algorithm changes, the prototype is different. In Section~\ref{subsec:OLS}, we first compare two prototypes. In Section~\ref{subsec:knockoff-OLS}, we further compare the associated versions of knockoff. 

In the orthodox knockoff, ranking algorithm is Lasso, augmented design is SDP-knockoff, and symmetric statistic is signed maximum. We re-name it {\it SDP-knockoff-Lasso}. If ranking algorithm is changed to least-squares (with the other two components unchanged), we call it {\it SDP-knockoff-OLS}. 
In each method, if augmented design is changed to CI-knockoff (with the other two components unchanged), we call them {\it CI-knockoff-Lasso} and {\it CI-knockoff-OLS}, respectively. 
SDP-knockoff-Lasso, CI-knockoff-Lasso, and their prototype, Lasso-path, have been studied in Section~\ref{sec:tamperDesign}. 
In this section, we study SDP-knockoff-OLS, CI-knockoff-OLS, and their prototype, least-squares, and compare the results with those in Section~\ref{sec:tamperDesign}. 

We consider the general design, where $G=X'X$ can be any positive definite matrix. We then restrict ourselves to the special case of $2\times 2$ blockwise design in \eqref{block}. The reason we can study general designs is that the least-squares solution has a simple and explicit form (but the Lasso solution does not). 


\subsection{The prototype, least-squares} \label{subsec:OLS}
Before studying SDP-knockoff-OLS and CI-knockoff-OLS, we first study their common prototype, the least-squares (see \eqref{ols}). 

\begin{thm} \label{thm:OLS}
Consider a linear regression model where \eqref{Gram}-\eqref{RWmodel2} hold and $n\geq 2p$. Let $\omega_j$ be the $j$-th diagonal element of $G^{-1}$. Suppose $\min_{1\leq j\leq p}\{\omega_j\}\leq C_0$, for a constant $C_0>0$.  
 Let $W_j^*$ be as in \eqref{ols}. 
For any constant $u>0$, let $\FP_p(u)$ and $\FN_p(u)$ be the expected numbers of false positives and false negatives, by selecting variables with $W^*_j>\sqrt{2u\log(p)}$. As $p\to\infty$, 
 As $p\to\infty$, 
\[
\FP_p(u) \leq L_p \sum_{j=1}^p p^{-\omega_j^{-1}u}, \qquad \FN_p(u) \leq L_p p^{-\vartheta}\sum_{j=1}^p p^{-\omega_j^{-1}(\sqrt{r}-\sqrt{u})_+^2}. 
\]
\end{thm}

\begin{figure}[tb]
\centering
\includegraphics[height=.2\textwidth]{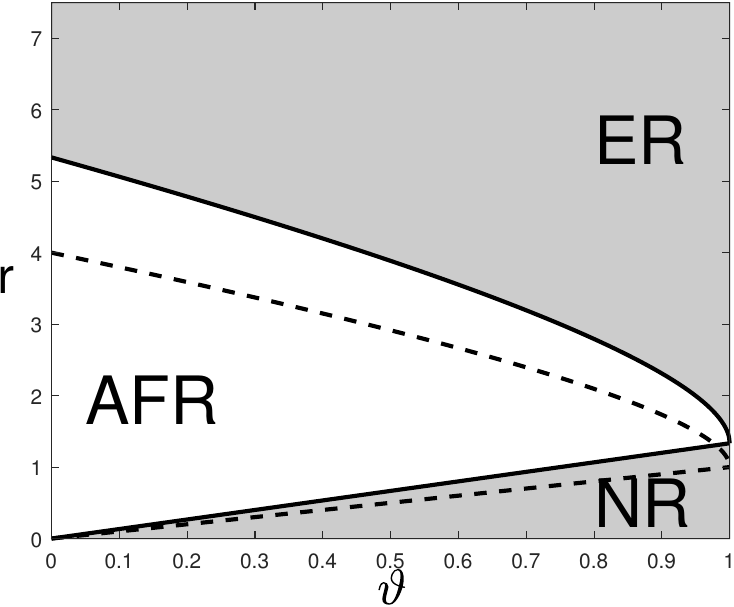}
\includegraphics[height=.2\textwidth]{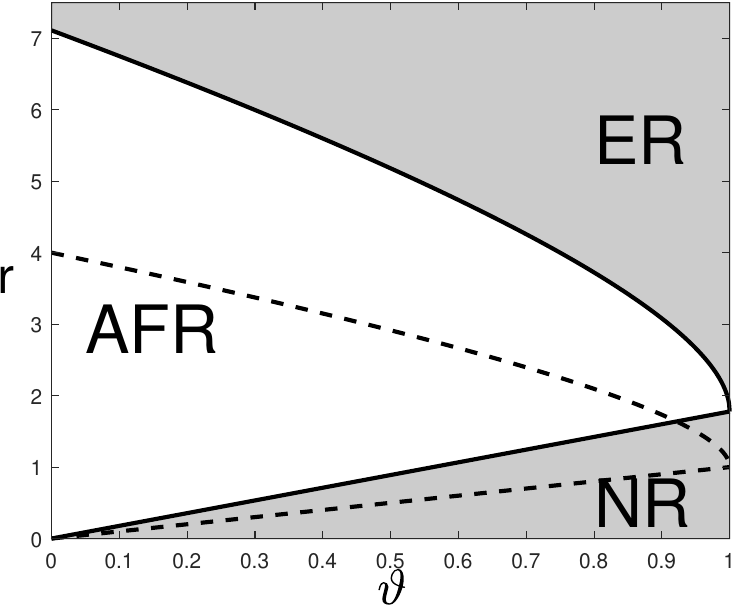}
\includegraphics[height=.2\textwidth]{Figure/phase_knockoff_update}
\includegraphics[height=.2\textwidth]{Figure/phase-05_zoomin}
\caption{The phase diagrams of least-squares (left) and CI-knockoff-OLS (middle left), for blockwise designs with $|\rho|=0.5$. For reference, we also plot the phase diagrams of Lasso-path for $\rho=0.5$ (middle right) and $\rho=-0.5$ (right), which are also the phase diagrams of CI-knockoff-Lasso.  
}\label{fig:OLSwithFDR-phase}
\end{figure}

\begin{cor}[Phase diagram of OLS] \label{cor:blockPD-OLS}
In the setting of Theorem~\ref{thm:OLS}, suppose $G$ is as in \eqref{block}, with a correlation parameter $\rho\in (-1,1)$. Then, $\omega_j=(1-\rho^2)^{-1}$ for $1\leq j\leq p$. As $p\to\infty$, $
\FP_p(u) = L_p p^{1-(1-\rho^2)u}$, and $\FN_p(u) = L_p p^{1-\vartheta - (1-\rho^2)(\sqrt{r}-\sqrt{u})_+^2}$. 
The phase diagram of least-squares is given by $h_{\mathrm{AFR}}(\vartheta) = \frac{\vartheta}{1-\rho^2}$ and $h_{\mathrm{ER}}(\vartheta) = \frac{(1+\sqrt{1-\vartheta})^2}{1-\rho^2}$. 
\end{cor}
\noindent 
Figure~\ref{fig:OLSwithFDR-phase} (left panel) shows the phase diagram of least-squares for $|\rho|= 0.5$; as a reference, in the right two panels, we plot again the phase diagrams of Lasso-path for $\rho=\pm 0.5$. For the comparison between least-squares and Lasso-path, we have the following observations: 
\begin{itemize} \itemsep 0pt
\item In terms of $h_{\mathrm{AFR}}(\vartheta)$, Lasso-path is always better than least-squares. To attain Almost Full Recovery, Lasso-path requires $r>\vartheta$, but least-squares requires $r>\vartheta/(1-\rho^2)$. 
\item In terms of $h_{\mathrm{ER}}(\vartheta)$, Lasso-path is better than least-squares when $\vartheta$ is relatively large (i.e., $\beta$ is comparably sparser), and least-squares is better than Lasso-path when $\vartheta$ is relatively small (i.e., $\beta$ is comparably denser). 
\item The sign of $\rho$ also matters. For small $\vartheta$, the advantage of least-squares over Lasso-path on $h_{\mathrm{ER}}(\vartheta)$ is much more obvious when $\rho$ is negative.   
 \end{itemize}
We give an intuitive explanation to the above phenomena. We say a signal variable (i.e., $\beta_j\neq 0$) is `isolated' if it is the only signal variable in the $2\times 2$ block, and we say two signals are `nested' if they are in the same $2\times 2$ block. In the sparser regime (i.e., $\vartheta$ is large), least-squares has a disadvantage because it is inefficient in discovering an `isolated' signal. In the less sparse regime (i.e., $\vartheta$ is small), Lasso-path has a disadvantage because it suffers from signal  cancellation when estimating a pair of `nested' signals (`signal cancellation' means a signal variable has a weak marginal correlation with $y$ due to the effect of other signals correlated with this one).
A more rigorous explanation is given in Section~\ref{sec:insight}, using geometry of solutions of least-squares and Lasso; see Lemma~\ref{lem:block-RejectRegion}, Figure~\ref{fig:reject-region-block}, and discussions therein. 

\subsection{Knockoff-OLS} \label{subsec:knockoff-OLS}
We now study SDP-knockoff-OLS and CI-knockoff-OLS. The next theorem provides a general result that applies to all augmented designs: 
\begin{thm} \label{thm:knockoffOLS-block}
Consider a linear model where \eqref{Gram}-\eqref{RWmodel2} hold. Suppose $n\geq 2p$. We construct $\tilde{X}$ in knockoff as in \eqref{knockoff1}, with some choice of $\mathrm{diag}(s)$. Write $G^*=[X, \tilde{X}]'[X, \tilde{X}]\in\mathbb{R}^{2p\times 2p}$. Suppose $\mathrm{diag}(s)$ is chosen such that $G^*$ is non-singular. 
Let $A_j\in\mathbb{R}^{2\times 2}$ be the submatrix of $(G^*)^{-1}$ restricted to the $j$th and $(j+p)$th rows and columns. Denote $\omega_{1j}=A_j(1,1)$ and $\omega_{2j}=A_j(1,2)$. Suppose $\max_{1\leq j\leq p}\{\omega_{1j}\}\leq C_0$, for a constanat $C_0>0$.
Let $Z_j$, $\tilde{Z}_j$ and $W_j$ be as in \eqref{knockoff6} and \eqref{Wj-knockoff}, where $f$ is the signed maximum in \eqref{knockoff4}. 
For any constant $u>0$, let $\FP_p(u)$ and $\FN_p(u)$ be the expected numbers of false positives and false negatives, by selecting variables with $W_j>\sqrt{2u\log(p)}$. As $p\to\infty$,  
\[
\FP_p(u) \leq L_p \sum_{j=1}^p p^{-\omega_{1j}^{-1}u}, \qquad \FN_p(u) \leq L_pp^{-\vartheta} \sum_{j=1}^p p^{-\omega_{1j}^{-1}\min\bigl\{(\sqrt{r}-\sqrt{u})_+^2,\;\;\frac{\omega_{1j}}{\omega_{1j}+|\omega_{2j}|}\cdot \frac{1}{2}r\bigr\}}.
\]
\end{thm}

The phase diagram of knockoff-OLS is governed by the quantities $\{\omega_{1j}\}_{1\leq j\leq p}$. We now consider the special case of the $2\times 2$ blockwise design in \eqref{block}, where the augmented design is such that $\mathrm{diag}(s)=(1-a)I_p$, where $a=\max\{0, 2|\rho|-1\}$ in SDP-knockoff and $a=\rho^2$ in CI-knockoff. We note that in SDP-knockoff, the matrix $G^*$ is singular when $|\rho|\geq 1/2$. In other words, SDP-knockoff-OLS is well defined only for $|\rho|<1/2$.

\begin{cor}[Phase diagram of knockoff-OLS]  \label{cor:blockPD-FDR-OLS}
In the same setting of Theorem~\ref{thm:knockoffOLS-block}, suppose $G$ is as in \eqref{block}, with a correlation parameter $\rho\in (-1,1)$.
\begin{itemize}
\item SDP-knockoff-OLS (only defined for $|\rho|<1/2$): $\omega_{1j}=(1-4\rho^2)^{-1}(1-2\rho^2)$ for $1\leq j\leq p$. 
The phase diagram is given by $h_{\mathrm{AFR}}(\vartheta) = \frac{\vartheta(1-2\rho^2)}{1-4\rho^2}$ and $h_{\mathrm{ER}}(\vartheta) = \frac{(1+\sqrt{1-\vartheta})^2(1-2\rho^2)}{1-4\rho^2}$. 
\item CI-knockoff-OLS: $\omega_{1j}=(1-\rho^2)^{-2}$ for $1\leq j\leq p$. 
The phase diagram is given by $h_{\mathrm{AFR}}(\vartheta) = \frac{\vartheta}{(1-\rho^2)^2}$ and $h_{\mathrm{ER}}(\vartheta) = \frac{(1+\sqrt{1-\vartheta})^2}{(1-\rho^2)^2}$. 
\end{itemize}
\end{cor}
Figure~\ref{fig:OLSwithFDR-phase} (second left panel) shows the phase diagram of CI-knockoff-OLS for $|\rho|=0.5$. In this figure, the right two panels are the phase diagrams of CI-knockoff-Lasso for $\rho=\pm 0.5$.

From Corollary~\ref{cor:blockPD-FDR-OLS} and Figure~\ref{fig:OLSwithFDR-phase}, we draw two conclusions: First, for both SDP-knockoff-OLS and and CI-knockoff-OLS, whenever $\rho\neq 0$, their phase diagrams are strictly inferior to the phase diagram of the least-squares (prototype). This is different from the case of using Lasso as ranking algorithm, where the phase diagrams of CI-knockoff-Lasso and Lasso-path (prototype) are the same in the blockwise design for all $\rho\in (-1,1)$. Second, the comparison of CI-knockoff-OLS and CI-knockoff-Lasso is largely similar to the comparison between the least-squares and Lasso-path (see Section~\ref{subsec:OLS}). 

\smallskip

{\bf Remark 6}. When we use the least-squares as the ranking algorithm, such a gap between knockoff and its prototype always exists, for a general design. 
To see this, note that by Theorem~\ref{thm:OLS} and Theorem~\ref{thm:knockoffOLS-block}, the phase diagrams of knockoff-OLS and its prototype are governed by the quantities $\{\omega_{1j}\}_{1\leq j\leq p}$ and $\{\omega_{j}\}_{1\leq j\leq p}$, respectively. 
Since $\omega_{1j}$ and $\omega_j$ are the $j$th diagonal elements of $G^{-1}$ and $(G^*)^{-1}$, respectively, and $G$ is a principal submatrix of $G^*$, it follows by elementary linear algebra that $
\omega_{j}\leq \omega_{1j}$ is always true (and this inequality is often strict). Unfortunately, it is impossible to mitigate this gap by using the augmented design in \eqref{knockoff1}, no matter how we choose $\mathrm{diag}(s)$. \cite{xing2019controlling} proposed a new idea of constructing an augmented design, called the Gaussian mirror, which is tailored to using the least-squares as the ranking algorithm. In a companion paper \citep{Paper2}, we show that the Gaussian mirror attains the same phase diagram as the least-squares.


{\bf Remark 7}. Besides Lasso and least-squares, we may consider other ranking algorithms, such as the thresholded Lasso, non-convex penalization methods, and  the forward-backward selection.
See \cite{ke2021comparison} about the phase diagrams of these methods.

\section{The proof ideas and some geometric insights} \label{sec:insight}

A key technical tool in the proof is the following lemma, which is proved in the Appendix. Recall that $L_p$ is a generic notation of multi-$\log(p)$ terms; see Definition~\ref{Multi-Log(p)}.
For a vector $v$, $\|v\|$ denotes the $\ell^2$-norm; for a matrix $M$, $\|M\|$ denotes the spectral norm. 
\begin{lem}  \label{lem:tool}
Fix an integer $d\geq 1$, a vector $\mu\in \mathbb{R}^d$, a covariance matrix $\Sigma\in\mathbb{R}^{d\times d}$, and an open set $S\subset\mathbb{R}^d$ such that $\mu\notin S$. The quantities $(d, \mu, \Sigma, S)$ do not change with $p$. 
Suppose $
b\equiv \inf_{x\in S}\{(x-\mu)'\Sigma^{-1}(x-\mu)\}<\infty$. 
Consider a sequence of random vectors $X_p\in\mathbb{R}^d$, indexed by $p$, satisfying that 
\[
X_p|(\mu_p, \Sigma_p)\sim {\cal N}_d\Bigl(\mu_p,\;\; \frac{1}{2\log(p)}\Sigma_p\Bigr),
\]
where $\mu_p\in\mathbb{R}^d$ is a random vector and $\Sigma_p\in\mathbb{R}^{d\times d}$ is a random covariance matrix. As $p\to\infty$, suppose for any fixed $\gamma>0$ and $L>0$, $\mathbb{P}(\|\mu_p-\mu\|>\gamma)\leq p^{-L}$ and $\mathbb{P}(\|\Sigma_p-\Sigma\|>\gamma)\leq p^{-L}$. Then, as $p\to\infty$,
\[
\mathbb{P}(X_p\in S)=L_pp^{-b}, 
\]
or equivalently, $p^{b+\delta}\mathbb{P}(X_p\in S)\to\infty$ and $p^{b+\delta}\mathbb{P}(X_p\in S)\to 0$ for any constant $\delta>0$. 
\end{lem}

\noindent
This lemma connects the rate of convergence of $\mathbb{P}(X_p\in S)$ with the geometric property of the set $S$. The exponent $b$ is the ``radius'' of the largest ellipsoid that centers at $\mu$ and is fully contained in the complement of $S$.  

\begin{figure}[tb!]
\includegraphics[width=.245\textwidth]{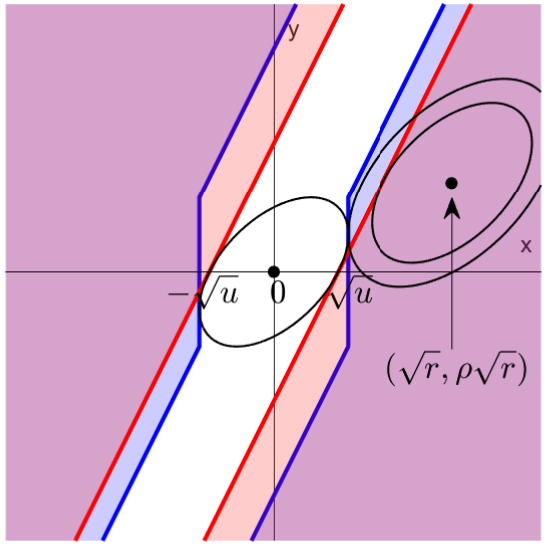} 
\includegraphics[width=.245\textwidth]{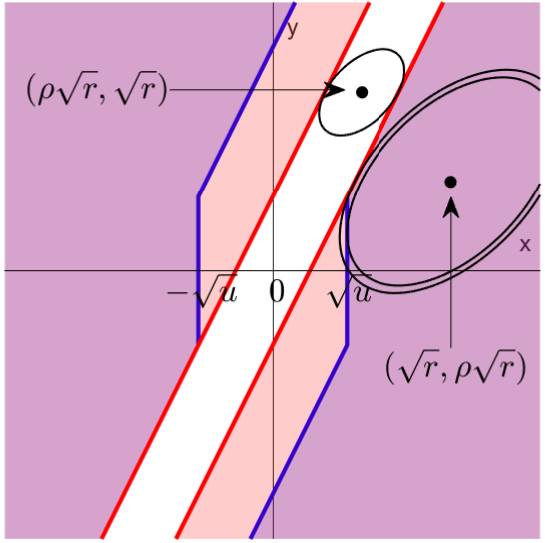}
\includegraphics[width=.245\textwidth]{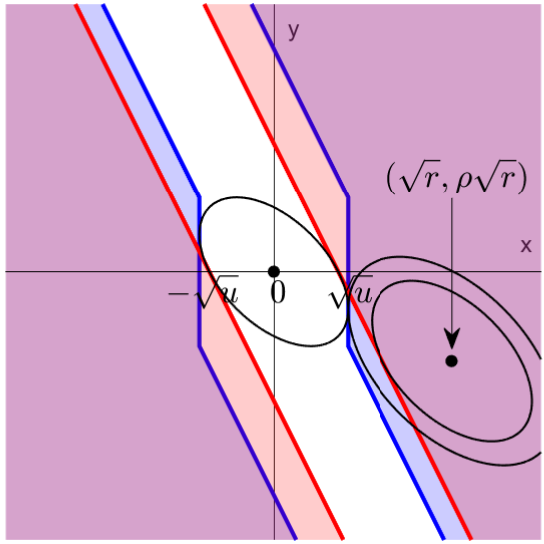} 
\includegraphics[width=.245\textwidth]{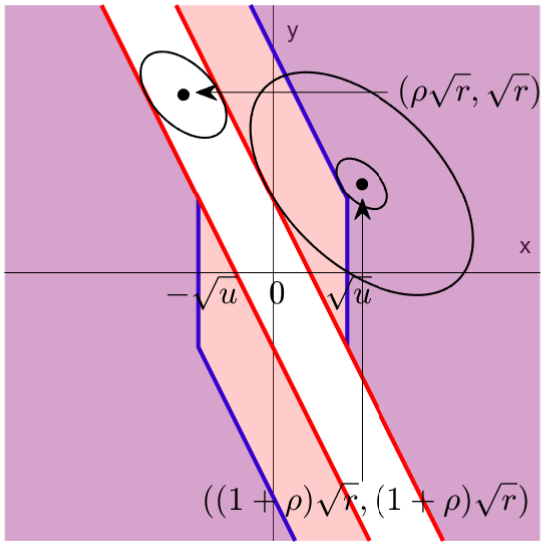}
\caption{Rejection regions and `most-likely' cases in block-wise diagonal designs (x-axis: $x_j'y/\sqrt{2\log(p)}$; y-axis: $x_{j+1}'y/\sqrt{2\log(p)}$). From left to right: (i) positive $\rho$ and large $\vartheta$, (ii) positive $\rho$ and small $\vartheta$, (iii) negative $\rho$ and large $\vartheta$, (iv) negative $\rho$ and small $\vartheta$. In each plot, the blue solid lines define rejection region of Lasso-path, and the red solid lines define rejection region of least-squares. For each method, $\FP_p$ is determined by the largest FP-ellipsoid in ${\cal R}^c$, and $\FN_p$ is determined by the largest FN-ellipsoid in ${\cal R}$, where the centers of these ellipsoids are determined by $(\beta_j, \beta_{j+1})$ in the `most-likely' case. In each plot, the largest FP-ellipsoid is controlled to be the same for both Lasso-path and least-squares, and so the method with a larger FN-ellipsoid is better.} \label{fig:reject-region-block}
\end{figure}

\paragraph{Proof sketch.} 
We illustrate how to use Lemma~\ref{lem:tool} to prove the theorems in Sections~\ref{sec:orthogonal}-\ref{sec:RankingAlg}. Take the proof of Theorem~\ref{thm:lasso} for example. Consider the block-wise design in \eqref{block}. Under this design, the objective of Lasso is separable, and it reduces to solving many 2-dimensional Lasso problems separably. Fix $j$ and suppose $\{j,j+1\}$ is a block. Let $W_j^*$ be as in \eqref{lasso}. Write
\beq \label{h-hat-temp}
\hat{h} = \bigl(x_j'y,\;  x_{j+1}'y\bigr)'/\sqrt{2\log(p)}\;\; \in\;\; \mathbb{R}^2. 
\eeq
Since the Lasso objective is separable,  
$(W_j^*, W_{j+1}^*)$ are purely determined by $\hat{h}$. Particularly, there exists ${\cal R}_u\subset\mathbb{R}^2$, such that $W_j^*>\sqrt{2u\log(p)}$ if and only if  $\hat{h}\in {\cal R}_u$. We call ${\cal R}_u$ the ``rejection region'' of Lasso-path. 
The probabilities of a false positive and a false negative occurring at $j$ are respectively 
\[
\mathbb{P}\bigl(\hat{h}\in {\cal R}_u,\, \beta_j=0\bigr)\qquad \mbox{and}\qquad \mathbb{P}\bigl(\hat{h}\in {\cal R}^c_u,\, \beta_j=\tau_p\bigr). 
\]
Conditioning on $\beta$, the random vector $\hat{h}$ has a bivariate normal distribution, whose mean is a constant vector  and whose covariance matrix is $\frac{1}{2\log(p)} B$, where $B$ is the same as in \eqref{block}. 
Applying Lemma~\ref{lem:tool}, we reduce the proof into two steps: In Step 1, we derive the rejection region ${\cal R}_u$. In Step 2,  for each possible realization of $\beta$ with $\beta_j=0$, we calculate $
b(\beta)\equiv \inf_{x\in {\cal R}_u}\{(x-\mu(\beta))'B^{-1}(x-\mu(\beta))\}$, and for each possible realization of $\beta$ with $\beta_j\neq 0$, we calculate $
b(\beta)\equiv \inf_{x\in {\cal R}^c_u}\{(x-\mu(\beta))'B^{-1}(x-\mu(\beta))\}$, where $\mu(\beta)\equiv \mathbb{E}[\hat{h}|\beta]$. 
Both steps can be carried out by direct calculations. We use a similar strategy to prove other theorems. The proof is sometimes complicated. For example, to analyze knockoff for block-wise diagonal designs, we have to consider the random vector $\hat{h}=(x_j'y, x_{j+1}'y, \tilde{x}_j'y,  \tilde{x}_{j+1}'y)'/\sqrt{2\log(p)}\in\mathbb{R}^4$. The proof requires deriving a 4-dimensional rejection region and calculating $b(\beta)$, for an arbitrary $\rho\in (-1,1)$. The calculations are very tedious. 

\paragraph{The geometric insight about two prototypes.} 
We use the geometric interpretation of our proofs to give more insights about Lasso-path versus least-squares (see Corollary~\ref{cor:blockPD-Lasso} and Corollary~\ref{cor:blockPD-OLS}). 
Under the blockwise design \eqref{block}, for each method, the objective is separable, so that the event $W_j^*>\sqrt{2u\log(p)}$ can be described via a 2-dimensional rejection region. 
The next lemma gives the rejection regions of Lasso-path and least-squares: 
\begin{lem}\label{lem:block-RejectRegion}
Consider a linear model, where the Gram matrix satisfies \eqref{block}, with a correlation parameter $\rho\in (-1,1)$.
Let $W_j^{*,path}$ and $W_j^{*,old}$ be as in \eqref{lasso} and \eqref{ols}, respectively. 
Suppose $\{j, j+1\}$ is a block. 
Write $\hat{h}=(x_j'y, x_{j+1}'y)'/\sqrt{2\log(p)}$.
Define
\begin{align*}
{\cal R}^{\mathrm{path}}_u(\rho) &= \{(h_1,h_2):h_1-\rho h_2 >(1-\rho)\sqrt{u},\, h_1>\sqrt{u}\}\cr
&\qquad \cup \{(h_1,h_2): h_1-\rho h_2 > (1+\rho)\sqrt{u}\}\cr&\qquad \cup  \{(h_1,h_2): h_1- \rho h_2 < -(1-\rho)\sqrt{u},\,  h_1<-\sqrt{u}\}\cr
&\qquad \cup \{(h_1,h_2):x-\rho y<  -(1+\rho)u\},
\quad \mbox{for }\rho\geq 0,\cr
{\cal R}^{\mathrm{path}}_u(\rho) & = \{(h_1, h_2): (h_1, -h_2)\in {\cal R}^{\mathrm{path}}_u(-\rho)\}, \qquad\mbox{for }\rho <0, \cr 
{\cal R}^{\mathrm{ols}}_u(\rho) & = \{(h_1,h_2): h_1-\rho h_2> (1-\rho^2)\sqrt{u}\}\cr
&\qquad \cup \{(h_1,h_2): h_1-\rho h_2<- (1-\rho^2)\sqrt{u}\}.  
\end{align*}
Then, for Lasso-path, $W_j^{*,path}>\sqrt{2u\log(p)}$ if and only if $\hat{h}\in {\cal R}^{\mathrm{path}}_u(\rho)$; for least-squares,  $W_j^{*,ols}>\sqrt{2u\log(p)}$ if and only if $\hat{h}\in {\cal R}^{\mathrm{ols}}_u(\rho)$. 
\end{lem}

These rejection regions are shown in Figure~\ref{fig:reject-region-block}. Their geometric properties are different for positive and negative $\rho$. Fix $j$. Let $\hat{h}$ be as in \eqref{h-hat-temp}, and write $\mu(\beta)=\mathbb{E}[\hat{h}|\beta]$. 
\begin{itemize}
\item The rate of convergence of $\FP_p(u)$ is determined by the largest ellipsoid that centers at $\mu(\beta)$ and is contained in $\mathcal{R}_u^c$. We call this ellipsoid the {\it FP-ellipsoid}. 
\item The rate of convergence of $\FN_p(u)$ is determined by the largest ellipsoid that centers at $\mu(\beta)$ and is contained in $\mathcal{R}_u$. We call this ellipsoid the {\it FN-ellipsoid}. 
\end{itemize}
By direct calculations, $
\mu(\beta)=\bigl(\beta_j+\rho\beta_{j+1},\; \rho\beta_j + \beta_{j+1}\bigr)'/\sqrt{2\log(p)}$. 
Under our model, $(\beta_j, \beta_{j+1})$ has 4 possible values  $\{(0, 0), (0,\tau_p),  (\tau_p,0),(\tau_p, \tau_p)\}$, where the first two correspond to a null at $j$ and the last two correspond to a non-null at $j$. The probability of having a selection error at $j$ thus splits into 4 terms, and which term is dominating depends on the values of $\vartheta$ and $\rho$. 
The realization of $(\beta_j, \beta_{j+1})$ that plays a dominating role is called the `most-likely' case. 
For example, when $\vartheta$ is large (i.e., $\beta$ is sparser), the most-likely case of a false positive occuring at $j$ is when $(\beta_j, \beta_{j+1})=(0,0)$; when $\vartheta$ is small (i.e., $\beta$ is less sparse), the most-likely case of a false positive is when $(\beta_j, \beta_{j+1})=(0,\tau_p)$. 
Table~\ref{tb:cases} summarizes the `most-likely' cases. We also visualize the `most-likely' cases for different $(\rho, \vartheta)$ in Figure~\ref{fig:reject-region-block}. In each plot of Figure~\ref{fig:reject-region-block}, we have coordinated the thresholds $u$ in two methods so that the FP-ellipsoid is exactly the same. It suffices to compare the FN-ellipsoid: The method with a larger FN-ellipsoid has a faster rate of convergence on the Hamming error. 
It is clear that, when $\vartheta$ is large,  the FN-ellipsoid of Lasso-path is larger; when  $\vartheta$ is small, the FN-ellipsoid of least-squares is larger. This explains the different performances of two methods. Moreover, when $\vartheta$ is small, comparing the case of a positive $\rho$ with the case of a negative $\rho$, we find that the difference between FN-ellipsoids of two methods are much more prominent in the case of a negative $\rho$. This explains why the sign of $\rho$ matters.

\begin{table}[hbt]
\centering
\scalebox{.88}{
\begin{tabular}{cccccc}
\toprule
Sparsity & Correlation & Error type &  Most-likely case  & Center of ellipsoid \\
\midrule 
\multirow{2}{*}{large $\vartheta$} & \multirow{2}{*}{positive/negative $\rho$} & FP & $\beta_j=0$, $ \beta_{j+1}=0$ &  $(0,\, 0)$\\
& & FN & $\beta_j=\tau_p$, $\beta_{j+1}=0$   & $(\sqrt{r},\, \rho\sqrt{r})$\\
\hline 
\multirow{2}{*}{small $\vartheta$} &  \multirow{2}{*}{positive $\rho$} & FP & $\beta_j=0$, $\beta_{j+1}=\tau_p$  & $(\rho \sqrt{r}, \, \sqrt{r})$\\
& &FN & $\beta_j=\tau_p$, $\beta_{j+1}=0$  & $(\sqrt{r},\, \rho\sqrt{r})$ \\
\hline 
\multirow{2}{*}{small $\vartheta$} &  \multirow{2}{*}{negative $\rho$}&FP & $\beta_j=0$, $\beta_{j+1}=\tau_p$  & $(\rho \sqrt{r}, \, \sqrt{r})$\\
& & FN & $\beta_j=\tau_p$, $\beta_{j+1}=\tau_p$ & $\bigl((1+\rho)\sqrt{r}, \, (1+\rho)\sqrt{r}\bigr)$ \\
\bottomrule
\end{tabular}}
\caption{The `most-likely' cases and the corresponding ellipsoid center $\mu(\beta)$} \label{tb:cases}
\end{table}

\section{Simulations} \label{sec:simu}

We use numerical experiments to support and exemplify the theoretical results in Sections~\ref{sec:orthogonal}-\ref{sec:RankingAlg}. In Experiments 1 and 2, we consider orthogonal designs and block-wise diagonal designs, respectively. In Experiments 3 and 4, we consider other design classes, including block-wise diagonal designs with larger blocks, factor models, exponentially decaying designs, and normalized Wishart designs. 
We consider four methods, Lasso-path (Lasso), least-squares (OLS), knockoff with Lasso-path ranking (KF.Lasso) and with least-squares ranking (KF.OLS). We use either the signed maximum or the difference as the symmetric statistic, and for KF we choose  $\mathrm{diag}(s)=\min\{1, \, 2\lambda_{\min}(G)\}\cdot I_p$, unless specified otherwise. It is called the equi-correlated knockoff (EC-KF), and is the same as the SDP-knockoff for orthogonal designs and the $2\times 2$ block-wise diagonal designs. In Experiments 1-3, this is the only $\mathrm{diag}(s)$ we use, and so we write EC-KF as KF for short. In Experiments 4, we also consider the conditional independence knockoff (CI-KF). For most experiments, fixing a parameter setting, we generate $200$ data sets and record the averaged Hamming selection error among these 200 repetitions. 

\vspace{0.1in}

\begin{figure}[tb]
\centering
\includegraphics[height=.25\textwidth,  width=.25\textwidth, trim=0 35 0 45, clip=true]{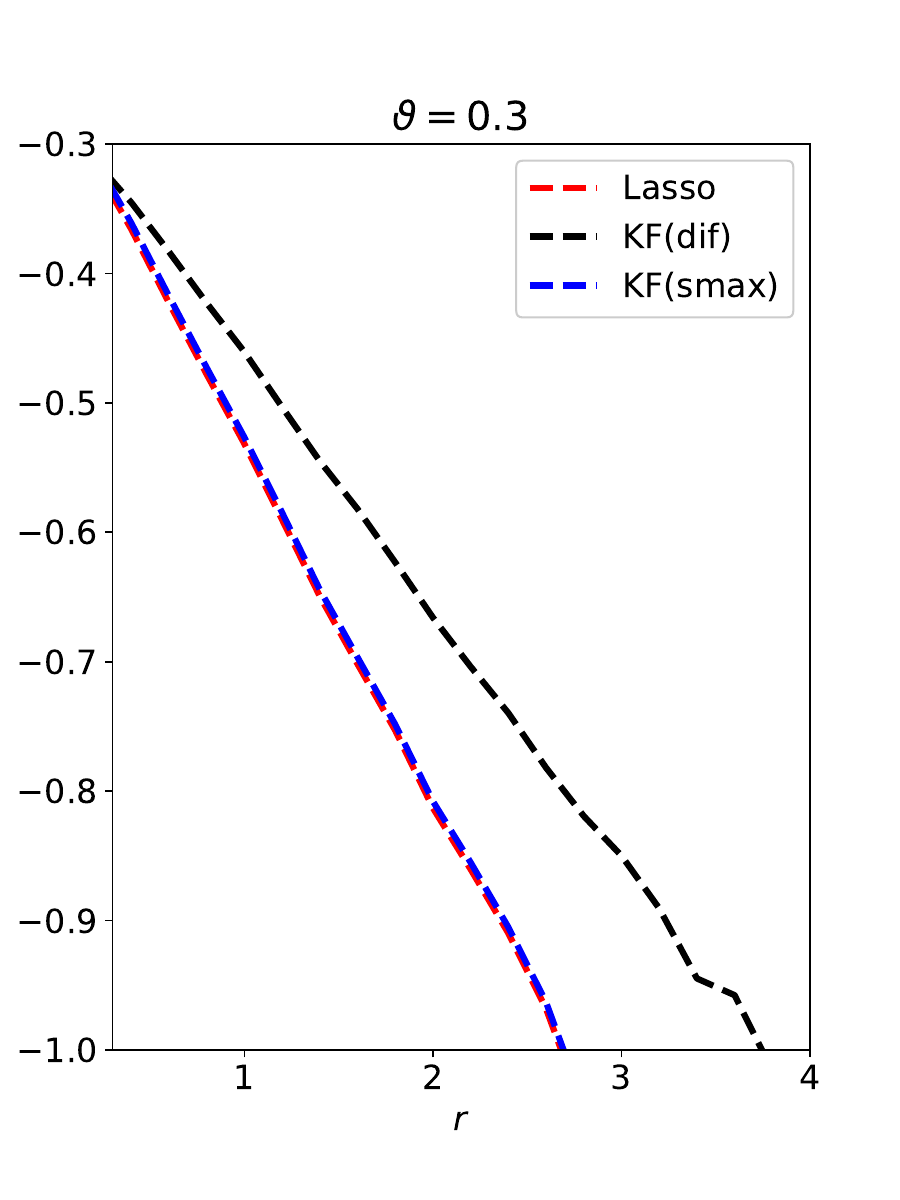} $\qquad$
\includegraphics[height=.25\textwidth,  width=.25\textwidth, trim=0 35 0 45, clip=true]{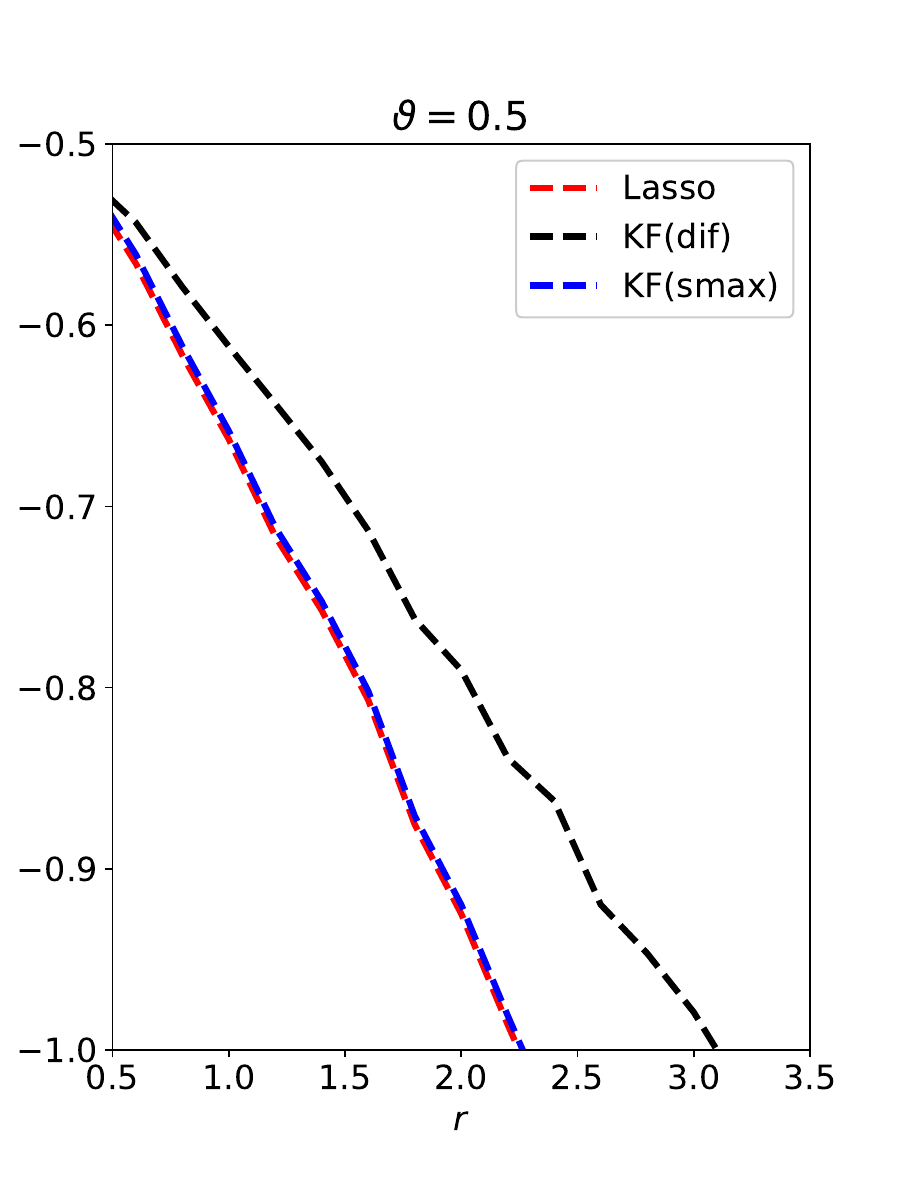} $\qquad$
\caption{Experiment 1 (orthogonal designs). The y-axis is $\log_p(H_p/p)$, where $H_p$ is the average Hamming error over 200 repetitions. 
}\label{fig:exp1}
\end{figure}

{\bf Experiment 1}. We investigate  the performance of different methods for orthogonal designs. Given $(n,p)=(2000,1000)$, $\vartheta\in \{0.3,0.5\}$ and $r$ ranging on a grid from $0$ to $6$ with step size  $0.2$, we generate data $y$ from $N(X\beta, I_n)$ where $X$ is an $n\times p$ matrix with unit length columns that are orthogonal to each other and $\beta$ is generated from (\ref{RWmodel1}). We implemented Lasso and KF.Lasso using both the signed maximum and the difference as the symmetric statistic. Under the orthogonal design, Lasso and OLS yield the same importance metric thus OLS and KF.OLS are neglected in this experiment.
Each method outputs $p$ importance statistics, and we threshold these importance statistics at $\sqrt{2u^*\log(p)}$ where $u^*$ minimizes $\FN_p(u)+\FN_p(u)$ in theory. 
The results are in Figure~\ref{fig:exp1}, where the y-axis is $\log_p(H_p/p)$, and $H_p$ is the averaged Hamming selection error over 200 repetitions.

The theory in Sections~\ref{sec:orthogonal}-\ref{sec:tamperDesign} suggests the following for orthogonal designs: (i) Regarding the choice of symmetric statistic for KF, the signed maximum outperforms the difference. (ii) With signed maximum as the symmetric statistic, KF.Lasso has a similar performance as Lasso. These theoretical results are perfectly validated by simulations (see Figure~\ref{fig:exp1}).

\vspace{0.1in}

\begin{figure}[tb]
\centering
\hspace*{-.8cm}
\includegraphics[height=.25\textwidth, width=.25\textwidth, trim=0 35 0 45, clip=true]{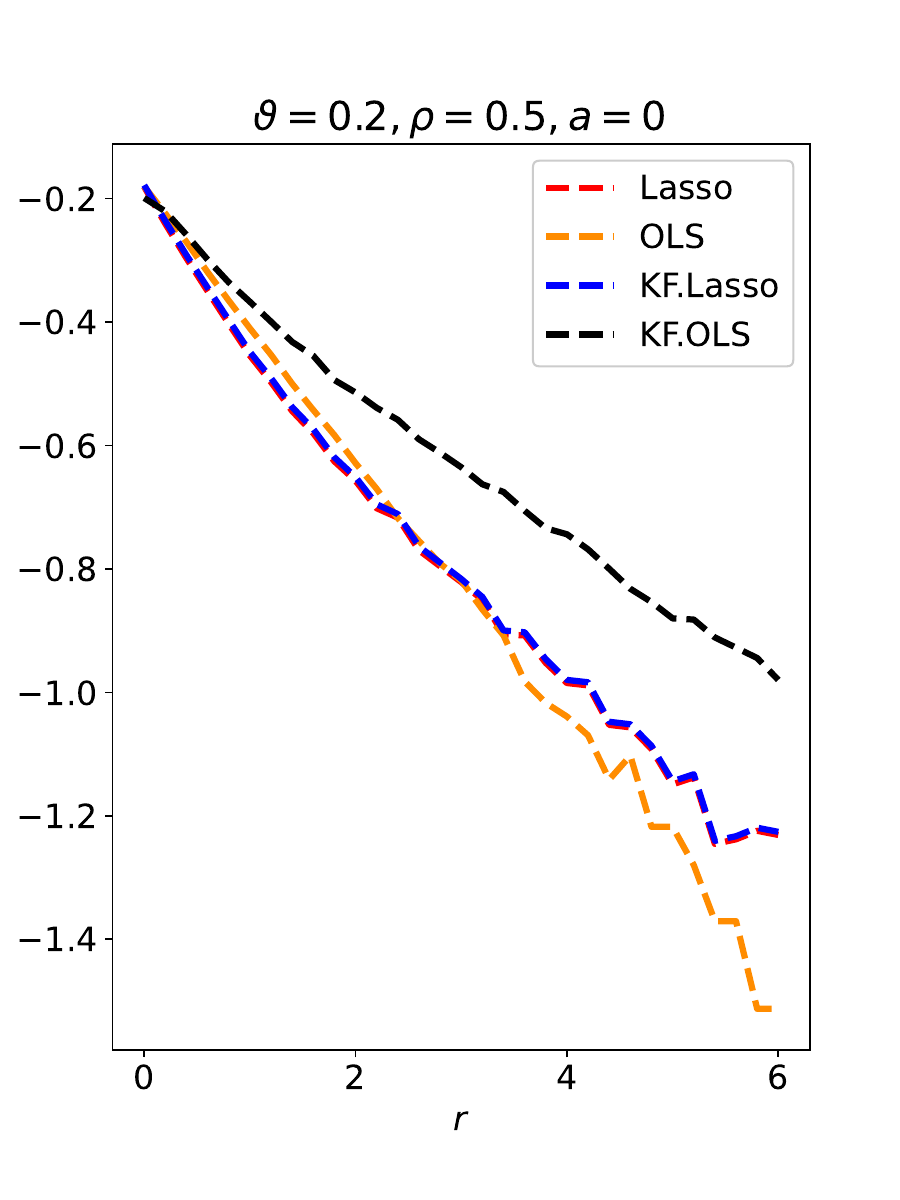} $\qquad$
\hspace*{-1cm}
\includegraphics[height=.25\textwidth,  width=.25\textwidth, trim=0 35 0 45, clip=true]{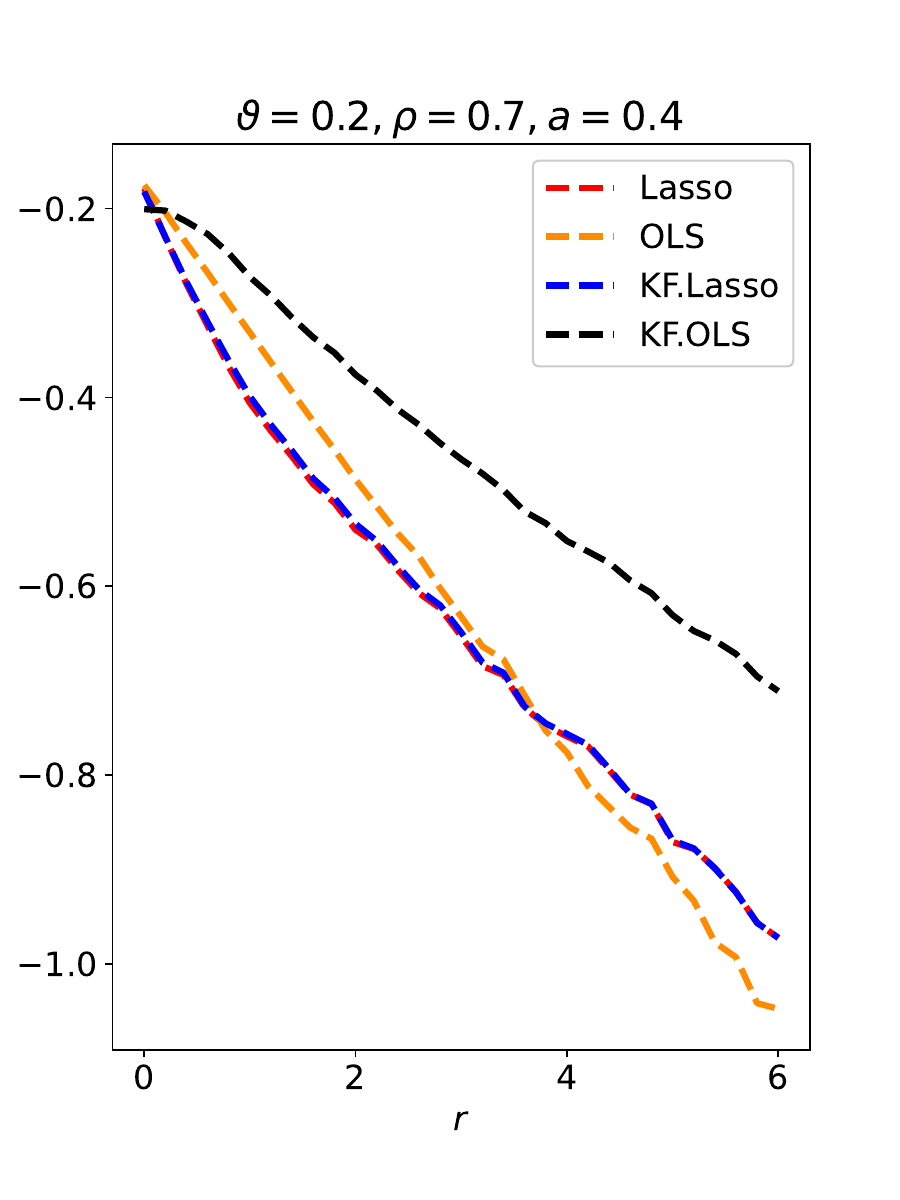} $\qquad$ 
\hspace*{-1cm}
\includegraphics[height=.25\textwidth,  width=.25\textwidth, trim=0 35 0 45, clip=true]{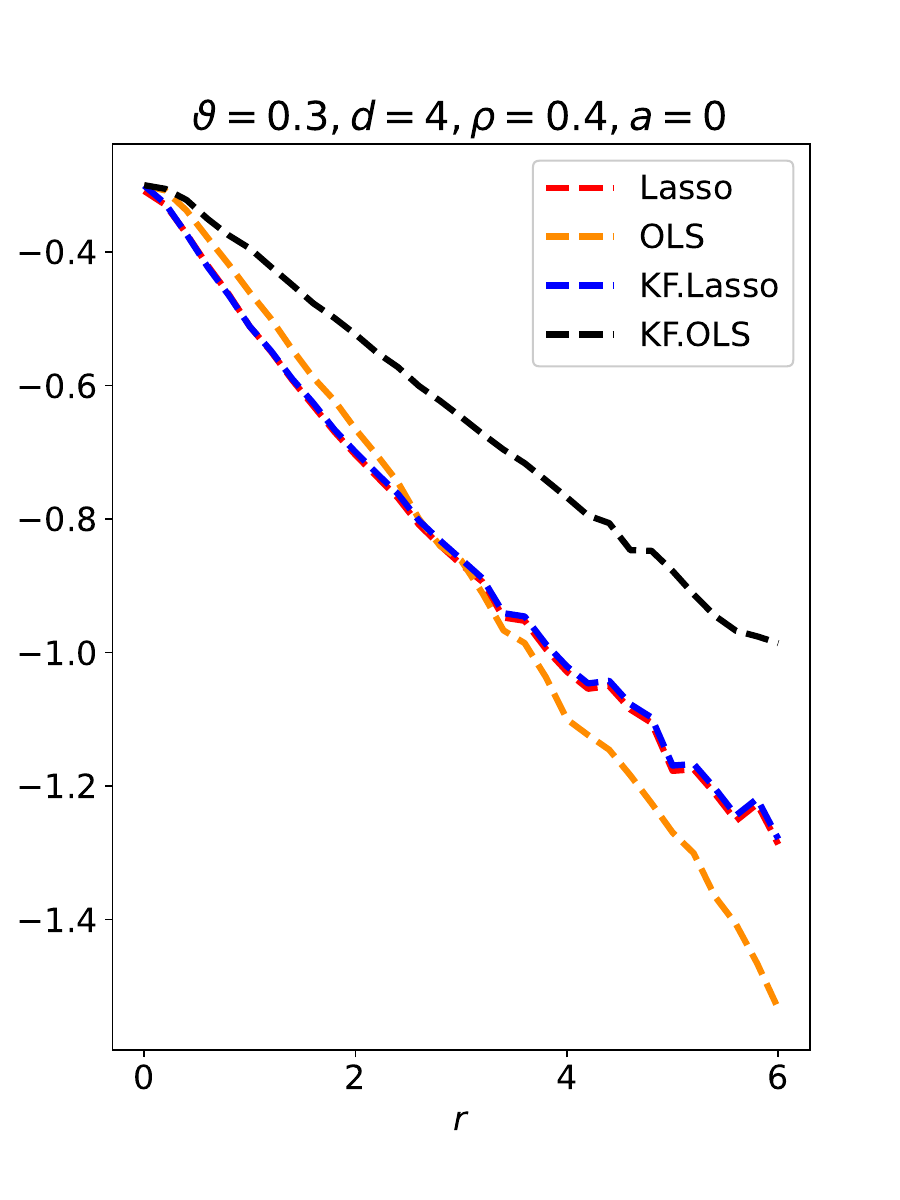} $\qquad$
\hspace*{-1cm}
\includegraphics[height=.25\textwidth,  width=.25\textwidth, trim=0 35 0 45, clip=true]{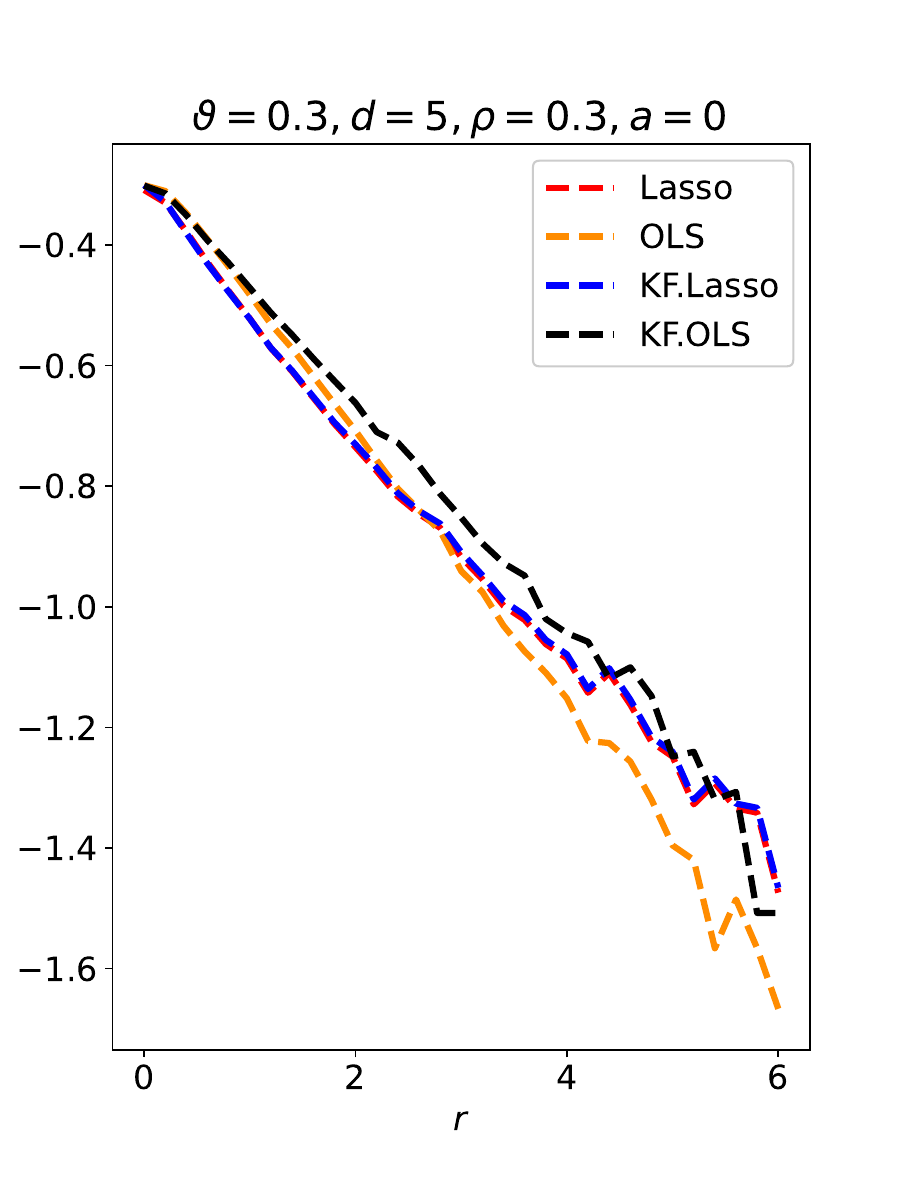}
\caption{Experiments 2 and 3 (block-wise diagonal designs, $d$: block size, $\rho$: off-diagonal entries). The y-axis is $\log_p(H_p/p)$, where $H_p$ is the average Hamming error. The parameter $a$ controls the construction of knockoff. 
}\label{fig:exp23}
\end{figure}

{\bf Experiment 2}. 
We consider the block-wise diagonal design with $2\times 2$ blocks, where we take $\rho=0.5$ and $\rho=0.7$. In the data generation, we fix an $n\times p$ matrix $X$ such that $X'X$ has the desirable form. We then generate $(\beta, y)$ in the same way as before. For each $\rho$, we 
fix $(n,p,\vartheta)=(2000,1000,0.2)$, and let $r$ range on a grid from $0$ to $8$ with a step size $0.2$. For KF.Lasso and KF.OLS, we now fix the symmetric statistic as signed maximum and the default choice of $\mathrm{diag}(s)$ yields that $\mathrm{diag}(s)=(1-a)I_p$ with $a=2\rho-1$. In this case, $G^*=[X,\tilde{X}]'[X,\tilde{X}]$ is degenerated, thus an $\epsilon=10^{-5}$ was subtracted from each elements of $\mathrm{diag}(s)$ to ensure KF.OLS is applicable.
The results are in the first two panels of Figure~\ref{fig:exp23}.

The theory in Section~\ref{sec:tamperDesign} suggests that since the two values of $\rho$ considered here are in $(\rho_0,1)$, KF.Lasso has a similar performance as its prototype, Lasso. While according to Section~\ref{sec:RankingAlg}, KF.OLS has a inferior performance comparing to its prototype, OLS. The simulation results are consistent with these theoretical predictions. 
Moreover, we can see that,  for the current $\vartheta$ value, OLS has a smaller Hamming error than that of Lasso when $r$ is large, and the opposite is true when $r$ is small. These also agree with our theory. 

\vspace{0.15in}

{\bf Experiment 3}. We further consider blockwise diagonal designs with larger-size blocks. Given $d\geq 2$ and $p$ that is a multiple of $d$, we generate $X\in\mathbb{R}^{n\times p}$ such that $X'X$ is block-wise diagonal with  $d\times d$ diagonal blocks, where the off-diagonal elements of each block are all equal to $\rho$. Other steps of the data generation are the same as in Experiment 2.  We consider $(d,\rho)=(4,0.4)$ and $(d, \rho)= (5,0.3)$.  For each choice of $(d, \rho)$, we set $(n,p,\vartheta)=(2000,1000,0.3)$ and let $r$ range on a grid from $0$ to $6$ with a step size $0.2$. We use signed maximum as symmetric statistic in KF and use the equi-correlated knockoff described above. The results are in the last two panels of Figure~\ref{fig:exp23}.

One noteworthy observation is that KF.Lasso still has a similar performance as its own prototype. Meanwhile, KF.OLS can get close to its prototype in the case where $\rho$ is close to 0. Another observation is that OLS outperforms Lasso when $r$ is large, and Lasso slightly outperforms OLS when $r$ is small. While our theory is only derived for $d=2$, the simulations suggest that similar insight continues to apply when the block size gets larger. 

\vspace{0.1in}

\begin{figure}[tb]
\centering
\hspace*{-.8cm}
\includegraphics[height=.25\textwidth, width=.25\textwidth, trim=0 35 0 45, clip=true]{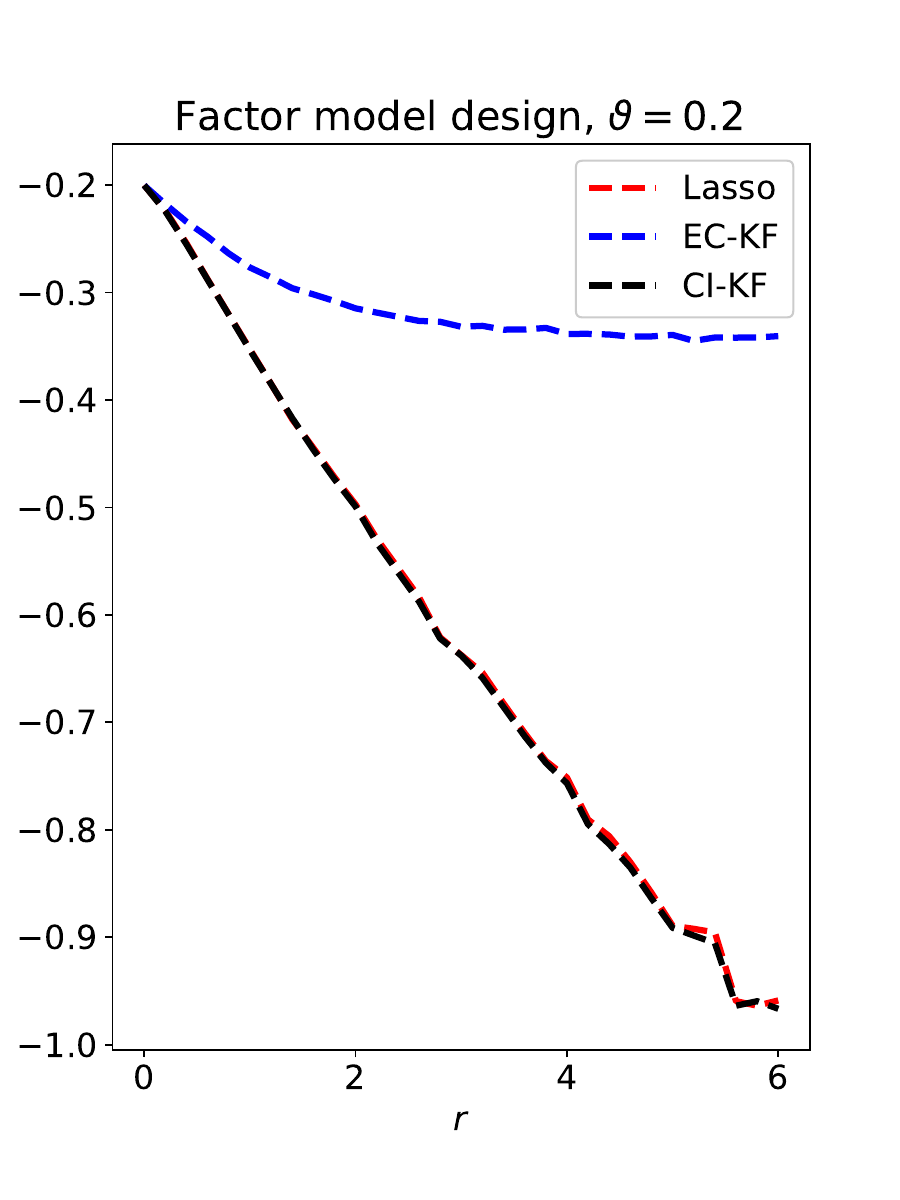} $\qquad$
\hspace*{-1cm}
\includegraphics[height=.25\textwidth,  width=.25\textwidth, trim=0 35 0 45, clip=true]{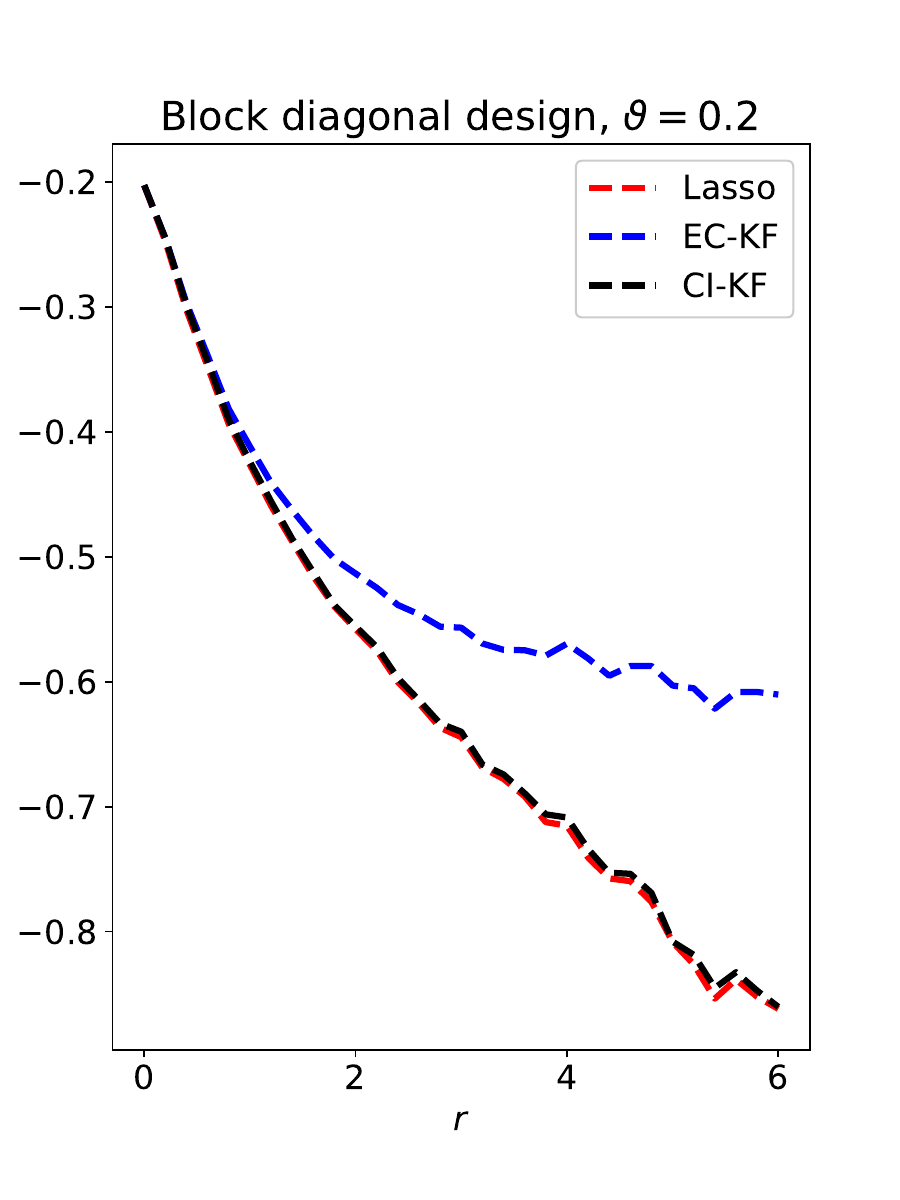} $\qquad$
\hspace*{-1cm}
\includegraphics[height=.25\textwidth,  width=.25\textwidth, trim=0 35 0 45, clip=true]{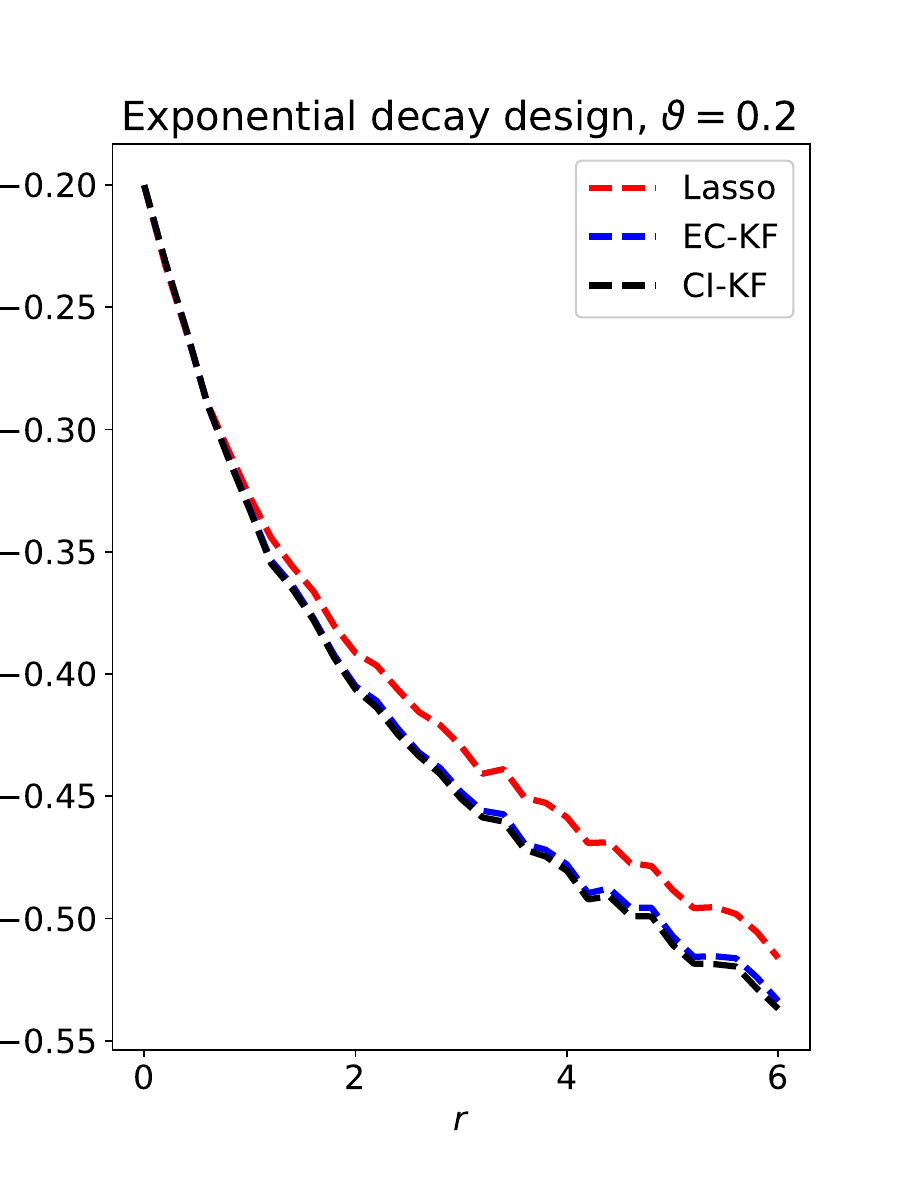} $\qquad$
\hspace*{-1cm}
\includegraphics[height=.25\textwidth,  width=.25\textwidth, trim=0 35 0 45, clip=true]{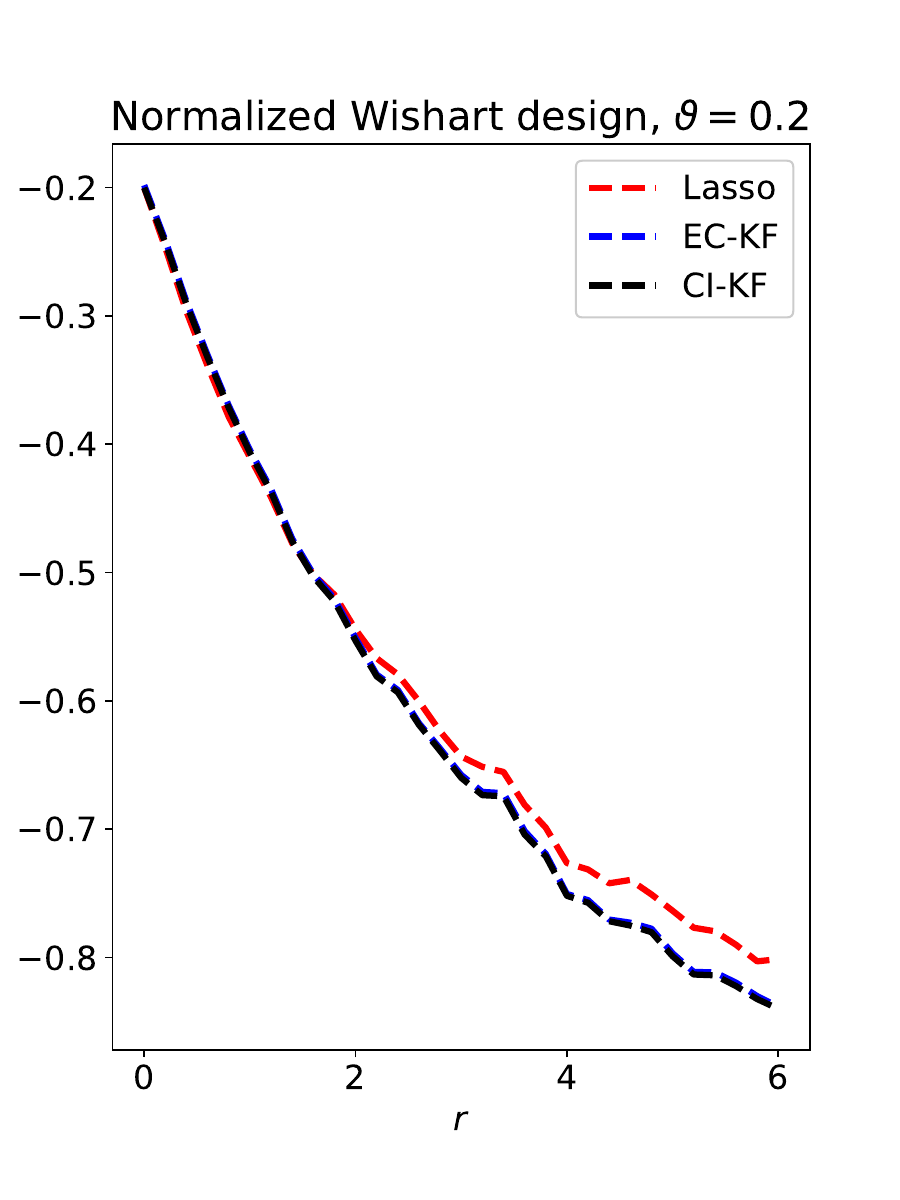}\\
\hspace*{-.8cm}
\includegraphics[height=.25\textwidth, width=.25\textwidth, trim=0 35 0 45, clip=true]{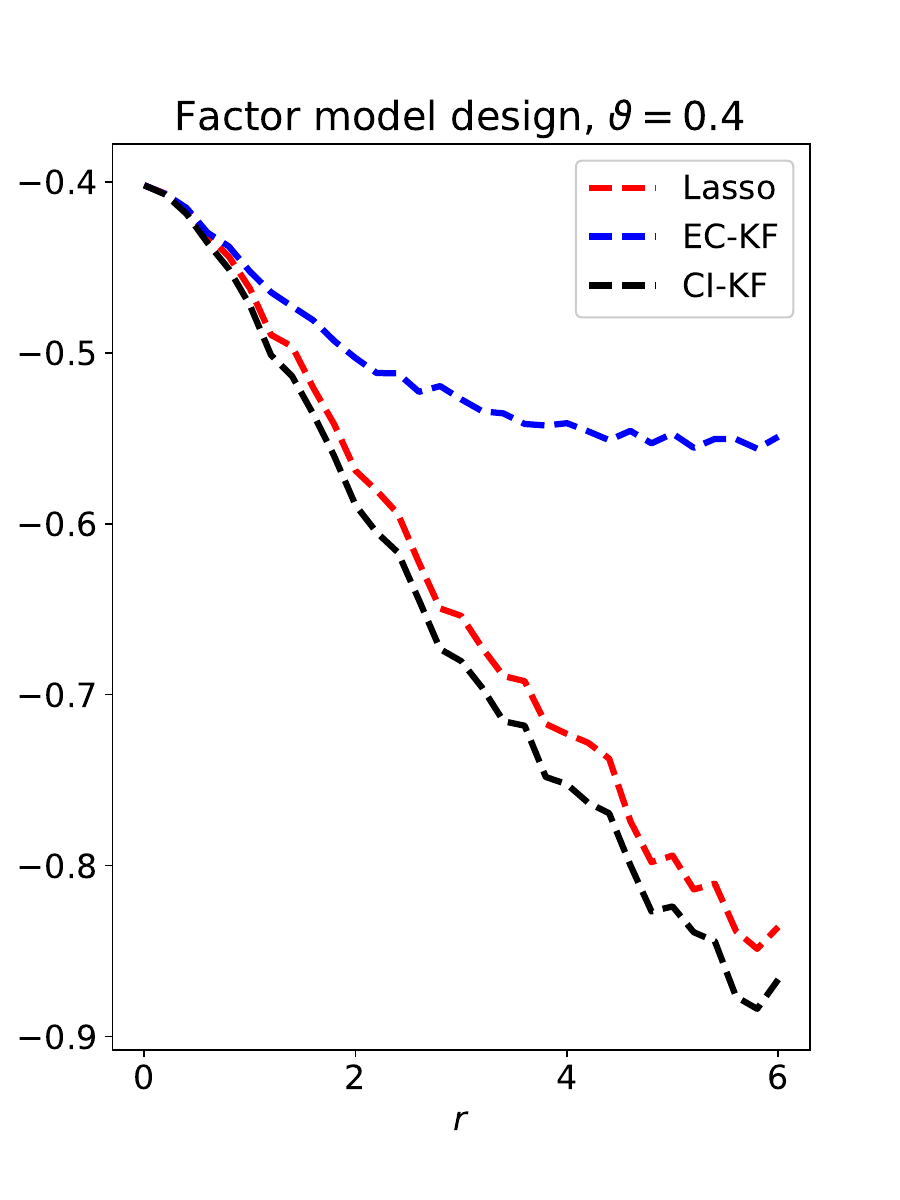} $\qquad$
\hspace*{-1cm}
\includegraphics[height=.25\textwidth,  width=.25\textwidth, trim=0 35 0 45, clip=true]{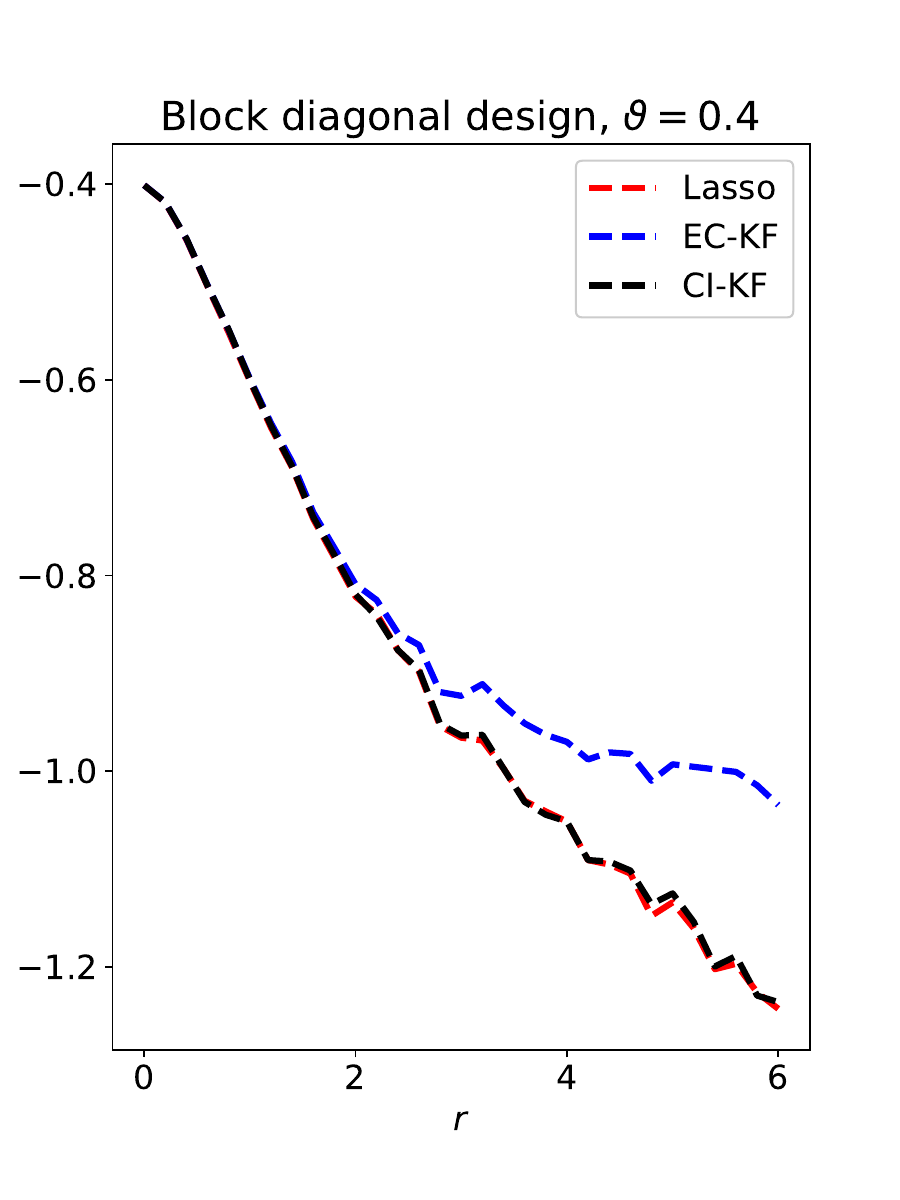} $\qquad$
\hspace*{-1cm}
\includegraphics[height=.25\textwidth,  width=.25\textwidth, trim=0 35 0 45, clip=true]{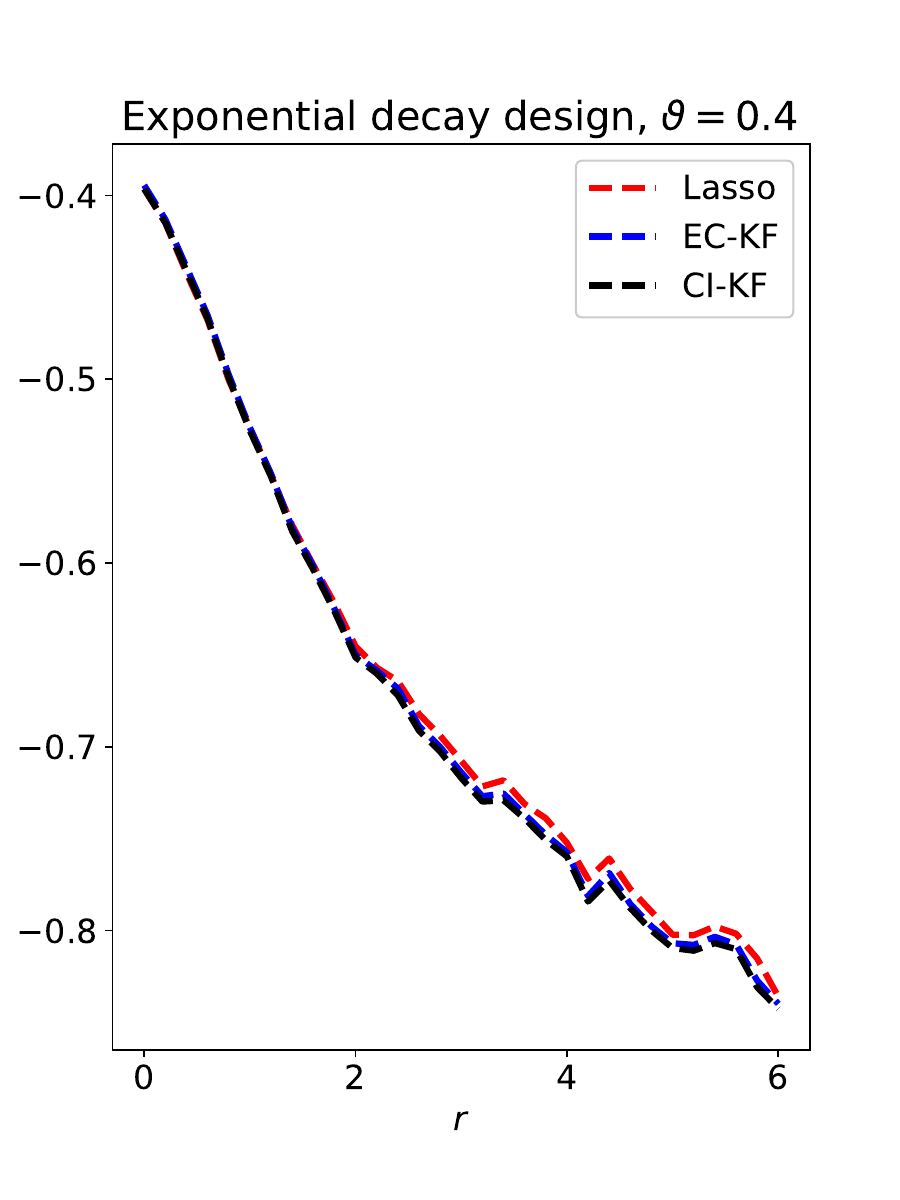} $\qquad$
\hspace*{-1cm}
\includegraphics[height=.25\textwidth,  width=.25\textwidth, trim=0 35 0 45, clip=true]{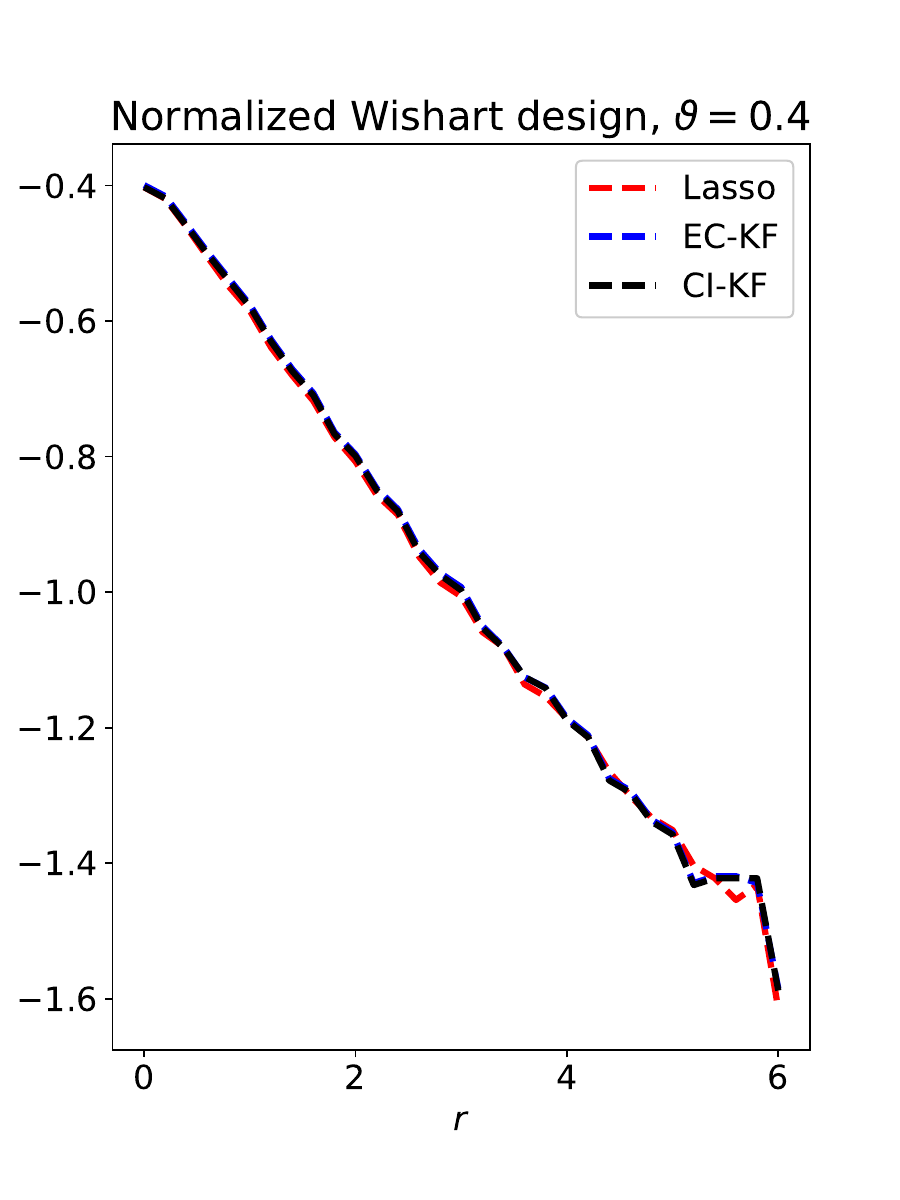}

\caption{Experiment 4 (general designs). The y-axis is $\log_p(H_p/p)$, where $H_p$ is the average Hamming error. 
We focus on comparing two constructions of knockoff's, EC-KF and CI-KF, and include Lasso as the benchmark. 
}\label{fig:exp4}
\end{figure}

{\bf Experiment 4}.
In Section~\ref{sec:tamperDesign}, we studied variants of knockoff with different augmented designs. The theory for $2\times 2$ block-wise designs suggests that using  CI-knockoff to construct $\tilde{X}$ yields a higher power than using EC-knockoff (for $2\times 2$ block-wise design, EC-knockoff is the same as SDP-knockoff). In this experiment, we investigate whether using CI-knockoff still yields a power boost for other design classes.   
We consider 4 types of designs:
\begin{itemize} \itemsep 0pt
    \item {\it Factor models}: $X'X=(BB'+I_p)/2$, where $B$ is a $p\times 2$ matrix whose $j$-th row is equal to $[\cos(\alpha_j), \sin(\alpha_j)]$ with $\{\alpha_j\}_{j=1,\cdots,p}$ $iid$ drawn from  $\mathrm{Uniform}[0,2\pi]$;
    \item {\it Block diagonal}: Same as in Experiment 2, where $\rho= 0.5$. 
    \item {\it Exponential decay}: The $(i,j)$-th element of $X'X$ is $0.6^{|i-j|}$, for $1\leq i,j\leq p$. 
    \item {\it Normalized Wishart}: $X'X$ is the sample correlation matrix of $n$ $iid$ samples of $N(0, I_p)$.
\end{itemize}
In the normalized Wishart design, the CI-knockoff in \eqref{CI-knockoff} may not satisfy $\mathrm{diag}(s)\preceq 2G$. We modify it to $\text{diag}(s)=\alpha[\text{diag}(G^{-1})]^{-1}$, where $\alpha$ is the maximum value in $[0,1]$ such that $\mathrm{diag}(s)\preceq 2G$. For each design, we fix $(n,p)=(1000,300)$, let $\vartheta$ take values in $ \{0.2,0.4\}$ and let $r$ range on a grid from $0$ to $6$ with a step size $0.2$. 
Different from previous experiments, we generate $\beta$ from $
\beta_j \overset{iid}{\sim} (1-\epsilon_p)\nu_0 + \frac{1}{2} \epsilon_p\nu_{\tau_p} + \frac{1}{2}\epsilon_p \nu_{-\tau_p}$, for $1\leq j\leq p$.
The motivation of using this model is to allow for negative entries in $\beta$. Even when $X'X$ contains only nonnegative elements, this signal model can still reveal the effect of having negative correlations in the design. We compare two versions of knockoff, EC-knockoff and CI-knockoff, along with the prototype, Lasso. The results are in Figure~\ref{fig:exp4}.

For the $2\times 2$ block-wise diagonal design, the simulations suggest that CI-KF significantly outperforms EC-KF, and that CI-KF has a similar performance as the prototype, Lasso. This is consistent with the theory in Section~\ref{subsec:FDR-lasso} and Section~\ref{subsec:CI-knockoff}. 
CI-KF also yields a significant improvement over EC-KF in the factor design, and the two methods perform similarly in the exponentially decaying design and the normalized Wishart design. We notice that the Gram matrix of the normalized Wishart design has uniformly small off-diagonal entries for the current $(n,p)$, which is similar to the orthogonal design and  explains why EC-KF and CI-KF do not have much difference. Combining these simulation results, we recommend CI-KF for practical use. Additionally, in some settings (e.g., factor design, $\vartheta=0.4$; exponentially decaying design, $\vartheta=0.2$), CI-KF even outperforms its prototype Lasso. One possible reason is that the ideal threshold we use is derived by ignoring the multi-$\log(p)$ term, but this term can have a non-negligible effect for a moderately large $p$, so the Hamming error of Lasso presented here may be larger than the actual optimal one.

\section{Discussions}  \label{sec:discuss}
How to maximize the power when controlling FDR at a targeted level is a problem of great interest. We focus on the FDR control method, knockoff, and point out that it has three key components: ranking algorithm, augmented design, and symmetric statistic. 
Since each component admits multiple choices, knockoff has many different variants. All the variants guarantee finite-sample FDR control. Our goal is to understand which variants enjoy good power. In a Rare/Weak signal model, for each variant of knockoff under consideration, we derive explicit forms of false positive rate and false negative rate, and obtain the theoretical phase diagram. 
The results provide useful guidelines of choosing the version of knockoff to use in practice. 
We also define the prototype of knockoff, which uses only one component, ranking algorithm, and has access to an ideal threshold. We compare the phase diagram of knockoff with the phase diagram of prototype. 
The results help us understand the extra price we pay for finding a data-driven threshold to control FDR.


We have several notable discoveries: (i) For the choice of symmetric statistic, signed maximum is better than difference, because the latter has an inferior phase diagram in the orthogonal design. (ii) For the choice of augmented design, CI-knockoff is better than SDP-knockoff, because the latter has an inferior phase diagram in a simple blockwise diagonal design. (ii) For the choice of ranking algorithm, roughly, Lasso is better than least-squares when the signals are extremely sparse and the design correlations are moderate; and least-squares is better than Lasso when the signals are only moderately sparse and the design correlations are more severe. (iv) In a simple blockwise diagonal design, when knockoff uses Lasso as ranking algorithm, with proper choices of two other components, knockoff has the same phase diagram as its prototype (i.e., we pay a negligible  price for finding a data-driven threshold). This is however not true when knockoff uses least-squares as ranking algorithm. 



There are several directions to extend our current results. First, we focus on the regime where FDR and TPR converge to either 0 or 1 and characterize the rates of convergence. The more subtle regime where FDR and TPR converge to constants between 0 and 1 is not studied. We leave it to future work. Second, the study of knockoff here is only for block-wise diagonal designs. For general designs, it is very tedious to derive the precise phase diagram, but some cruder results may be less tedious to derive, such as an upper bound for the Hamming error. This kind of results will help shed more insights on how to construct the knockoff variables (e.g., how to choose $\mathrm{diag}(s)$). Third, we only investigate Lasso-path or the least-squares as options of the ranking algorithm. It is interesting to study the power of FDR control methods based on other ranking algorithms, such as the marginal screening and iterative sure screening \citep{fan2008sure} and the covariance assisted screening \citep{ke2014covariance,ke2017covariate}. The covariance assisted screening was shown to yield optimal phase diagrams for a broad class of sparse designs; whether it can be developed into an FDR control method with ``optimal'' power remains unknown and is worth future study. Last, some FDR control methods may not fit exactly the unified framework here. For instance, the multiple data splits \citep{dai2020false} is a method that controls FDR through data splitting. We can similarly assess its power using the Rare/Weak signal model and phase diagram, except that we need to assume the rows of $X$ are $i.i.d.$ generated. We leave such study to future work.

\acks{The research of Jun S. Liu was partially supported by the NSF grant DMS-201541 and the NIH R01grant HG011485-01.  The research of Zheng T. Ke was partially supported by the NSF CAREER grant DMS-1943902.}


\newpage

\appendix


\section{Proof of Lemma~\ref{lem:tool}}
By definition of the multi-$\log(p)$ term, it suffices to show that, for every $\epsilon >0$, as $p\to\infty$, 
\beq \label{lem-tailprob-0}
p^{-\epsilon+b}\ \mathbb{P}(X_p \in S) \to 0, \qquad\mbox{and}\qquad p^{\epsilon+b}\ \mathbb{P}(X_p \in S) \to \infty.
\eeq

We introduce two sets $\underline{S}$ and $\overline{S}$ such that
\[
\underline{S} \subset S\subset \overline{S}. 
\]
Define $m(x)=(x-\mu)'\Sigma^{-1}(x-\mu)$ for any $x\in\mathbb{R}^d$. By definition,  $b=\inf_{x\in S}m(x)$. As a result, $m(x)\geq b$ for all $x\in S$. Define
\beq \label{lem-tailprob-upS}
 \overline{S} =\{x\in\mathbb{R}^p: m(x)\geq b\}. 
\eeq
Then, $S\subset  \overline{S}$. Furthermore, since $m(x)$ is a quadratic function and $b=\inf_{x\in S}m(x)$, given any $\epsilon>0$, there exists $x_0\in S$ such that 
\beq \label{lem-tailprob-1}
m(x_0)\leq b+\epsilon/8.
\eeq 
Note that \eqref{lem-tailprob-1} guarantees that $\|x_0-\mu\|$ is bounded. For any $x\in S$ and $\|x-x_0\|\leq 1$, 
\begin{align*}
|m(x)-m(x_0)| & \leq 2|(x-\mu)'\Sigma^{-1}(x-x_0)| + |(x-x_0)'\Sigma^{-1}(x-x_0)|\cr
&\leq 2\|x-\mu\|\|\Sigma^{-1}\|\cdot \|x-x_0\| + \|\Sigma^{-1}\|\|x-x_0\|^2\cr
&\leq C_1\|x-x_0\| + C_2\|x-x_0\|^2, 
\end{align*}
where $C_1$ and $C_2$ are positive constants that only depend on $(\mu,\Sigma, b,\epsilon)$. It follows that there exists a constant $\delta_1>0$ such that 
\beq \label{lem-tailprob-2}
x\in S, \quad \|x-x_0\|\leq \delta_1 \qquad \Longrightarrow\qquad  |m(x)-m(x_0)|\leq \epsilon/8. 
\eeq
Additionally, since $S$ is an open set and $x_0\in S$, there exists $\delta_2>0$, such that 
\[
\{x\in \mathbb{R}^d: \|x-x_0\|\leq \delta_2\}\subset S.
\]
Define
\beq \label{lem-tailprob-downS}
 \underline{S} = \{x\in \mathbb{R}^d: \|x-x_0\|\leq \delta\}, \qquad \mbox{where}\quad \delta = \min\{\delta_1, \delta_2\}. 
\eeq
It is easy to see that $ \underline{S}\subset S$. Additionally, in light of \eqref{lem-tailprob-1} and \eqref{lem-tailprob-2}, 
\beq \label{lem-tailprob-3}
m(x)\leq b+\epsilon/4, \qquad\mbox{for all }x\in  \underline{S}. 
\eeq
Since $\underline{S}\subset S\subset\overline{S}$, to show \eqref{lem-tailprob-0}, it suffices to show that
\beq \label{lem-tailprob-goal1}
p^{\epsilon+b}\ \mathbb{P}\bigl(X_p \in  \underline{S}\bigr)\to\infty
\eeq
and 
\beq \label{lem-tailprob-goal2}
p^{-\epsilon+b}\ \mathbb{P} \bigl(X_p \in \overline{S}\bigr)\to 0. 
\eeq

First, we show \eqref{lem-tailprob-goal1}. Let $f_p(x)$ denote the density of ${\cal N}_d\bigl(\mu_p, \; \frac{1}{2\log(p)}\Sigma_p\bigr)$. 
Write $m_p(x)=(x-\mu_p)'\Sigma_p^{-1}(x-\mu_p)$. It is seen that  
\beq \label{lem-tailprob-fp}
f_p (x) = \frac{[2\log(p)]^{d/2}}{(2\pi)^{d/2} |\det(\Sigma_d)|^{1/2}} \cdot p^{-m_p(x)}. 
\eeq
By direct calculations, 
\begin{align}  \label{lem-tailprob-4}
\mathbb{P}\bigl(X_p \in \underline{S}\; |\; \mu_p, \Sigma_p \bigr)
&= \frac{[2\log(p)]^{d/2}}{(2\pi)^{d/2} |\det(\Sigma_p)|^{1/2}}  \int_{x\in \underline{S}} p^{-m_p(x)} dx \cr
&\geq \frac{[\log(p)]^{d/2}}{\pi^{d/2} |\det(\Sigma_p)|^{1/2}}  \cdot \mathrm{Volume}(\underline{S}) \cdot p^{-\sup_{x\in \underline{S}}\{m_p(x)\}}. 
\end{align}
The assumptions on $(\mu_p, \Sigma_p)$ imply that, for any constant $\gamma>0$,
\[
\lim_{p\to\infty}\mathbb{P}\bigl(\|\mu_p-\mu\|>\gamma \mbox{ or }\|\Sigma_p-\Sigma\|>\gamma\bigr)= 0. 
\]
Let $E$ be the event that $\|\mu_p-\mu\|\leq \gamma_*$ and $\|\Sigma_p-\Sigma\|\leq \gamma_*$, for some $\gamma_*$ to be decided. On this event, for any $x\in \underline{S}$, 
\begin{align*}
|m(x)-m_p(x)| &\leq |(x-\mu)'\Sigma^{-1}(x-\mu)-(x-\mu)'\Sigma_p^{-1}(x-\mu)| \cr
&\qquad + |(x-\mu)'\Sigma_p^{-1}(x-\mu)-(x-\mu_p)'\Sigma_p^{-1}(x-\mu_p)|  \cr
&\leq |(x-\mu)'(\Sigma^{-1}-\Sigma_p^{-1})(x-\mu)| + 2|(x-\mu)'\Sigma_p^{-1}(\mu-\mu_p)|\cr
&\qquad + (\mu-\mu_p)'\Sigma_p^{-1}(\mu-\mu_p)\cr
&\leq \|x-\mu\|^2\|\Sigma^{-1}\|\|\Sigma_p^{-1}\|\cdot \|\Sigma_p-\Sigma\| + 2\|x-\mu\|\|\Sigma_p^{-1}\|\cdot \|\mu-\mu_p\| \cr
&\qquad + \|\Sigma_p^{-1}\|\cdot \|\mu-\mu_p\|^2\cr
&\leq C_3\gamma_* + C_4\gamma_*^2,
\end{align*}
where $C_3$ and $C_4$ are positive constants that do not depend on $\gamma_*$, and  in the last line we have used the fact that $\underline{S}$ is a bounded set so that $\|x-\mu\|$ is bounded. It follows that we can choose an appropriately small $\gamma_*$ such that
\begin{align}  \label{lem-tailprob-5}
|m(x)-m_p(x)|\leq \epsilon/4, \qquad \mbox{for all }x\in \underline{S}. 
\end{align}
Combining \eqref{lem-tailprob-5} with \eqref{lem-tailprob-3} gives 
\[
\sup_{x\in \underline{S}} m_p(x)\leq b+\epsilon/2, \qquad \mbox{on the event }E. 
\]
Moreover, since $\underline{S}$ is a ball with radius $\delta$,
\[
\mathrm{Volume}(\underline{S}) =\delta^d \cdot \mathrm{Volume}(B_d), 
\]
where $B_d$ is the unit ball in $\mathbb{R}^d$, whose volume is a constant. We plug the above results into \eqref{lem-tailprob-4} and notice that $|\det(\Sigma_p)|\geq |\det(\Sigma)|-C_5\delta$ on the event $E$, for a constant $C_5>0$. It yields that, when $(\mu_p,\Sigma_p)$ satisfies the event $E$, 
\beq \label{lem-tailprob-6}
\mathbb{P}\bigl(X_p \in \underline{S}\; |\; \mu_p, \Sigma_p \bigr)\geq c_0 [\log(p)]^{d/2}\cdot p^{-(b+\epsilon/2)},
\eeq
for some constant $c_0>0$. 
It follows that 
\[
\mathbb{P}\bigl(X_p \in \underline{S}\bigr)\geq  \mathbb{P}(E)\cdot  c_0 [\log(p)]^{d/2} p^{-(b+\epsilon/2)}.
\]
We plug it into the left hand side of \eqref{lem-tailprob-goal1} and note that $\mathbb{P}(E)\to 1$ as $p\to\infty$. This gives the desirable claim in \eqref{lem-tailprob-goal1}. 

Next, we show \eqref{lem-tailprob-goal2}. We define a counterpart of the set $\overline{S}$ by 
\[
\overline{S}_p = \{x\in\mathbb{R}^d: m_p(x)\geq b\}.  
\]
Define $Y_p=\sqrt{2\log(p)}\cdot \Sigma_p^{-1/2}( X_p-\mu_p)$. Then, $Y_p\sim {\cal N}_d(0, I_d)$ and 
\[
X_p \in \overline{S}_p \qquad\mbox{if and only if}\qquad \|Y_p\|^2 \geq 2b\log(p). 
\]
The distribution of $\|Y_p\|^2$ is a $\chi^2_d$ distribution, which does not depend on $(\mu_p,\Sigma_p)$. We have 
\begin{align} \label{lem-tailprob-7}
\mathbb{P}\bigl(X_p \in \overline{S}_p\bigr)&=\mathbb{E}\bigl[
\mathbb{P}\bigl(X_p \in \overline{S}_p\; |\;\mu_p, \Sigma_p \bigr)\bigr]\cr
&= \mathbb{E}\bigl[  \mathbb{P}\bigl(\|Y_p\|^2 \geq 2b\log(p)\bigr) \bigr]\cr
&= \mathbb{P}\bigl(\chi^2_d \geq 2b\log(p)\bigr). 
\end{align} 
For chi-square distribution, the tail probability has an explicit form:
\[
\mathbb{P}\bigl(\chi^2_d \geq 2b\log(p)\bigr) = \frac{\Gamma(d/2,\ b\log(p))}{\Gamma(d/2)}, 
\]
where $\Gamma(s,x)\equiv\int_x^{\infty} t^{s-1}\exp(-t)dt$ is the upper incomplete gamma function and $\Gamma(s)\equiv\Gamma(s,0)$ is the ordinary gamma function. By property of the upper incomplete gamma function, $\Gamma(s,x)/(x^{s-1}\exp(-x))\to 1$ as $x \to \infty$.  It follows that 
\[
 \frac{\Gamma(d/2,\ b\log(p))}{[b\log(p)]^{d/2-1}p^{-b}}\to 1, \qquad \mbox{as}\quad p\to\infty. 
\]
In particular, when $p$ is sufficiently large, the left hand side is $\geq 1/2$. We plug these results into \eqref{lem-tailprob-7} to get
\beq    \label{lem-tailprob-8}
\mathbb{P}\bigl(X_p \in \overline{S}_p\bigr)\geq \frac{[b\log(p)]^{d/2-1}}{2\Gamma(d/2)}\cdot p^{-b}. 
\eeq
It remains to study the difference caused by replacing $\overline{S}_p$ by $\overline{S}$. Let 
\[
U_p = (\overline{S}\backslash\overline{S}_p)\cup(\overline{S}_p\backslash\overline{S}).
\] Then,
\beq \label{lem-tailprob-8(2)}
\bigl|\mathbb{P}\bigl(X_p \in \overline{S}\bigr)- \mathbb{P}\bigl(X_p \in \overline{S}_p\bigr)\bigr|\leq \mathbb{P}\bigl(X_p \in U_p\bigr). 
\eeq
Similar to \eqref{lem-tailprob-4}, we have 
\begin{align}  \label{lem-tailprob-9}
\mathbb{P}\bigl(X_p  \in U_p \; |\;\mu_p, \Sigma_p \bigr)&
= \frac{[2\log(p)]^{d/2}}{(2\pi)^{d/2} |\det(\Sigma_p)|^{1/2}}  \int_{x\in U_p } p^{-m_p(x)} dx \cr
&\leq \frac{[\log(p)]^{d/2}}{\pi^{d/2} |\det(\Sigma_p)|^{1/2}}  \cdot \mathrm{Volume}(U_p) \cdot p^{-\inf_{x\in U_p}\{m_p(x)\}}. 
\end{align}
For a constant $\gamma>0$ to be decided, let $F$ be the event that 
\beq  \label{lem-tailprob-10}
\|\mu_p-\mu\|\leq \gamma, \qquad \mbox{and}\quad \|\Sigma_p- \Sigma\|\leq \gamma.
\eeq
On this event, we study both $\mathrm{Volume}(U_p)$ and $\inf_{x\in U_p}m_p(x)$. Re-write 
\[
U_p=(\overline{S}^c\backslash\overline{S}^c_p)\cup(\overline{S}^c_p\backslash\overline{S}^c).
\]
By definition, $\overline{S}^c=\{x\in\mathbb{R}^d: m(x) \leq b\}=\{x\in\mathbb{R}^d: \|\Sigma^{-1/2}(x-\mu)\| \leq \sqrt{b}\}$, and $\overline{S}_p^c=\{x\in\mathbb{R}^d: \|\Sigma_p^{-1/2}(x-\mu_p)\| \leq \sqrt{b}\}$. On the event $F$, for any $x\in\overline{S}^c_p$, 
\begin{align*}
\|\Sigma^{-1/2}(x-\mu)\| & \leq \sqrt{b} + \|\Sigma^{-1/2}(x-\mu)-\Sigma_p^{-1/2}(x-\mu_p)\|\cr
&\leq \sqrt{b} + \|\Sigma^{-1/2}(\mu_p-\mu)\| + \|(\Sigma^{-1/2}-\Sigma_p^{-1/2})(x-\mu_p)\|\cr
&\leq  \sqrt{b} + \|\Sigma^{-1/2}\|\cdot \|\mu_p-\mu\| + \|\Sigma^{1/2}\Sigma_p^{-1/2}-I_d\|\cdot \|\Sigma_p^{-1/2}(x-\mu_p)\|\cr
&\leq \sqrt{b}+\|\Sigma^{-1/2}\|\cdot \|\mu_p-\mu\| + \sqrt{b}\cdot \|\Sigma^{1/2}\Sigma_p^{-1/2}-I_d\|\cr
&\leq \sqrt{b}+C_5\gamma, 
\end{align*}
for a constant $C_5>0$ that does not depend on $\gamma$. Choosing $\gamma<C_5^{-1}\sqrt{b}$, we have $\|\Sigma^{-1/2}(x-\mu)\|\leq 2\sqrt{b}$ for all $x\in \overline{S}_p^c$. Additionally, by definition, $\|\Sigma^{-1/2}(x-\mu)\|\leq\sqrt{b}$ for all $x\in\overline{S}^c$. Combining the above gives
\[
U_p\; \subset\;  (\overline{S}^c \cup \overline{S}^c_p)
\; \subset\; \bigl\{x\in\mathbb{R}^d: \|\Sigma^{-1/2}(x-\mu)\|\leq 2\sqrt{b}\bigr\}.  
\]
Recall that $B_d$ is the unit ball in $\mathbb{R}^d$. It follows immediately that 
\beq   \label{lem-tailprob-11}
\mathrm{Volume}(U_p)\leq (2\sqrt{b})^d\cdot \mathrm{Volume}(B_d), \qquad\mbox{on the event $F$}. 
\eeq
At the same time, for any $x\in \overline{S}$, on the event $F$,
\begin{align*}
\|\Sigma_p^{-1/2}(x-\mu_p)\| & \geq \|\Sigma^{-1/2}(x-\mu)\| - \|\Sigma_p^{-1/2}(x-\mu_p)-\Sigma^{-1/2}(x-\mu)\|\cr
&\geq  \|\Sigma^{-1/2}(x-\mu)\| - \|\Sigma_p^{-1/2}(\mu_p-\mu)\| - \|(\Sigma^{-1/2}-\Sigma_p^{-1/2})(x-\mu)\|\cr
&\geq  \|\Sigma^{-1/2}(x-\mu)\| - \|\Sigma_p^{-1/2}\|\cdot \|\mu_p-\mu\| - \|\Sigma_p^{-1/2}\Sigma^{1/2}-I_d\|\cdot \|\Sigma^{-1/2}(x-\mu)\|\cr
&= \|\Sigma^{-1/2}(x-\mu)\|\bigl(1-\|\Sigma_p^{-1/2}\Sigma^{1/2}-I_d\|\bigr)   -\|\Sigma^{-1/2}\|\cdot \|\mu_p-\mu\|\cr
&\geq \|\Sigma^{-1/2}(x-\mu)\|(1-C_6\gamma) - \|\Sigma^{-1/2}\|\gamma\cr
&\geq \sqrt{b}(1-C_6\gamma) - \|\Sigma^{-1/2}\|\gamma, 
\end{align*}
where $C_6>0$ is a constant that does not depend on $\gamma$ and in the last line we have used the fact that $\|\Sigma^{-1/2}(x-\mu)\|\geq \sqrt{b}$ for $x\in \overline{S}$. 
We choose $\gamma$ properly small so that $\sqrt{b}(1-C_6\gamma)-\|\Sigma^{-1/2}\|\gamma\geq \sqrt{b-\epsilon/2}$. It follows that
\beq  \label{lem-tailprob-12}
m_p(x)= \|\Sigma_p^{-1/2}(x-\mu_p)\|^2\geq b-\epsilon/2, \qquad\mbox{for all }x\in \overline{S}.  
\eeq
Additionally, the definition of $\overline{S}_p$ already guarantees that $m_p(x)\geq b$ for all $x\in \overline{S}_p$. Consequently,
\beq \label{lem-tailprob-13}
\inf_{x\in U_p}m_p(x) \geq \inf_{x\in \overline{S}\cup \overline{S}_p}\{ m_p(x)\} \geq b-\epsilon/2, \qquad\mbox{on the event $F$}. 
\eeq
We plug \eqref{lem-tailprob-11} and  \eqref{lem-tailprob-13} into \eqref{lem-tailprob-9}. It yields that, on the event $F$, 
\beq 
\mathbb{P}\bigl(X_p \in U_p \;|\;\mu_p, \Sigma_p \bigr)\leq C_7[\log(p)]^{d/2}\cdot p^{-(b-\epsilon/2)},
\eeq 
for a constant $C_7>0$. Then, 
\[
\mathbb{P}\bigl(X_p\in U_p\bigr) \leq \mathbb{P}(F)\cdot C_7[\log(p)]^{d/2}\cdot p^{-(b-\epsilon/2)} + \mathbb{P}(F^c). 
\]
By our assumption, for any $\gamma>0$ and $L>0$, $\mathbb{P}(\|\mu_p-\mu\|>\gamma)\leq p^{-L}$ and $\mathbb{P}(\|\Sigma_p-\Sigma\|>\gamma)\leq p^{-L}$. In particular, we can choose $L=b$. It gives 
\[
\mathbb{P}(F^c)\leq p^{-b}.
\]
We combine the above results and plug them into \eqref{lem-tailprob-8(2)}. It follows that 
\beq  \label{lem-tailprob-14}
\bigl|\mathbb{P}\bigl(X_p \in \overline{S}\bigr)- \mathbb{P}\bigl(X_p \in \overline{S}_p\bigr)\bigr|\leq C_7[\log(p)]^{d/2}\cdot p^{-(b-\epsilon/2)} + p^{-b}. 
\eeq
Combining \eqref{lem-tailprob-8} and \eqref{lem-tailprob-14} gives
\[
\mathbb{P}\bigl(X_p/ \in \overline{S}\bigr)\leq [1+o(1)]\cdot C_7[\log(p)]^{d/2}\cdot p^{-(b-\epsilon/2)}. 
\]
This gives the claim in \eqref{lem-tailprob-goal2}. The proof of this lemma is complete.

\section{Proof of Lemma~\ref{lem:block-RejectRegion}}
First, we study the least-squares. Note that $\hat{\beta}$ has an explicit solution: $\hat{\beta}=G^{-1}X^Ty$. Since $G$ is a block-wise diagonal matrix, we immediately have
\[
\begin{bmatrix}
\hat{\beta}_{j}\\
 \hat{\beta}_{j+1}\end{bmatrix}
  = \matB^{-1}\begin{bmatrix}x_j^Ty\\
   x_{j+1}^Ty\end{bmatrix}=\frac{1}{1-\rho^2}\begin{bmatrix}
   x_j^Ty-\rho x_{j+1}^Ty\\
   x_{j+1}^Ty-\rho x_{j+1}^Ty\end{bmatrix}. 
\]
Recall that $\tilde{y}=X'y/\sqrt{2\log(p)}$. Then, $|\hat{\beta}_j|>\sqrt{2u\log(p)}$ if and only if
\[
\frac{1}{1-\rho^2}|\tilde{y}_j-\rho \tilde{y}_{j+1}|>\sqrt{u}. 
\] 
It immediately gives the rejection region for least-squares. 

Next, we study the Lasso-path. We write $W_j^{*,\text{path}}$ as $W_j^*$ for notation simplicity. The lasso estimate $\hat{\beta}(\lambda)$ minimizes the objective
\[
Q(b) =\frac{1}{2}\|y-Xb\|^2+\lambda\|b\|_1 = \frac{1}{2}\|y\|^2 -y^TXb + \frac{1}{2}b^TG b + \lambda\|b\|_1. 
\]
When $G$ is a block-wise diagonal matrix, the objective $Q(b)$ is separable, and we can optimize over each pair of $(b_j, b_{j+1})$ separately. It reduces to solving many bi-variate problems:
\begin{equation}\label{bivarLasso}
    (\hat{\beta}_j(\lambda),\hat{\beta}_{j+1}(\lambda))^T =\mathrm{argmin}_b\Big\{\frac{1}{2}||y-[x_j,x_{j+1}]b||_2^2+\lambda ||b||_1\Big\}.
\end{equation}
Write $\hat{b}=(\hat{\beta}_j(\lambda),\hat{\beta}_{j+1}(\lambda))^T$ and let 
\[
B= \begin{bmatrix} 1 \;\;\; & \rho\\\rho\;\;\; & 1\end{bmatrix}\qquad \mbox{and}\qquad 
h = \begin{bmatrix}
x_j^Ty \\ x_{j+1}^Ty
\end{bmatrix}. 
\]
Then, the optimization \eqref{bivarLasso} can be written as
\beq \label{bivarLasso2}
\hat{b} = \mathrm{argmin}_b \bigl\{ - h^Tb + b^TBb/2 + \lambda\|b\|_1\bigr\}. 
\eeq
Recall that $W_j^*$ is the value of $\lambda$ at which $\hat{b}_1$ becomes nonzero for the first time. Our goal is to find a region of $(h_1, h_2)$ such that $W_j^*>t_p(u)\equiv \sqrt{2u\log(p)}$.

It suffices to consider the case of $\rho\geq 0$. To see this, 
we consider changing $\rho$ to $-\rho$ in the matrix $B$. The objective remains unchanged if we also change $b_2$ to $-b_2$ and $h_2$ to $-h_2$. Note that the change of $b_2$ to $-b_2$ has no impact on $W_j^*$; this means $W_j^*$ is unchanged if we simultaneously flip the sign of $\rho$ and $h_2$. Consequently, once we know the rejection region for $\rho>0$, we can immediately obtain that for $\rho<0$ by a reflection of the region with respect to the x-axis. 

\begin{figure}[t]
\includegraphics[height=.5\textwidth, trim=300 0 100 0, clip=true]{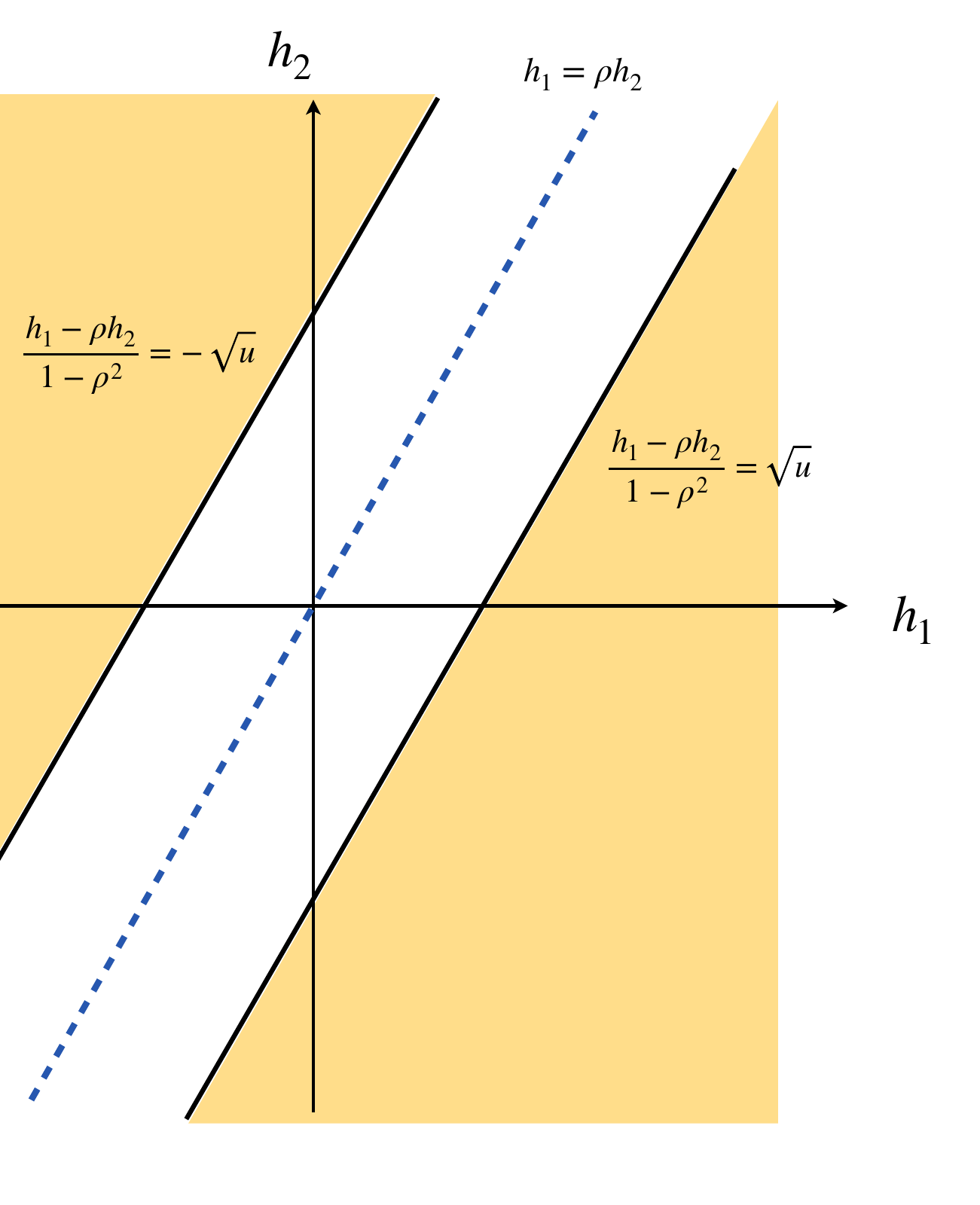}
\includegraphics[height=.5\textwidth,  trim=100 0 100 0, clip=true]{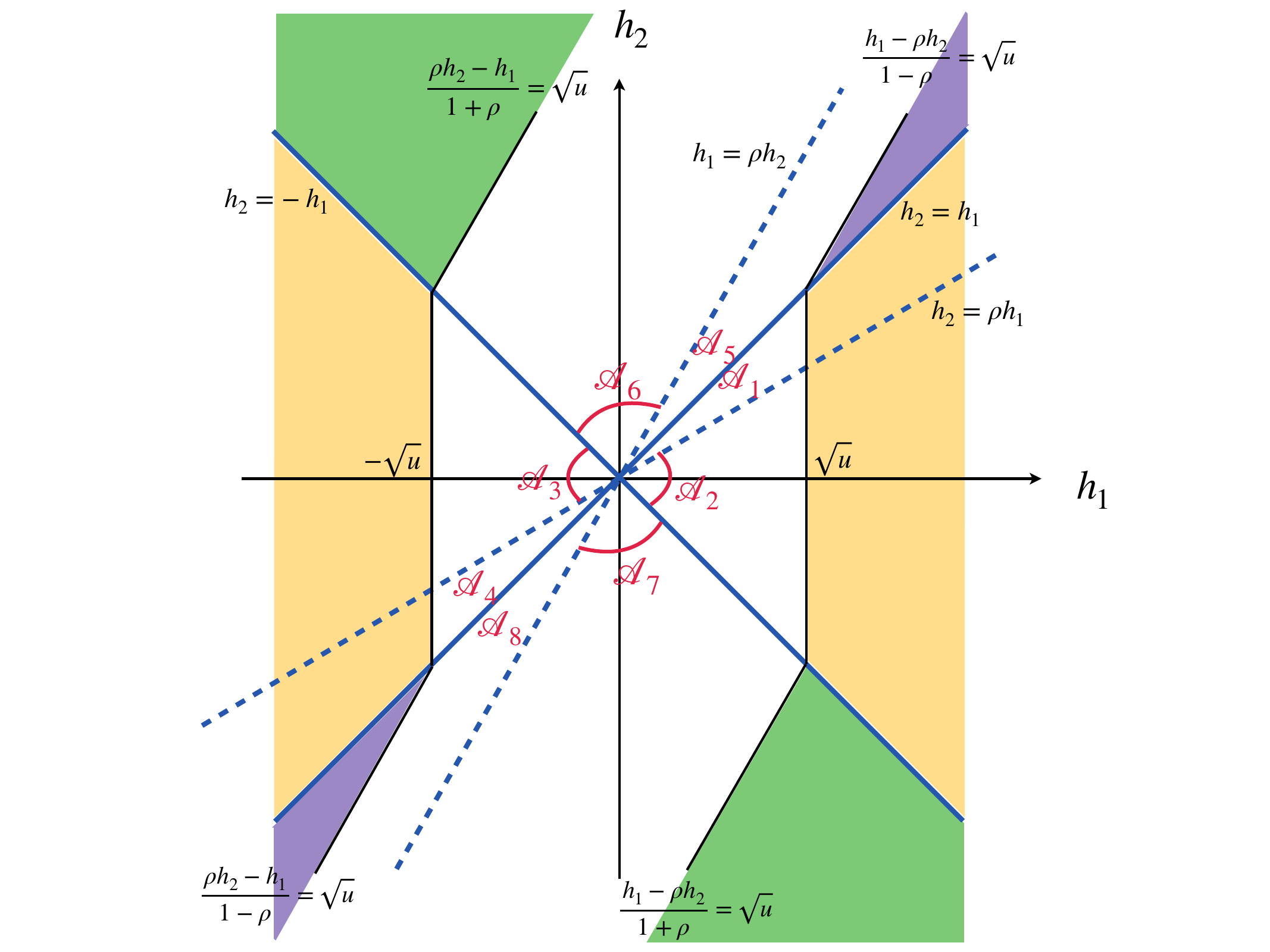}
\caption{The rejection region of least-squares (left) and Lasso-path (right). On the right panel, the regions ${\cal A}_1$-${\cal A}_8$ are the same as those defined in the proof. In the regions ${\cal A}_1$-${\cal A}_4$, $W_j^*=|h_1|$, and the rejection region is colored by yellow. In the regions ${\cal A}_5$ and ${\cal A}_8$, $W_j^*=|h_1-\rho h_2|/(1-\rho)$, and the rejection region is colored by purple. In the regiions ${\cal A}_6$ and ${\cal A}_7$, $W_j^*=|h_1-\rho h_2|/(1+\rho)$, and the rejection region is colored by green.} \label{fig:region-proof}
\end{figure}

Below, we fix $\rho\geq 0$. We first derive the explicit form of the whole solution path and then use it to decide the rejection region. Taking sub-gradients of \eqref{bivarLasso}, we find that $\hat{b}$ has to satisfy 
\beq \label{lem-region-1}
\begin{bmatrix} 1 \;\;\; & \rho\\\rho\;\;\; & 1\end{bmatrix} 
\begin{bmatrix} \hat{b}_1\\\hat{b}_2 
\end{bmatrix} + \lambda
\begin{bmatrix}
\sgn(\hat{b}_1) \\ \sgn(\hat{b}_2)
\end{bmatrix} =\begin{bmatrix}h_1\\h_2\end{bmatrix}, 
\eeq
where $\sgn(x)=1$ if $x>0$, $\sgn(x)=-1$ if $x<0$, and $\sgn(x)$ can be equal to any value in $[-1,1]$ if $x=0$. Let $\lambda_1>\lambda_2>0$ be the values at which variables enter the solution path. When $\lambda\in (\lambda_1,\infty)$, $\hat{b}_1=0$ and $\hat{b}_2=0$. Plugging them into \eqref{lem-region-1} gives $\mathrm{sgn}(\hat{b}_1)=\lambda^{-1}h_1$. The definition of $\mathrm{sgn}(\hat{b}_1)$ implies that $|h_1|\leq\lambda$, for any $\lambda>\lambda_1$. We then have $|h_1|\leq \lambda_1$. Similarly, it is true that $|h_2|\leq \lambda_1$. It gives 
\beq \label{lem-region-2}
\lambda_1 = \max\{|h_1|, |h_2|\}. 
\eeq
We first assume $|h_1|>|h_2|$. 
By \eqref{lem-region-1} and continuity of solution path, there exists a sufficiently small constant $\delta>0$ such that, for $\lambda\in (\lambda_2-\delta,\lambda_2)$, the following equation holds. 
\beq \label{lem-region-3}
\matB \begin{bmatrix}\hat{b}_1(\lambda)\\\hat{b}_2(\lambda) \end{bmatrix}+\lambda\begin{bmatrix}\mathrm{sgn}(\hat{b}_1)\\\mathrm{sgn}(\hat{b}_2)\end{bmatrix} = \begin{bmatrix}h_1\\h_2\end{bmatrix}. 
\eeq
The sign vector of $\hat{b}$ for $\lambda\in (\lambda_2-\delta,\lambda_2)$ has four different cases: $(1,1)^T$, $(1, -1)^T$, $(-1,1)^T$, $(-1,-1)^T$. For these four different cases, we can use \eqref{lem-region-3} to solve $\hat{b}$. The solutions in four cases are respectively 
\begin{align*}
& \frac{1}{1-\rho^2}\begin{bmatrix}(h_1-\rho h_2)-(1-\rho)\lambda\\(h_2-\rho h_1)-(1-\rho)\lambda \end{bmatrix},\qquad \frac{1}{1-\rho^2}\begin{bmatrix}(h_1-\rho h_2)-(1+\rho)\lambda\\(h_2-\rho h_1)+(1+\rho)\lambda \end{bmatrix}, \cr
&\frac{1}{1-\rho^2}\begin{bmatrix}(h_1-\rho h_2)+(1+\rho)\lambda\\(h_2-\rho h_1)-(1+\rho)\lambda \end{bmatrix}, \qquad \frac{1}{1-\rho^2}\begin{bmatrix}(h_1-\rho h_2)+(1-\rho)\lambda\\(h_2-\rho h_1)+(1-\rho)\lambda \end{bmatrix}. 
\end{align*}
The solution $\hat{b}$ has to match the sign assumption on $\hat{b}$. For each of the four cases, the requirement becomes
\begin{itemize}
\item Case 1:\; $(h_1-\rho h_2)-(1-\rho)\lambda>0$, \; $(h_2-\rho h_1)-(1-\rho)\lambda >0$.
\item Case 2:\; $(h_1-\rho h_2)-(1+\rho)\lambda>0$, \; $(h_2-\rho h_1)+(1+\rho)\lambda <0$.
\item Case 3:\; $(h_1-\rho h_2)+(1+\rho)\lambda<0$, \; $(h_2-\rho h_1)-(1+\rho)\lambda >0$.
\item Case 4:\; $(h_1-\rho h_2)+(1-\rho)\lambda<0$, \; $(h_2-\rho h_1)+(1-\rho)\lambda <0$.
\end{itemize} 
Note that we have assumed $|h_1|>|h_2|$. Then, Case $k$ is possible only in the region ${\cal A}_k$, where
\begin{align*}
& {\cal A}_1 = \{(h_1, h_2): h_1>0,\;\; \rho h_1 < h_2  < h_1\}, \quad
{\cal A}_2 = \{(h_1, h_2): h_1>0,\;\; -h_1 < h_2  <\rho h_1\},\cr
&{\cal A}_3 = \{(h_1, h_2): h_1<0, \;\; \rho h_1 <h_2<- h_1\}, \quad
{\cal A}_4 =  \{(h_1, h_2): h_1<0, \;\; h_1 < h_2  <\rho h_1\}. 
\end{align*}
In each case, $\lambda_1=|h_1|$. To get the value of $\lambda_2$, we use the continuity of the solution path. It implies that $\hat{b}_2(\lambda)=0$ at $\lambda=\lambda_2$. As a result, the value of $\lambda_2$ in Case $k$ is 
\beq \label{lem-region-lambda-1to4}
\lambda_2^{(1)}=\frac{h_2-\rho h_1}{1-\rho}, \qquad \lambda_2^{(2)}=\frac{\rho h_1-h_2}{1+\rho}, \qquad \lambda_2^{(3)}= \frac{h_2-\rho h_1}{1+\rho}, \qquad \lambda_2^{(4)}= \frac{\rho h_1-h_2}{1-\rho}.  
\eeq
It is easy to verify that $\lambda_2<\lambda_1$ in each case. We also need to check that in the region ${\cal A}_k$, the KKT condition \eqref{lem-region-1} can be satisfied with $\hat{b}_2=0$ for all $\lambda\in (\lambda^{(k)}_2,\lambda_1)$. For example, in Case 1, \eqref{lem-region-1} becomes
\[
\matB \begin{bmatrix}\hat{b}_1\\0 \end{bmatrix}+\lambda\begin{bmatrix}1\\c \end{bmatrix} = \begin{bmatrix}h_1\\h_2\end{bmatrix}, \qquad\mbox{for some }|c|\leq 1. 
\]
We can solve the equations to get $\hat{b}_1=h_1-\lambda$ and $\lambda c=h_2-\rho\hat{b}_1=(h_2-\rho h_1)-\lambda $. It can be verified that $|(h_2-\rho h_1)-\lambda|\leq \lambda$ for $(h_1,h_2)\in {\cal A}_1$ and $\lambda\in (\lambda^{(1)}_2, \lambda_1)$. The verification for other cases is similar and thus omitted. We then assume $|h_2|>|h_1|$. By symmetry, we will have the same result, except that $(h_1, h_2)$ are switched in the expression of ${\cal A}$ and $(\lambda_1, \lambda_2)$. This gives the other four cases:
\begin{align*}
& {\cal A}_5 = \{(h_1, h_2): h_2>0,\;\; \rho h_2 < h_1  < h_2\}, \quad
{\cal A}_6 = \{(h_1, h_2): h_2>0,\;\; -h_2 < h_1  <\rho h_2\},\cr
&{\cal A}_7 = \{(h_1, h_2): h_2<0, \;\; \rho h_2 <h_1<- h_2\}, \quad
{\cal A}_8  =  \{(h_1, h_2): h_2<0, \;\; h_2 < h_1  <\rho h_2\}. 
\end{align*}
In these four cases, we similarly have $\lambda_1=|h_2|$ and 
\beq  \label{lem-region-lambda-5to8}
\lambda_2^{(5)}=\frac{h_1-\rho h_2}{1-\rho}, \qquad \lambda_2^{(6)}=\frac{\rho h_2-h_1}{1+\rho}, \qquad \lambda_2^{(7)}= \frac{h_1-\rho h_2}{1+\rho}, \qquad \lambda_2^{(8)}= \frac{\rho h_2-h_1}{1-\rho}.  
\eeq
These eight regions are shown in Figure~\ref{fig:region-proof}. 

We then compute $W_j^*$ and the associated rejection region. 
Note that $W_j^*=\lambda_1$ in Case 1-Case 4, and $W_j^*=\lambda_2$ in Case 5-Case 8. It follows directly that 
\beq
W_j^*= \begin{cases}
|h_1|, & \mbox{if } (h_1, h_2)\in  {\cal A}_1\cup {\cal A}_2\cup {\cal A}_3\cup {\cal A}_4,\\
|h_1-\rho h_2|/(1-\rho), &\mbox{if }(h_1, h_2)\in {\cal A}_5\cup {\cal A}_8,\\
|h_1-\rho h_2|/(1+\rho), &\mbox{if }(h_1, h_2)\in {\cal A}_6\cup {\cal A}_7. 
\end{cases}
\eeq
As a result, the region $W_j^*>\sqrt{2u\log(p)}$ if and only if the vector $(x_j^Ty, x_{j+1}^Ty)/\sqrt{2\log(p)}$ is in the following set:
\begin{align*}
{\cal R} &= \bigl[({\cal A}_1\cup {\cal A}_2\cup {\cal A}_3\cup {\cal A})\cap \{|h_1|>\sqrt{u}\}\bigr]\cr
&\qquad \cup \bigl[({\cal A}_5\cup {\cal A}_8)\cap \{ |h_1-\rho h_2|>(1-\rho)\sqrt{u} \}\bigr]\cr
&\qquad \cap \bigl[({\cal A}_6\cup {\cal A}_7)\cap \{ |h_1-\rho h_2|>(1+\rho)\sqrt{u} \}\bigr]. 
\end{align*}
In Figure~\ref{fig:region-proof}, the 3 subsets are colored by yellow, purple, and green, respectively. 
This gives the rejection region for Lasso-path.

\section{Proof of Theorem~\ref{thm:knockoff}}
By definition of $(\FP_p, \FN_p)$ and the Rare/Weak signal model \eqref{RWmodel1}-\eqref{RWmodel2}, we have
\beq \label{knockoff-ortho-start}
\FP_p=\sum_{j=1}^p (1-\epsilon_p) \mathbb{P}(W_j>t_p(u)|\beta_j=0), \quad 
\FN_p =\sum_{j=1}^p \epsilon_p\ \mathbb{P}(W_j<t_p(u)|\beta_j=\tau_p),  
\eeq
where $\epsilon_p=p^{-\vartheta}$, $\tau_p=\sqrt{2r\log(p)}$, and $t_p(u)=\sqrt{2u\log(p)}$. Therefore, it suffices to study $\mathbb{P}(W_j>t_p(u)|\beta_j=0)$ and $\mathbb{P}(W_j<t_p(u)|\beta_j=\tau_p)$. 

Fix $1\leq j\leq p$. The knockoff filter applies Lasso to the design matrix $[X,\tilde{X}]$. This design is belongs to the block-wise diagonal design \eqref{block} with a dimension $2p$ and $\rho=a$. The variable $j$ and its own knockoff are in one block.  Fix $j$ and write
\beq \label{thm-knockoff-ortho-00}
h_1 = x_j'y/\sqrt{2\log(p)}, \qquad \mbox{and}\qquad h_2=\tilde{x}_j'y/\sqrt{2\log(p)}.
\eeq
It is easy to see that $(x_j'y, \tilde{x}_j'y)'$ follows a distribution ${\cal N}_2({\bf0}_2, \Sigma)$ when $\beta_j=0$, and it follows a distribution ${\cal N}_2(\mu\sqrt{2\log(p)}, \; \Sigma)$, when $\beta_j=\tau_p$, where 
\[
\mu =\begin{bmatrix}\sqrt{r}\\a\sqrt{r}\end{bmatrix}, \qquad \Sigma = 
\matBa. 
\]
Let ${\cal R}$ be the region of $(h_1, h_2)$ corresponding to the event that $\{W_j>t_p(u)\}$. It follows from Lemma~\ref{lem:tool} that
\begin{align} \label{thm-knockoff-ortho-0}
\mathbb{P}(W_j>t_p(u)|\beta_j=0)& =L_pp^{-\inf_{h\in {\cal R}}\{ h'\Sigma^{-1}h\}}, \cr \mathbb{P}(W_j<t_p(u)|\beta_j=\tau_p)& =L_pp^{-\inf_{h\in {\cal R}^c} \{(h-\mu)'\Sigma^{-1}(h-\mu)\}}. 
\end{align}
Below, we first derive the rejection region ${\cal R}$, and then compute the exponents in \eqref{thm-knockoff-ortho-0}. 

Recall that $Z_j$ and $\tilde{Z}_j$ are the same as in \eqref{Wj}. They are indeed the values of $\lambda$ at which the variable $j$ and its knockoff enter the solution path of a bivariate lasso as in \eqref{bivarLasso}. We can apply the solution path derived in the proof of Lemma~\ref{lem:block-RejectRegion}, with $\rho=a$. Before we proceed to the proof, we  argue that it suffices to consider the case of $a\geq 0$. If $a<0$, we can simultaneously flip the signs of $a$ and $h_2$, so that the objective \eqref{bivarLasso} remains unchanged; as a result, the values of $(Z_j,\tilde{Z}_j)$ remain unchanged, so does the symmetric statistic $W_j$. It implies that, if we flip the sign of $a$,  the rejection region is reflected with respect to the x-axis. At the same time, in light of the exponents in \eqref{thm-knockoff-ortho-0}, we consider two ellipsoids
\beq \label{ellipsoids}
{\cal E}_{\mathrm{FP}}(t)= \{h\in\mathbb{R}^2: h'\Sigma^{-1}h\leq t\}, \qquad {\cal E}_{\mathrm{FN}}(t)= \{h\in\mathbb{R}^2: (h-\mu)'\Sigma^{-1}(h-\mu)'\leq t\}. 
\eeq
Similarly, if we simultaneously flip the signs of $a$ and $h_2$, these ellipsoids remain unchanged. It implies that, if we flip the sign of $a$, these ellipsoids are reflected with respect to the x-axis. Combining the above observations, we know that the exponents in \eqref{thm-knockoff-ortho-0} are unchanged with a sign flip of $a$, i.e., they only depend on $|a|$. We assume $a\geq 0$ without loss of generality. 

\begin{figure}[t]
\includegraphics[height=.34\textwidth, trim=80 100 100 50, clip=true]{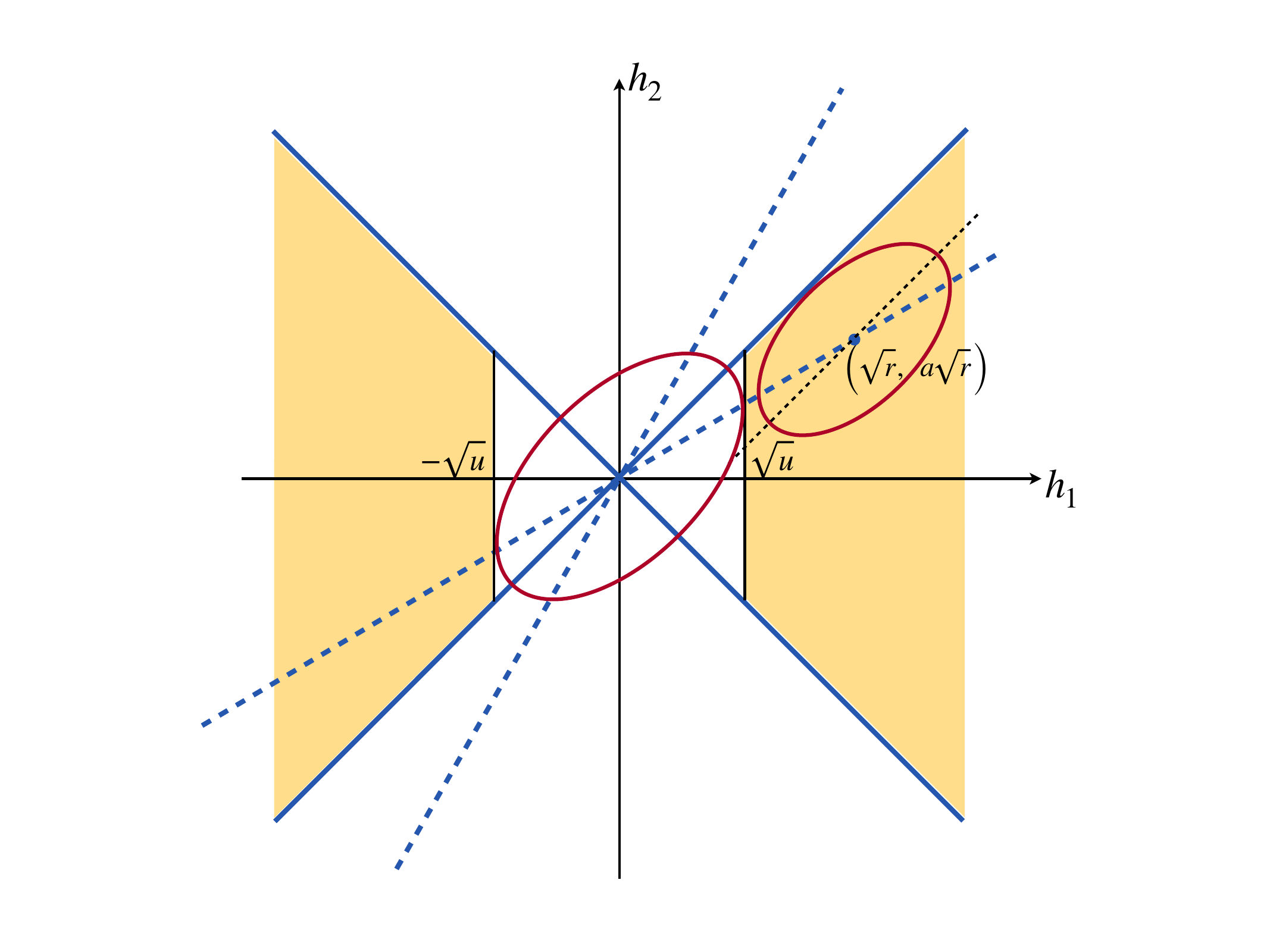}
\includegraphics[height=.34\textwidth,  trim=80 100 100 50, clip=true]{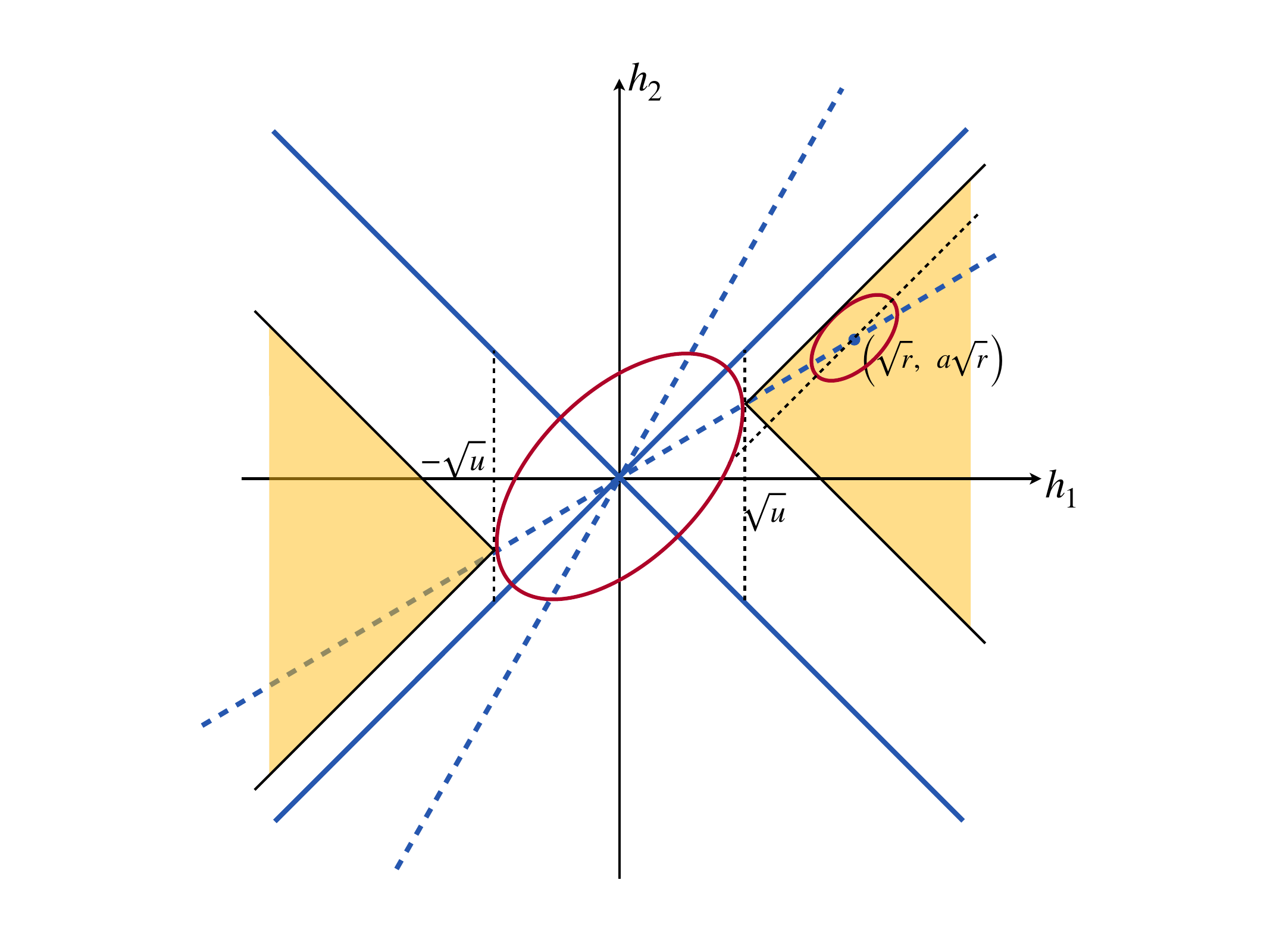}
\caption{The rejection region of knockoff in the orthogonal design, where the symmetric statistic is signed maximum (left) and difference (right). The rate of convergence of $\FP_p$ is captured by an ellipsoid centered at $(0,0)$, and the rate of convergence of $\FN_p$ is captured by an ellipsoid centered at $(\sqrt{r}, a\sqrt{r})$.} \label{fig:region-proof-2}
\end{figure}

Fix $a\geq 0$. Write $z=Z_j/\sqrt{2\log(p)}$ and $\tilde{z}=\tilde{Z}_j/\sqrt{2\log(p)}$. The symmetric statistics in \eqref{Wj} can be re-written as 
\[
W_j^{\mathrm{sgm}} = (z \vee \tilde{z})\sqrt{2\log(p)}\cdot
\begin{cases}
+1,&\text{if } z >\tilde{z} \\
-1, &\text{if } z \leq \tilde{z}
\end{cases}, 
\qquad W_j^{\mathrm{dif}} = (z-\tilde{z}) \sqrt{2\log(p)}. 
\]
Recall that $h_1$ and $h_2$ are as in \eqref{thm-knockoff-ortho-00}. Let $\lambda_1>\lambda_2>0$ be the values of $\lambda$ at which variables enter the solution path of a bivariate lasso. In the proof of Lemma~\ref{lem:block-RejectRegion}, we have derived the formula of $(\lambda_1, \lambda_2)$; see \eqref{lem-region-lambda-1to4} and \eqref{lem-region-lambda-5to8} (with $\rho$ replaced by $a$). It follows that
\[
(z, \tilde{z})=\begin{cases}
(\lambda_1, \lambda_2), & \mbox{in the regions ${\cal A}_1$-${\cal A}_4$},\\
(\lambda_2, \lambda_1), & \mbox{in the regions ${\cal A}_5$-${\cal A}_8$},  
\end{cases}
\]
where regions ${\cal A}_1$-${\cal A}_8$ are the same as those on the right panel of Figure~\ref{fig:region-proof} (with $\rho$ replaced by $a$). Plugging in \eqref{lem-region-lambda-1to4} and \eqref{lem-region-lambda-5to8} gives the following results: 
\begin{itemize}
\item Region ${\cal A}_1$:\; $z=h_1$, \;  $\tilde{z}=\frac{h_2-\rho h_1}{1-a}$, \;  $W_j^{\mathrm{sgm}}=h_1\sqrt{2\log(p)}$, \;  $W_j^{\mathrm{dif}}=\frac{h_1-h_2}{1-a}\sqrt{2\log(p)}$. 
\item Region ${\cal A}_2$: \;  $z=h_1$,\; $\tilde{z}=\frac{\rho h_1-h_2}{1+a}$, \;$W_j^{\mathrm{sgm}}=h_1\sqrt{2\log(p)}$, \;  $W_j^{\mathrm{dif}}=\frac{h_1+h_2}{1+a}\sqrt{2\log(p)}$. 
\item Region ${\cal A}_3$: \;  $z=-h_1$,\; $\tilde{z}=\frac{h_2-\rho h_1}{1+a}$, \; $W_j^{\mathrm{sgm}}=-h_1\sqrt{2\log(p)}$, \; $W_j^{\mathrm{dif}}=-\frac{h_1+h_2}{1+a}\sqrt{2\log(p)}$. 
\item Region ${\cal A}_4$: \; $z=-h_1$,\;  $\tilde{z}=\frac{\rho h_1-h_2}{1-a}$, \; $W_j^{\mathrm{sgm}}=-h_1\sqrt{2\log(p)}$, \; $W_j^{\mathrm{dif}}=\frac{h_2-h_1}{1-a}\sqrt{2\log(p)}$. 
\item Regions ${\cal A}_5$-${\cal A}_8$: \; $|Z_j|<|\tilde{Z}_j|$,\;  $W_j^{\mathrm{sgm}}<0$,\;  $W_j^{\mathrm{dif}}<0$. 
\end{itemize}
The event that $W_j^{\mathrm{sgm}}>\sqrt{2u\log(p)}$ corresponds to that $(h_1, h_2)$ is in the region of 
\begin{align} \label{lem-knockoff-ortho-1}
{\cal R}_u^{\mathrm{sgm}} &= ({\cal A}_1\cup {\cal A}_2\cup {\cal A}_3\cup {\cal A}_4)\cap \{|h_1|>\sqrt{u}\}\cr
&=\{|h_1|>|h_2|, \; |h_1|>\sqrt{u}\}. 
\end{align}
The event that $W_j^{\mathrm{dif}}>\sqrt{2u\log(p)}$ corresponds to that $(h_1, h_2)$ is in the region of 
\begin{align} \label{lem-knockoff-ortho-2}
{\cal R}_u^{\mathrm{dif}} &= \bigl({\cal A}_1\cap \{h_1-h_2>(1-a)\sqrt{u}\}\bigr)\cup \bigl({\cal A}_2\cap \{h_1+h_2>(1+a)\sqrt{u}\}\bigr)\cr
&\qquad \cup \bigl({\cal A}_3\cap \{h_1+h_2<-(1+a)\sqrt{u}\}\bigr)\cup \bigl({\cal A}_4\cap \{h_1-h_2<-(1-a)\sqrt{u}\}\bigr). 
\end{align}
These two regions are shown in Figure~\ref{fig:region-proof-2}. 

We are now ready to compute the exponents in \eqref{thm-knockoff-ortho-0}. 
First, we compute $\inf_{h\in {\cal R}}\{h'\Sigma^{-1}h\}$. Let ${\cal E}_{\mathrm{FP}}(t)$ be the same as in \eqref{ellipsoids}. Then, 
\[
\inf_{h\in {\cal R}}\{h'\Sigma^{-1}h\}=\sup\bigl\{t>0: {\cal E}_{\mathrm{FP}}(t)\cap {\cal R}\neq \emptyset\bigr\}. 
\]
When the rejection region is ${\cal R}_u^{\mathrm{sgm}}$, from Figure~\ref{fig:region-proof-2}, we can increase $t$ until ${\cal E}_{\mathrm{FP}}(t)$ intersects with the line of $h_1=\pm \sqrt{u}$. For any $h$ on the surface of this ellipsoid, the perpendicular vector of its tangent plane is proportional to $'\Sigma^{-1}h$. When the ellipsoid intersects with the line of $h_1=\pm \sqrt{u}$, the perpendicular vector should be proportional to $(1,0)'$. Therefore, we need to find $h$ such that 
\[
h_1=\pm \sqrt{u}, \quad h'\Sigma^{-1}h=t, \quad \mbox{and}\quad \Sigma^{-1}h\propto (1,0)'. 
\] 
The second equation requires that $h_2=a h_1$. Combining it with the first equation gives $h=(\pm \sqrt{u},\pm a\sqrt{u})$. We then plug it into the second equation to obtain $t=u$. This gives
\beq  \label{lem-knockoff-ortho-3}
\inf_{h\in {\cal R}_u^{\mathrm{sgm}}}\{h'\Sigma^{-1}h\} = u. 
\eeq
When the rejection region is ${\cal R}_u^{\mathrm{dif}}$, there are 3 possible cases: 
\begin{itemize}
\item[(i)] The ellipsoid intersects with the line $h_1-h_2=(1-a)\sqrt{u}$, 
\item[(ii)] The ellipsoid intersects with the line $h_1+h_2=(1+a)\sqrt{u}$, 
\item[(iii)] The ellipsoid intersects with the point $h=(\sqrt{u}, a\sqrt{u})$.
\end{itemize}
In Case (i), we can compute the intersection point by solving $h$ for $h_1-h_2=(1-a)\sqrt{u}$ and $\Sigma^{-1}h\propto (1,-1)'$. The second relationship gives $h_2=-h_1$. Together with the first relationship, we have $h=(\frac{1-a}{2}\sqrt{u}, \frac{1-a}{2}\sqrt{u})$. It is not in ${\cal R}_u^{\mathrm{dif}}$. Similarly, for Case (ii), we can show that the intersection point is $h=(\frac{1+a}{2}\sqrt{u}, \frac{1+a}{2}\sqrt{u})$, which is not in ${\cal R}_u^{\mathrm{dif}}$ either. The only possible case is Case (iii), where the intersection point is $(\sqrt{u}, a\sqrt{u})$ and the associated $t=h'\Sigma^{-1}t=u$.  
We have proved that
\beq  \label{lem-knockoff-ortho-4}
\inf_{h\in {\cal R}_u^{\mathrm{dif}}}\{h'\Sigma^{-1}h\} = u. 
\eeq
Next, we compute $\inf_{h\in {\cal R}^c}\{(h-\mu)'\Sigma^{-1}(h-\mu)\}$. Let ${\cal E}_{\mathrm{FN}}(t)$ be the same as in \eqref{ellipsoids}. Then, 
\[
\inf_{h\in {\cal R}^c}\{(h-\mu)'\Sigma^{-1}(h-\mu)\}=\sup\bigl\{t>0: {\cal E}_{\mathrm{FN}}(t)\cap {\cal R}^c \neq \emptyset\bigr\}. 
\]
Note that the center of the ellipsoid is $\mu=(\sqrt{r}, a\sqrt{r})$. When either ${\cal R}={\cal R}_u^{\mathrm{sgm}}$ or ${\cal R}={\cal R}_u^{\mathrm{dif}}$, $\mu\notin {\cal R}^c$ if and only if $r>u$. In other words, the above is well defined only if $r>u$. We now fix $r>u$. 
When the rejection region is ${\cal R}_u^{\mathrm{sgm}}$, the ellipsoid intersects with either the line of $h_1=\sqrt{u}$ or the line of $h_1=h_2$. Since the perpendicular vector of the tangent plane of the ellipsoid at $h$ is proportional to $'\Sigma^{-1}(h-\mu)$, we can solve the intersection points from
\[
\begin{cases}
h_1 = \sqrt{u},\\
\Sigma^{-1}(h-\mu)\propto (1,0)',
\end{cases}\qquad\mbox{and}\qquad  \begin{cases}
h_1 = h_2,\\
\Sigma^{-1}(h-\mu)\propto (1,-1)'. 
\end{cases}
\]
By calculations, the two intersection points are $h=(\sqrt{u},\; a\sqrt{u})$ and $h=(\frac{1+a}{2}\sqrt{r}, \; \frac{1+a}{2}\sqrt{r})$. The associated value of $(h-\mu)'\Sigma^{-1}(h-\mu)$ is $t=(\sqrt{r}-\sqrt{u})^2$ and $t=(1-a)r/2$, respectively. When we increase the ellipsoid until it interacts with $({\cal R}_u^{\mathrm{sgm}})^c$, the corresponding $t$ is the smaller of the above two values. This gives
\beq  \label{lem-knockoff-ortho-5}
\inf_{h\in ({\cal R}_u^{\mathrm{sgm}})^c}\{(h-\mu)'\Sigma^{-1}(h-\mu)\} = \min\Bigl\{ (\sqrt{r}-\sqrt{u})_+^2,\; \frac{1-a}{2}r\Bigr\}. 
\eeq
When the rejection region is ${\cal R}_u^{\mathrm{dif}}$, the ellipsoid intersects with either the line of $h_1-h_2=(1-a)\sqrt{u}$ or the line of $h_1+h_2=(1+a)\sqrt{u}$. We can solve the intersection points from
\[
\begin{cases}
h_1-h_2=(1-a)\sqrt{u},\\
\Sigma^{-1}(h-\mu)\propto (1,-1)',
\end{cases}\qquad\mbox{and}\qquad  \begin{cases}
h_1+h_2=(1+a)\sqrt{u},\\
\Sigma^{-1}(h-\mu)\propto (1,1)'.
\end{cases}
\]
Solving these equations gives the two intersection points: $h=(\frac{1+a}{2}\sqrt{r}+\frac{1-a}{2}\sqrt{u},\;  \frac{1+a}{2}\sqrt{r}-\frac{1-a}{2}\sqrt{u})$ and $h=(\frac{1-a}{2}\sqrt{r}+\frac{1+a}{2}\sqrt{u},\;  -\frac{1-a}{2}\sqrt{r}+\frac{1+a}{2}\sqrt{u})$. The corresponding value of $(h-\mu)'\Sigma^{-1}(h-\mu)$ is $t=\frac{1-a}{2}(\sqrt{r}-\sqrt{u})^2$ and $t=\frac{1+a}{2}(\sqrt{r}-\sqrt{u})^2$, respectively. The smaller of these two values is $\frac{1-a}{2}(\sqrt{r}-\sqrt{u})^2$. We have proved that
\beq  \label{lem-knockoff-ortho-6}
\inf_{h\in ({\cal R}_u^{\mathrm{dif}})^c}\{(h-\mu)'\Sigma^{-1}(h-\mu)\} = \frac{1-a}{2}(\sqrt{r}-\sqrt{u})^2_+. 
\eeq
We plug \eqref{lem-knockoff-ortho-3}-\eqref{lem-knockoff-ortho-6} into \eqref{thm-knockoff-ortho-0}, and we further plug it into \eqref{knockoff-ortho-start}. This gives the claim for $a\geq 0$. As we have argued, the results for $a<0$ only requires replacing $a$ by $|a|$.


\section{Proof of Theorem~\ref{thm:lasso}}
Without loss of generality, we assume $p$ is even. Then, for block-wise diagonal designs as in \eqref{block}, the Lasso objective is separable. Therefore, for each $W_j^*$, it is not affected by any $\beta_k$ outside the  block. Additionally, by symmetry, the distribution of $W_j^*$ is the same for all $1\leq j\leq p$.   
It follows that 
\begin{align} \label{thm-lasso-FP}
\FP_p(u) &= L_pp\cdot \mathbb{P}\bigl\{W^*_j>t_p(u) \,\big|\, (\beta_j,\beta_{j+1})=(0,0)\bigr\}\cr
&\qquad + L_pp^{1-\vartheta}\cdot \mathbb{P}\bigl\{W^*_j>t_p(u) \,\big|\, (\beta_j,\beta_{j+1})=(0,\tau_p)\bigr\},
\end{align}
where $j$ can be odd index. 
Similarly, we can derive that
\begin{align} \label{thm-lasso-FN}
\FN_p(u) &= L_pp^{1-\vartheta}\cdot \mathbb{P}\bigl\{W^*_j<t_p(u) \,\big|\, (\beta_j,\beta_{j+1})=(\tau_p,0)\bigr\}\cr
&\qquad + L_pp^{1-2\vartheta}\cdot \mathbb{P}\bigl\{W^*_j<t_p(u) \,\big|\, (\beta_j,\beta_{j+1})=(\tau_p,\tau_p)\bigr\}.
\end{align}
Fix variables $\{j, j+1\}$, and consider the random vector $\hat{h}=(x_j'y, x_{j+1}'y)'/\sqrt{\log(p)}$. Then,
\[
\hat{h}\sim {\cal N}_2\Bigl(\mu, \; \frac{1}{\log(p)}\Sigma\Bigr), \qquad \mbox{where}\quad \Sigma=\matB. 
\]
The vector $\mu$ is equal to 
\beq \label{thm-lasso-mu}
\mu^{(1)}\equiv \begin{bmatrix}0\\0\end{bmatrix}, \quad \mu^{(2)}\equiv 
\begin{bmatrix}\rho\sqrt{r}\\ \sqrt{r}\end{bmatrix}
, \quad \mu^{(3)}\equiv \begin{bmatrix}\sqrt{r}\\ \rho\sqrt{r}\end{bmatrix}, \quad \mu^{(4)}\equiv \begin{bmatrix}(1+\rho)\sqrt{r}\\ (1+\rho)\sqrt{r})\end{bmatrix}, 
\eeq
in the four cases where $(\beta_j,\beta_{j+1})'$ is $(0,0)'$, $(0, \tau_p)'$, $(\tau_p, 0)'$, and $(\tau_p, \tau_p)'$, respectively. Let ${\cal R}_u$ be the rejection region induced by Lasso-path, given explicitly in Lemma~\ref{lem:block-RejectRegion}. By  Lemma~\ref{lem:tool}, the probabilities in \eqref{thm-lasso-FP} and \eqref{thm-lasso-FN} are related to the following quantities:
\[
\alpha_k=\begin{cases}
\inf \inf_{h\in {\cal R}_u}\{(h-\mu^{(k)})'\Sigma^{-1}(h-\mu^{(k)})\},& k=1,2,\cr
\inf_{h\in {\cal R}^c_u}\{(h-\mu^{(k)})'\Sigma^{-1}(h-\mu^{(k)})\}, & k=3,4. 
\end{cases}
\]
 and plug it into \eqref{thm-lasso-FP} and \eqref{thm-lasso-FN}. It gives
\beq \label{thm-lasso-FP+FN}
\FP_p(u) =L_p p^{ 1-\min\{\alpha_1,\; \vartheta + \alpha_2\}},\quad
\FN_p(u) =L_p p^{ 1-\min\{\vartheta+ \alpha_3,\; 2\vartheta+ \alpha_4\}}. 
\eeq
It remains to compute the exponents $\alpha_1$-$\alpha_4$. 

First, we consider the case that $\rho\geq 0$. The rejection region in Figure~\ref{fig:region-proof} is defined by the following lines:
\begin{itemize}
\item Line 1: $h_1-\rho h_2=(1-\rho)\sqrt{u}$. 
\item Line 2: $h_1=\sqrt{u}$. 
\item Line 3: $h_1-\rho h_2 = (1+\rho)\sqrt{u}$.
\item Line 4: $h_1-\rho h_2=-(1-\rho)\sqrt{u}$.
\item Line 5: $h_1=-\sqrt{u}$.
\item Line 6: $h_1-\rho h_2=-(1+\rho)\sqrt{u}$. 
\end{itemize}
Consider a general ellipsoid:
\[
{\cal E}(t; \mu)= \{h\in\mathbb{R}^2: (h-\mu)'\Sigma^{-1}(h-\mu)'\leq t\}. 
\]
Given any line $h_1+bh_2=c$, 
as $t$ increases, this ellipsoid eventually intersects with this line. The intersection point is computed by the following equations:
\[
h_1+bh_2=c, \qquad \Sigma^{-1}(h-\mu)\propto (1, b)'. 
\] 
The second equation (it is indeed a linear equation on $h$) says that the perpendicular vector of the tangent plane is orthogonal to the line. Solving the above equations gives the intersection point and the value of $t$: As long as $b^2\neq 1$, we have
\beq \label{thm-lasso-intersection}
h^*=\mu+ \frac{c-(\mu_1+b\mu_2)}{1+b^2+2b\rho}
\begin{bmatrix}
1+b\rho\\
b+\rho
\end{bmatrix}, \qquad t^* = \frac{[c-(\mu_1+b\mu_2)]^2}{1+b^2+2b\rho}. 
\eeq
Using the expressions of lines 1-6, we can obtain the corresponding $t^*$ for 6 lines: 
\begin{align*}
&t^*_1= \frac{[(1-\rho)\sqrt{u}-(\mu_1-\rho\mu_2)]^2}{1-\rho^2}, \qquad t^*_2= (\sqrt{u}-\mu_1)^2, \qquad t_3^*=\frac{[(1+\rho)\sqrt{u}-(\mu_1-\rho\mu_2)]^2}{1-\rho^2},\cr
&t_4^*= \frac{[(1-\rho)\sqrt{u}+(\mu_1-\rho\mu_2)]^2}{1-\rho^2}, \qquad t_5^*=(\sqrt{u}+\mu_1)^2, \qquad t_6^*=\frac{[(1+\rho)\sqrt{u}+(\mu_1-\rho\mu_2)]^2}{1-\rho^2}. 
\end{align*}
We first look at the ellipsoid ${\cal E}(t;\mu^{(1)})$ and study when it intersects with ${\cal R}_u$. Note that $\mu^{(1)}=(0, 0)'$. The above $t^*$ values become
\[
t_2^*=t_5^*=u, \qquad t_1^*=t_4^*=\frac{u}{1+\rho}, \qquad t_3^*=t_6^*=\frac{u}{1-\rho}. 
\]
Therefore, as we increase $t$, this ellipsoid first intersects with line 1 and line 4. For line 1, the intersection point is $((1-\rho)\sqrt{u}, 0)'$, but it is outside the rejection region (see Figure~\ref{fig:region-proof}); the situation for line 4 is similar. We then further increase $t$, and the ellipsoid intersects with line 2 and line 5, where the intersection point is $(\sqrt{u}, \rho\sqrt{u})'$; this point is indeed on the boundary of the rejection region. We thus conclude that 
\beq \label{thm-lasso-exponent1}
\inf_{h\in {\cal R}_u}\{(h-\mu^{(1)})'\Sigma^{-1}(h-\mu^{(1)})\}=u.
\eeq
We then look at the the ellipsoid ${\cal E}(t;\mu^{(2)})$, with $\mu^{(2)}=(\rho\sqrt{r}, \sqrt{r})'$. The $t^*$ values for 6 lines are:
\[
t_1^*=t_4^*=\frac{1-\rho}{1+\rho}u, \qquad t_2^*=(\sqrt{u}-\rho\sqrt{r})^2, \qquad t_3^*=t_6^*=\frac{1+\rho}{1-\rho}u, \qquad t_5^*=(\sqrt{u}+\rho\sqrt{r})^2. 
\]
The smallest $t^*$ is among $\{t_1^*, t_2^*, t_4^*\}$. 
Since $\mu^{(2)}$ is in the positive orthant, the intersection point of the ellipsoid with line 4 must be outside the rejection region, so we further restrict to $t_1^*$ and $t_2^*$. The ellipsoid intersects with line 1 at $(\rho\sqrt{r}+(1-\rho)\sqrt{u},\; \sqrt{r})'$. This point is on the boundary of ${\cal R}_u$ if and only if its second coordinate is $\geq\sqrt{u}$ (see Figure~\ref{fig:region-proof}), i.e., $u\leq r$. 
The ellipsoid intersects with line 2 at $(\sqrt{u},\; \rho\sqrt{u}+(1-\rho^2)\sqrt{r})'$. This point is on the boundary of ${\cal R}_u$ if and only if its second coordinate is $\leq \sqrt{u}$ (see Figure~\ref{fig:region-proof}), i.e., $u\geq (1+\rho)^2r$. In the range of $r<u<(1+\rho)^2 r$, the ellipsoid intersects with ${\cal R}_u$ at the corner point $(\sqrt{u}, \sqrt{u})'$, with the corresponding 
\[
t^*=r+\frac{2}{1+\rho}u-2\sqrt{ru}=\begin{cases}
\frac{1-\rho}{1+\rho}u+(\sqrt{u}-\sqrt{r})^2,\cr
(\sqrt{u}-\rho\sqrt{r})^2+\frac{1-\rho}{1+\rho}\bigl(\sqrt{u}-(1+\rho)\sqrt{r}\bigr)^2.\end{cases} 
\]
This $t^*$ has two equivalent expressions. 
Comparing them with $t_1^*$ and $t_2^*$, we can see that the smallest $t^*$ is a continuous function of $u$, given $(\rho, r)$. 
It follows that 
\begin{align} \label{thm-lasso-exponent2}
& \inf_{h\in {\cal R}_u}\{(h-\mu^{(2)})'\Sigma^{-1}(h-\mu^{(2)})\}\cr
=\;\; &\frac{1-\rho}{1+\rho}u+(\sqrt{u}-\sqrt{r})_+^2 - \frac{1-\rho}{1+\rho}\bigl(\sqrt{u}-(1+\rho)\sqrt{r}\bigr)_+^2.
\end{align}
We plug \eqref{thm-lasso-exponent1} and \eqref{thm-lasso-exponent2} into \eqref{thm-lasso-FP+FN}. It gives the expression of $\FP_p(u)$ for $\rho\geq 0$. 

We then look at the ellipsoid ${\cal E}(t;\mu^{(3)})$, with $\mu^{(3)}=(\sqrt{r}, \rho\sqrt{r})'$. Note that we now investigate its distance to the complement of ${\cal R}_u$. In order for $\mu^{(3)}$ to outside ${\cal R}_u^c$ (i.e., in the interior of ${\cal R}_u)$, we require that $u<r$; furthermore, when $u<r$, the ellipsoid can only intersect with lines 1-2 (see Figure~\ref{fig:region-proof}). Using the formula of $t^*$ in the equation below \eqref{thm-lasso-intersection}, we have
\[
t_1^*=\frac{1-\rho}{1+\rho}\bigl((1+\rho)\sqrt{r}-\sqrt{u}\bigr)^2, \qquad t_2^*=(\sqrt{r}-\sqrt{u})^2. 
\]
By \eqref{thm-lasso-intersection}, the ellipsoid intersects with line 1 at $\bigl(\sqrt{r}-(1-\rho)[(1+\rho)\sqrt{r}-\sqrt{u}],\; \rho\sqrt{r}\bigr)'$. To guarantee that this point is on the boundary of ${\cal R}_u$, we need its second coordinate to be $\geq \sqrt{u}$ (see Figure~\ref{fig:region-proof}), i.e., $u\leq \rho^2 r$; furthermore, when $u>\rho^2r$, it can be easily seen from Figure~\ref{fig:region-proof} that the ellipsoid must have already crossed line 2. By \eqref{thm-lasso-intersection} again, the ellipsoid intersects with line 2 at $(\sqrt{u}, \rho \sqrt{u})'$. This point is always on the boundary of ${\cal R}_u$. It follows that
\beq \label{thm-lasso-exponent3}
\inf_{h\in {\cal R}^c_u}\{(h-\mu^{(3)})'\Sigma^{-1}(h-\mu^{(3)})\}=\min\Bigl\{ \frac{1-\rho}{1+\rho}\bigl((1+\rho)\sqrt{r}-\sqrt{u}\bigr)^2,\;\, (\sqrt{r}-\sqrt{u})_+^2 \Bigr\}.
\eeq
We then look at the ellipsoid ${\cal E}(t;\mu^{(4)})$, with $\mu^{(4)}=\bigl((1+\rho)\sqrt{r}, (1+\rho)\sqrt{r}\bigr)'$. It follows from figure~\ref{fig:region-proof} that $\mu^{(4)}$ is in the interior of the ellipsoid if and only if $(1+\rho)\sqrt{r}>\sqrt{u}$. We restrict to $(1+\rho)\sqrt{r}>\sqrt{u}$. Then, this ellipsoid can only touch lines 1-2 first. The $t^*$ values are
\[
t_1^* = \frac{1-\rho}{1+\rho}\bigl((1+\rho)\sqrt{r}-\sqrt{u}\bigr)^2, \qquad t_2^*=\bigl((1+\rho)\sqrt{r}-\sqrt{u}\bigr)^2. 
\]
Since $t_1^*<t_2^*$, the ellipsoid touches line 1 first, at the intersection point $\bigl((1-\rho)\sqrt{u} +\rho (1+\rho)\sqrt{r},\; (1+\rho)\sqrt{r}\bigr)'$. In order for this point to be on the boundary of ${\cal R}_u$, we need that its second coordinate is $\geq\sqrt{u}$, which translates to $\sqrt{u}\leq (1+\rho)\sqrt{r}$. This is always true when $r>u$ and $\rho>0$. It follows that 
\beq \label{thm-lasso-exponent4}
\inf_{h\in {\cal R}^c_u}\{(h-\mu^{(4)})'\Sigma^{-1}(h-\mu^{(4)})\}= \frac{1-\rho}{1+\rho}\bigl((1+\rho)\sqrt{r}-\sqrt{u}\bigr)_+^2.
\eeq
We plug \eqref{thm-lasso-exponent3} and \eqref{thm-lasso-exponent4} into \eqref{thm-lasso-FP+FN}. It gives the expression of $\FN_p(u)$ for $\rho\geq 0$. 

Next, we consider the case that $\rho<0$. By Lemma~\ref{lem:block-RejectRegion}, ${\cal R}_u(\rho)$ is a reflection of ${\cal R}_u(|\rho|)$ with respect to the x-axis. As a result, if we re-define $\hat{h}=(x_j'y,\; - x_{j+1}'y)/\sqrt{2\log(p)}$, then the rejection region becomes ${\cal R}_u(|\rho|)$, which has the same shape as that in Figure~\ref{fig:region-proof}. At the same time, the distribution of $\hat{h}$ becomes
\[
\hat{h}\sim {\cal N}_2\Bigl(\mu, \; \frac{1}{\log(p)}\Sigma\Bigr), \qquad \mbox{where}\quad \Sigma=\matBabs. 
\]
The vector $\mu$ is equal to 
\beq \label{thm-lasso-mu-negative}
\mu^{(1)}\equiv \begin{bmatrix}0\\0\end{bmatrix}, \quad \mu^{(2)}\equiv 
\begin{bmatrix}-|\rho|\sqrt{r}\\ -\sqrt{r}\end{bmatrix}
, \quad \mu^{(3)}\equiv \begin{bmatrix}\sqrt{r}\\ |\rho|\sqrt{r}\end{bmatrix}, \quad \mu^{(4)}\equiv \begin{bmatrix}(1-|\rho|)\sqrt{r}\\ -(1-|\rho|)\sqrt{r})\end{bmatrix}, 
\eeq
when $(\beta_j,\beta_{j+1})'$ is $(0,0)'$, $(0, \tau_p)'$, $(\tau_p, 0)'$, and $(\tau_p, \tau_p)'$, respectively. Therefore, the calculations are similar, except that the expressions of $\mu^{(1)}$ to $\mu^{(4)}$ have changed to \eqref{thm-lasso-mu-negative}. 

Below, for a negative $\rho$, we calculate the exponents in \eqref{thm-lasso-FP+FN} as follows: We pretend that $\rho>0$ and calculate the exponents using the same ${\cal R}_u$ and $\Sigma$ as before, with $\mu^{(1)}$ to $\mu^{(4)}$ replaced by those in \eqref{thm-lasso-mu-negative}. Finally, we replace $\rho$ by $|\rho|$ in all four exponents.

We now pretend that $\rho>0$. Then, for each ellipsoid ${\cal E}(t;\mu^{(k)})$, its intersection point with a line $h_1+bh_2=c$ still obeys the formula in \eqref{thm-lasso-intersection}, and the corresponding $t_*$ values associated with line 1-line 6 are still the same as those in the equation below \eqref{thm-lasso-intersection} (but the vector $\mu$ has changed). Comparing \eqref{thm-lasso-mu-negative} with \eqref{thm-lasso-mu}, we notice that $\mu^{(1)}$ and $\mu^{(3)}$ are unchanged. Therefore, the expressions of exponents in \eqref{thm-lasso-exponent1} and \eqref{thm-lasso-exponent3} are still correct. The current $\mu^{(2)}$ is a sign flip (on both x-axis and y-axis) of the $\mu^{(2)}$ in \eqref{thm-lasso-mu}; also, it can be seen from Figure~\ref{fig:region-proof} that the rejection region remains unchanged subject to a sign flip. Therefore, the expression in \eqref{thm-lasso-exponent3} is also valid. We only need to re-calculate the exponent in \eqref{thm-lasso-exponent4}. The current $\mu^{(4)}$ is in the 4-th orthant. It is in the interior of ${\cal R}_u$ only if $(1-\rho)\sqrt{r}>\sqrt{u}$, i.e., $u<(1-\rho)^2 r$. As we increase $t$, the ellipsoid ${\cal E}(t;\mu^{(4)})$ will first intersect with either line 2 or line 3. Using the formula of $t^*$ in the equation below \eqref{thm-lasso-intersection}, we have
\[
t_2^* = \bigl(\sqrt{u}-(1-\rho)\sqrt{r}\bigr)^2, \qquad t_3^* = \frac{1+\rho}{1-\rho} \bigl((1-\rho)\sqrt{r}-\sqrt{u}\bigr)^2. 
\]
While $t_2^*$ is the smaller one, the intersection point of the ellipsoid with line 2 is $(\sqrt{u}, -(1-\rho)\sqrt{r})'$, which by Figure~\ref{fig:region-proof} is in the interior of ${\cal R}_u$. Hence, the ellipsoid hits line 3 first. We conclude that
\beq \label{thm-lasso-exponent4-negative}
\inf_{h\in {\cal R}^c_u}\{(h-\mu^{(4)})'\Sigma^{-1}(h-\mu^{(4)})\}= \frac{1+\rho}{1-\rho}\bigl((1-\rho)\sqrt{r}-\sqrt{u}\bigr)_+^2.
\eeq
Finally, we plug \eqref{thm-lasso-exponent1}, \eqref{thm-lasso-exponent2}, \eqref{thm-lasso-exponent3} and \eqref{thm-lasso-exponent4-negative} into \eqref{thm-lasso-FP+FN}, and then change $\rho$ to $|\rho|$. This gives the expressions of $\FP_p(u)$ and $\FN_p(u)$ for a negative $\rho$.

\section{Proof of Theorem~\ref{thm:knockoff-block}}

We assume $\rho\geq 1/2$ throughout the proof. The calculation for the case where $\rho\leq -1/2$ is similar. By the design of the gram matrix $X^TX$ and the construction of the knockoff variables, we know Lasso regression problem with $2p$ variables can be reduced to $(p/2)$ independent four-variate Lasso regression problems:
\begin{equation}\label{Lasso3}
    (\hat{\beta}_j,\hat{\beta}_{j+1},\hat{\beta}_{j+p},\hat{\beta}_{j+p+1})(\lambda)=\mathrm{argmin}_b\Big\{\frac{1}{2}||y-(x_j,x_{j+1},\tilde{x}_j,\tilde{x}_{j+1})b||_2^2+\lambda ||b||_1\Big\}
\end{equation}
for $j=1,3,\cdots,p-1$. By taking the sub-gradients of the objective function in (\ref{Lasso3}), we know $(\hat{\beta}_j,\hat{\beta}_{j+1},\hat{\beta}_{j+p},\hat{\beta}_{j+p+1})$ should satisfy:
\begin{equation}\label{Lasso4}
\begin{split}
      (\hat{\beta}_j,\hat{\beta}_{j+1},\hat{\beta}_{j+p},\hat{\beta}_{j+p+1})G +\lambda (\text{sgn}(\hat{\beta}_j),\text{sgn}(\hat{\beta}_{j+1}),\text{sgn}(\hat{\beta}_{j+p}),\text{sgn}(\hat{\beta}_{j+p+1}))\\
      =(y^T x_j,y^T x_{j+1},y^T \tilde{x}_{j},y^T \tilde{x}_{j+1})  
\end{split}
\end{equation}
where $G=((1,\rho,2\rho-1,\rho)^T,(\rho,1,\rho,2\rho-1)^T,(2\rho-1,\rho,1,\rho)^T,(\rho,2\rho-1,\rho,1)^T)$ and sgn$(x)=1$ if $x>0$; $-1$ if $x<0$; any value in $[-1,1]$ if $x=0$. We have choose the correlation between a true variable and its knockoff to be $2\rho-1$, which is the smallest value such that $(X,\tilde{X})^T(X,\tilde{X})$ is semi-positive definite. In this case, $G$ is degenerated and has rank 3. As $\lambda$ is decreasing from infinity,
 we recognize that the first two variables (assume these two features are linear independent) entering the model will not leave before the third variable enters the model, which is obviously true from the close form solution of the bi-variate Lasso problem. We then show that the first two variables enter the Lasso path, individually. Furthermore, if the first two variables are a true variable and its knockoff variable, then the third and fourth variable enter the Lasso path simultaneously. 
 
Since $(y^T x_j,y^T x_{j+1},y^T \tilde{x}_{j},y^T \tilde{x}_{j+1})^T\sim \mathcal{N}(G(\beta_j,\beta_{j+1},0,0)^T ,G)$ is a degenerated normal random variable, we reparametrize it as $(m+d_1,m+d_2,m-d_1,m-d_2)$ with $(m,d_1,d_2)^T\sim \mathcal{N}((\rho\beta_j+\rho\beta_{j+1},(1-\rho)\beta_j,(1-\rho)\beta_{j+1})^T,\text{diag}(\rho,1-\rho,1-\rho))$. We intend to give the Lasso solution path (or $Z_j, \tilde{Z}_j$) as a function of $m,d_1$ and $d_2$. We only present the result in the case where $d_1>d_2>0$. Results from other cases are immediate by permuting the rows in equation set~(\ref{Lasso4}) and transforming to the $d_1>d_2>0$ case. Lasso solution path are obtained by the KKT condition (\ref{Lasso4}) and summarized in the table below.
\begin{table}[!htb]
\centering
\def~{\hphantom{0}}
\scalebox{0.85}{
\begin{tabular}{c|cccccc}
 range of $m$ &  $\lambda_1$  & $\text{sign}_1$ & $\lambda_2$  & $\text{sign}_2$ & $\lambda_3$  & $\text{sign}_3$   \\ \hline
$(-\infty,\frac{\rho}{1-\rho}(d_2-d_1))$ & $-m+d_1$  & $(0,0,0^-,0)$  & $-m-\frac{\rho}{1-\rho}d_1+\frac{1}{1-\rho}d_2$ & $(0,0,-,0^-)$ &  &   \\
$(\frac{\rho}{1-\rho}(d_2-d_1),0)$ & $-m+d_1$ & $(0,0,0^-,0)$ & $\frac{1-\rho}{\rho}m+d_1$ & $(0^+,0,-,0)$ & $d_2$ &  $(+,0^+,-,0^-)$ \\ 

$(0,\frac{\rho}{1-\rho}(d_1-d_2))$ & $m+d_1$ & $(0^+,0,0,0)$ & $\frac{\rho-1}{\rho}m+d_1$ & $(+,0,0^-,0)$ & $d_2$ &  $(+,0^+,-,0^-)$ \\ 

$(\frac{\rho}{1-\rho}(d_1-d_2),\infty)$ & $m+d_1$  & $(0^+,0,0,0)$ & $m-\frac{\rho}{1-\rho}d_1+\frac{1}{1-\rho}d_2$ & $(+,0^+,0,0)$ &  &  \\ \hline
\end{tabular}}
\label{tb:lassopath2}
\caption{Summary of solution path of the Lasso problem (\ref{Lasso3}). $\lambda_i$ record the critical value of $\lambda$ where a new variable enters the model and $\text{sign}_i$ records the sign and the limiting behavior of $(\hat{\beta}_j,\hat{\beta}_{j+p})$ as $\lambda \to \lambda_i^-$. Value of $\lambda_3$ is omitted in row 1 and 4 since it will not affect the value of $W_j$ and $W_{j+1}$.}
\end{table}

Here we explain the third row of the table as an example, $b_1=(\epsilon,0,0,0)^T$ is a solution of the KKT condition (\ref{Lasso4}) when $\lambda=m+d_1-\epsilon$ for $\epsilon\in (0,\frac{m}{\rho}]$, so $\text{sign}_1$ is expressed as $(0^+,0,0,0)$. By property of the Lasso solution, if $b_1$ and $b_2$ are both Lasso solutions, then $G(b_1-b_2)=0$ and $||b_1||_1=||b_2||_1$. $G(b_1-b_2)=0$ implies $b_1-b_2=\delta \times (1,-1,1,-1)^T$ for some $\delta\neq 0$. Therefore, $b_2=(\epsilon-\delta,\delta,-\delta,\delta)^T$ and $||b_2||_1\geq ||b_1||_1+2|\delta|$. This means the Lasso solution is unique with $\lambda=m+d_1-\epsilon$ and variable 1 is the only one entering the model when $\lambda$ gets below $\lambda_1$. When $\lambda=\frac{\rho-1}{\rho}m+d_1-\epsilon$ for $\epsilon\in(0,\frac{\rho-1}{\rho}m+d_1-d_2]$, $b_1=(\frac{m}{\rho}+\frac{\epsilon}{2-2\rho},0,-\frac{\epsilon}{2-2\rho},0)^T$ is a solution of the KKT conditions. If there is another Lasso solution $b_2$, then $b_2=(\frac{m}{\rho}+\frac{\epsilon}{2-2\rho}-\delta,\delta,-\frac{\epsilon}{2-2\rho}-\delta,\delta)^T$ and $||b_2||_1\geq ||b_1||_1+2|\delta|$. So $b_2$ does not exist and variable 3 is the only one entering the model when $\lambda$ gets below $\lambda_2$. When $\lambda=d_2-\epsilon$ for sufficient small positive $\epsilon$, $b_1=(\frac{m}{2\rho}+\frac{d_1}{2-2\rho},\frac{\epsilon}{2-2\rho},\frac{m}{2\rho}-\frac{d_1}{2-2\rho},-\frac{\epsilon}{2-2\rho})^T$ satisfies the KKT condition, thus variable 2 and 4 enters the model simultaneously. At this point, the Lasso solution is not unique and all solutions can be expressed as $b_1-\delta \times (1,-1,1,-1)^T$ with $\delta\in [-\frac{\epsilon}{2-2\rho},\frac{\epsilon}{2-2\rho}]$. Other rows from the table can be analyzed similarly.

Table~\ref{tb:lassopath2} implicitly expresses $Z_j, Z_{j+1},\tilde{Z}_j$ and $\tilde{Z}_{j+1}$ as a function of $d_1,d_2$ and $m$. By examining all possible ordinal relationship of $d_1,d_2$ and 0, we record the region in the space of $(d_1,d_2,m)$ such that $\hat{\beta}_j(u)>0$ and denote it as $R(u)$. $R(u)$ is the union of 4 disjoint sub-regions $\{R_i(u)\}_{i=1,\cdots,4}$, defined as following:
\begin{equation}\label{R}
\begin{split}
    R_1(u)=&\{(x,y,z):x>0,y>0,x>y,z>0,x+z>T\}\\
    & \cup\frac{1}{2}\{(x,y,z):x>0,y>0, x<y, z<0,z>x-y, x>T\}\\
    & \cup\frac{1}{2}\{(x,y,z):x>0, y>0, x<y, z>0, z<\frac{\rho}{1-\rho}(y-x), x>T\}\\
    & \cup \{(x,y,z):x>0, y>0, x<y, z>0,z>\max\big(\frac{\rho}{1-\rho}(y-x),T+\frac{\rho}{1-\rho}y-\frac{1}{1-\rho}x\big)\},
\end{split}
\end{equation}
$R_2(u)=\{(x,y,z):(-x,y,-z)\in R_1(u)\}$, $R_3(u)=\{(x,y,z):(x,-y,z)\in R_1(u)\}$ and $R_4(u)=\{(x,y,z):(-x,-y,-z)\in R_1(u)\}$,
where $T=\sqrt{2u\log(p)}$ and the $\frac{1}{2}$ ahead of a certain region means when $(d_1,d_2,m)$ is in this region, $\hat{\beta}_j(u)>0$ happens with $1/2$ probability. Let the four disjoint regions that composes $R_1(u)$ in (\ref{R}) be denoted by $R_{1,j}(u)$ for $j=1,\cdots,4$. We can similarly define $R_{i,j}(u)$ for $i=2,3,4$. By Lemma 1, as $p \to \infty$,
\begin{equation}\label{KnockoffFP}
\begin{split}
      \mathbb{P}(\beta_j=0, \hat{\beta}_j(u)\neq 0)
      =&\mathbb{P}(\hat{\beta}_j(u)\neq 0|\beta_j=0,\beta_{j+1}=0)\times \mathbb{P}(\beta_j=0,\beta_{j+1}=0)\\
      &+\mathbb{P}(\hat{\beta}_j(u)\neq 0|\beta_j=0,\beta_{j+1}=\tau_p)\times \mathbb{P}(\beta_j=0,\beta_{j+1}=\tau_p)\\
      =&L_p p^{-\inf_{R(u)} [(z^2/\rho+x^2/(1-\rho)+y^2/(1-\rho))/(2\log(p))]}\\
      &+L_p p^{-\vartheta-\inf_{R(u)} [((z-\rho\tau_p)^2/\rho+x^2/(1-\rho)+(y-(1-\rho)\tau_p)^2/(1-\rho))/(2\log(p))]},
\end{split}
\end{equation}
\begin{equation}\label{KnockoffFN}
\begin{split}
      \mathbb{P}(\beta_j\neq 0, \hat{\beta}_j(u)=0)
      =&\mathbb{P}(\hat{\beta}_j(u)= 0|\beta_j=\tau_p,\beta_{j+1}=0)\times \mathbb{P}(\beta_j=\tau_p,\beta_{j+1}=0)\\
      &+\mathbb{P}(\hat{\beta}_j(u)=0|\beta_j=\tau_p,\beta_{j+1}=\tau_p)\times \mathbb{P}(\beta_j=\tau_p,\beta_{j+1}=\tau_p)\\
      =&L_p p^{-\vartheta-\inf_{R(u)^C} [((z-\rho\tau_p)^2/\rho+(x-(1-\rho)\tau_p)^2/(1-\rho)+y^2/(1-\rho))/(2\log(p))]}\\
      &+L_p p^{-2\vartheta-\inf_{R(u)^C} [((z-2\rho\tau_p)^2/\rho+(x-(1-\rho)\tau_p)^2/(1-\rho)+(y-(1-\rho)\tau_p)^2/(1-\rho))/(2\log(p))]}.
\end{split}
\end{equation}
Define the $\rho$-distance function of two sets $A$ and $B$ in $\mathbb{R}^3$ as 
\[d_{\rho}(A,B)=\inf_{a\in A, b\in B}[(a_1-b_1)^2/(1-\rho)+(a_2-b_2)^2/(1-\rho)+(a_3-b_3)^2/\rho]\]
where $a_k, b_k$ denote the $k$-th coordinate of vector $a$ and $b$. An immediate property of the $\rho$-distance function would be
\[d_{\rho}(\cup_{i=1,\cdots,M}A_i,\cup_{j=1,\cdots,N}B_j)=\min_{i,j}d_{\rho}(A_i,B_j).\]

Utilizing the symmetry of the regions, we can compute the region distances involved in (\ref{KnockoffFP}) and (\ref{KnockoffFN}) explicitly. Take the second exponent in (\ref{KnockoffFP}) as an example, it can be simplified as
\[
\begin{split}
-\vartheta-&d_{\rho}(R(u),\{(0,(1-\rho)\tau_p,\rho\tau_p)\})/(2\log(p))\\
=-\vartheta-&d_{\rho}(R_1(u)\cup R_2(u)\cup R_3(u)\cup R_4(u),\{(0,(1-\rho)\tau_p,\rho\tau_p)\})/(2\log(p))\\
=-\vartheta-&d_{\rho}(R_{1,1}(u)\cup R_{1,3}(u)\cup R_{1,4}(u)\cup R_{2,2}(u),\{(0,(1-\rho)\tau_p,\rho\tau_p)\} ))/(2\log(p)).
\end{split}
\]
Define $\tilde{R}_{1,2}(u)=\{(x,y,z):x>0,y>0,z>0, x<y, x>T, z<y-x\}$, $\tilde{R}_{1,3}(u)=\{(x,y,z):x>0,y>0,z>0, x<y, x>T\}$ and $\tilde{R}_{1,4}(u)=\{(x,y,z):x>0,y>0,z>0, x<y, x<T, z>T+\frac{\rho}{1-\rho}y-\frac{1}{1-\rho}x\}$. Then $\tilde{R}_{1,2}(u)\subset \tilde{R}_{1,3}(u)$ and $R_{1,3}(u)\cup R_{1,4}(u)=\tilde{R}_{1,3}(u)\cup \tilde{R}_{1,4}(u)$. Since $\tilde{R}_{1,2}$(u) and $R_{2,2}(u)$ are symmetric about the plane $x=0$, we know
\[d_{\rho}(R_{2,2}(u),\{(0,(1-\rho)\tau_p,\rho\tau_p)\})=d_{\rho}(\tilde{R}_{1,2}(u),\{(0,(1-\rho)\tau_p,\rho\tau_p)\}).\]

Therefore, 
\[
\begin{split}
d_{\rho}(R&(u),\{(0,(1-\rho)\tau_p,\rho\tau_p)\})\\
=\min\{&d_{\rho}(R_{1,1}(u),\{(0,(1-\rho)\tau_p,\rho\tau_p)\}),d_{\rho}(\tilde{R}_{1,3}(u),\{(0,(1-\rho)\tau_p,\rho\tau_p)\}),\\
&d_{\rho}(\tilde{R}_{1,4}(u),\{(0,(1-\rho)\tau_p,\rho\tau_p)\})\}\\
=\min\Big\{&\frac{1-\rho}{2}\times \tau_p^2+\frac{2}{1+\rho}\times [(T-(1+\rho)\tau_p/2)_+]^2-\frac{1-\rho}{1+\rho}\times[(T-(1+\rho)\tau_p)_+]^2,\\
&\frac{1}{1-\rho}\times T^2+\frac{1}{1-\rho}\times[(T-(1-\rho)\tau_p)_+]^2,\ \ \frac{1-\rho}{1+\rho}\times T^2+((T-\tau_p)_+)^2\Big\}\\
=(T-&\rho\tau_p)^2+(\xi_\rho \tau_p-\eta_\rho T)_+^2-(\tau_p-T)_+^2,
\end{split}
\]
where $\xi_\rho=\sqrt{1-\rho^2}$ and $\eta_\rho=\sqrt{(1-\rho)/(1+\rho)}$.

Let $\tau_p=0$, we know $d_{\rho}(R(u),\{(0,0,0)\})=T^2$.
By (\ref{KnockoffFP}) we immediately have
\begin{equation}
\mathbb{P}(\beta_j=0, \hat{\beta}_j(u)\neq 0)
      =L_p p^{-\min\left\{u,\;\, \vartheta+(\sqrt{u}-\rho\sqrt{r})^2+(\xi_\rho\sqrt{r}-\eta_\rho \sqrt{u})_+^2-(\sqrt{r}-\sqrt{u})_+^2\right\}}.\end{equation}
We can see the false positive rate is exactly the same when using the Lasso filter and the Knockoff filter when $\rho>0$. For $\rho\geq 1/2$, we can similarly compute $d_{\rho}(R(u)^C,\{((1-\rho)\tau_p,0,\rho\tau_p)\})$ to be
\[
[(\tau_p-T)_+-((1-\xi_\rho )\tau_p-(1-\eta_\rho )T)_+-(\lambda_\rho\tau_p-\eta_\rho\ T)_+]^2,
\]
and $d_{\rho}(R(u)^C,\{((1-\rho)\tau_p,(1-\rho)\tau_p,2\rho\tau_p)\})$ to be
\[
[(\xi_\rho\ \tau_p-\eta_\rho T)_+-(\lambda_\rho \tau_p-\eta_\rho T)_+]^2,
\]
where $\xi_\rho=\sqrt{1-\rho^2}$, $\eta_\rho=\sqrt{(1-\rho)/(1+\rho)}$, and $\lambda_\rho=\sqrt{1-\rho^2}-\sqrt{1-\rho}$. 

Plug these results in to (\ref{KnockoffFN}), we have
\begin{equation}
      \mathbb{P}(\beta_j\neq 0, \hat{\beta}_j(u)=0)=L_p p^{- \vartheta-\left\{ (\sqrt{r}-\sqrt{u})_+ - [(1-\xi_\rho)\sqrt{r}-(1-\eta_\rho)\sqrt{u}]_+ - (\lambda_\rho\sqrt{r}-\eta_\rho \sqrt{u})_+\right\}^2 }.
\end{equation}
From here we have prove the result for $\rho\geq 1/2$ case.

In the case where $\rho\leq -1/2$, the exponent of false negative rate is additionally lower bounded by $-2\vartheta$. One can verify the rate given in the theorem through similar calculations. This is somehow more straight forwards since in the case where $\beta_j=\beta_{j+1}=\tau$, $(y^T x_j,y^T x_{j+1},y^T \tilde{x}_{j},y^T \tilde{x}_{j+1})^T\sim \mathcal{N}((1+\rho)\tau\cdot (1,1,-1,-1)^T,G)$, meaning there is no way to distinguish the true variable from its knockoff variable.  

\section{Proof of Theorem~\ref{thm:knockoff-block2}}
In the following proofs, we only consider $\rho \geq 0$ case, since $\rho<0$ case can be transformed to the positive $|\rho|$ case by flipping the sign of either $\beta_j$ or $\beta_{j+1}$ for $j=1,3,\cdots,p-1$. By the block diagonal structure of the gram matrix, the Lasso problem with $2p$ features can be reduced to $(p/2)$ independent four-variate Lasso regression problems:
\begin{equation}\label{Lasso-EC}
    \hat{b}(\lambda)=\mathrm{argmin}_b\Big\{\frac{1}{2}||y-(x_j,x_{j+1},\tilde{x}_j,\tilde{x}_{j+1})b||_2^2+\lambda ||b||_1\Big\}
\end{equation}
for $j=1,3,\cdots,p-1$.
Before we turn to the proof of the theorem, we first analysis the solution path of the following four-variate Lasso problem:
\beq\label{4-Lasso}
\hat{b} = \mathrm{argmin}_b \bigl\{ - h^Tb + b^TBb/2 + \lambda\|b\|_1\bigr\}. 
\eeq
with $B=((1,\rho,a,\rho)^T,(\rho,1,\rho,a)^T,(a,\rho,1,\rho)^T,(\rho,a,\rho,1)^T)$ and $a\in[2|\rho|-1,1]$.
By taking the sub-gradients, we know $\hat{b}$ should satisfy
\begin{equation}\label{4-Lasso-KKT}
     B\ \hat{b} +\lambda \ \text{sgn}(\hat{b}) =h.
\end{equation}
Let $\hat{b}_i$ and $h_i$ denotes the $i$-th coordinate of $\hat{b}$ and $h$. Let $\lambda_1 > \lambda_2> \lambda_3 > \lambda_4$ be the values at which variables enter the solution path. 
As discussed in the proof of Lemma~\ref{lem:block-RejectRegion}, $\lambda_1=\max\{|h_1|,|h_2|,|h_3|,|h_4|\}$. Without loss of generality, assume $\lambda_1=|h_1|$ and variable 1 is the first variable entering the model in solution path. We know for one variate Lasso problem, the only feature will not leave the model after its entry as $\lambda$ is decreasing. So in the four-variate Lasso (\ref{4-Lasso}), variable 1 will stay in the model until the second variable enters the model. Consider three bi-variate Lasso problems ($k=2,3,4$):
\beq\label{bi-Lassos}
\hat{b}^{(k)} = \mathrm{argmin}_{b^{(k)}}  \bigl\{ - (h^{(k)}) ^T b^{(k)}  + (b^{(k)}) ^TB^{(k)} b^{(k)} /2 + \lambda\|b^{(k)} \|_1\bigr\}
\eeq
with
\[
B^{(2)}=B^{(4)}= \begin{bmatrix} 1 \;\;\; & \rho\\\rho\;\;\; & 1\end{bmatrix}\qquad \mbox{and}\qquad 
B^{(3)}= \begin{bmatrix} 1 \;\;\; & a\\a\;\;\; & 1\end{bmatrix},
\]
$h^{(2)}=(h_1,h_2)$, $h^{(3)}=(h_1,h_3)$ and $h^{(4)}=(h_1,h_4)$.
Now, we claim $\lambda_2=\max_{k}\{\lambda_2^{(k)}\}$ where $\lambda_2^{(k)}$ is the value at which the second variables enter the solution path in the $k$-th bi-variate Lasso problems. Suppose $\lambda_2^{(i)}> \lambda_2^{(k)}$ for $i\neq k\in \{2,3,4\}$, when $\lambda\in [\lambda_2^{(i)},\lambda_1]$, we know the KKT condition (\ref{4-Lasso-KKT}) is satisfied with $h_2=h_3=h_4=0$ by looking at the KKT conditions of the bi-variate Lasso problems. When $\lambda\in[\lambda_2^{(i)}-\epsilon,\lambda_2^{(i)})$, a second variable $i$ must have entered the four-variate Lasso path, since the objective function of (\ref{4-Lasso}) is smaller when including variable $1$ and $i$ than including variable $1$ alone (this is because the second variable have entered the model in the $i$-th bi-variate Lasso path when $\lambda\in[\lambda_2^{(i)}-\epsilon,\lambda_2^{(i)})$). We are ready to prove the theorem now, using what we have shown regarding $\lambda_1$ and $\lambda_2$. We next compute the false positive rate and false negative rate given $(\beta_j,\beta_{j+1})=(0,0), (0,\tau_p), (\tau_p,0), (\tau_p,\tau_p), (-\tau_p,\tau_p)$ by deriving upper and lower bounds for those rates. 

We first establish some noatations. For the four-variate Lasso problem (\ref{Lasso-EC}), let $A_i$ denotes the event that variable $i$ is the first one entering the model, $A_{i_1,i_2}$ denotes the event that variable $i_1$ and $i_2$ are the first two entering the model (ignoring the order between $i_1$ and $i_2$) and $A_{i_1 \to i_2}$ denotes the event that variable $i_1$ is the first one and variable $i_2$ is the second one entering the model. Let $L_{i_1,i_2}$ denote the bi-variate Lasso problem with $y$ as the response and $x_{i_1}$, $x_{i_2}$ as the variables. Let $h\equiv (y^T x_j,y^T x_{j+1},y^T \tilde{x}_{j},y^T \tilde{x}_{j+1})$, then $h\sim \mathcal{N}(\mu,G)$ with $\mu=G(\beta_j,\beta_{j+1},0,0)^T$ and $G=((1,\rho,0,\rho)^T,(\rho,1,\rho,0)^T,(0,\rho,1,\rho)^T,(\rho,0,\rho,1)^T)$. When not causing any confusing, we write $t_p$ in place of $t_p(u)$ for simplicity.
\begin{itemize}
    \item When $(\beta_j,\beta_{j+1})=(0,0)$,
        \vspace{0.2cm}
    \beq \label{EC-1-1}
    \mathbb{P}\bigl\{W_j>t_p\big|(\beta_j,\beta_{j+1})=(0,0)\bigr\}=L_p p^{-u}.\eeq
    
     To derive a lower bound for $\mathbb{P}\bigl\{W_j>t_p\big|(\beta_j,\beta_{j+1})=(0,0)\bigr\}$, we look for a point in the region (or on the boundary of the region) that choose variable $j$ as a signal and apply Lemma~\ref{lem:tool}. The point we choose is $p_1=(t_p,\rho t_p,0,\rho t_p)^T$ where $t_p=\sqrt{2u\log(p)}$. It's obvious that when $h=p_1$, variable $j$ is the first one entering the Lasso path. Though $h=p_1$ is in the rejection region, it is also on the boundary of the region that choose variable $j$ as a signal because slight increasing the first coordinate will result in variable $j$ being selected. Since $h\sim \mathcal{N}(\mu_1,G)$ with $\mu_1=\textbf{0}$, by Lemma~\ref{lem:tool},
     \[\mathbb{P}\bigl\{W_j>t_p\big|(\beta_j,\beta_{j+1})=(0,0)\bigr\}\geq L_p p^{-(p_1-\mu_1)^T G^{-1}(p_1-\mu_1)/2\log(p)}=L_p p^{-u}.\]
     The upper bound is straight forward by considering the first variable-$i$ entering the model and notice that $W_i\sim \mathcal{N}(0,1)$:
    \beq \label{EC-1-2}
    \begin{split}
           \mathbb{P}\bigl\{W_j>t_p\big|(\beta_j,\beta_{j+1})=(0,0)\bigr\}
           =&\sum_{i}\mathbb{P}\bigl\{W_j>t_p,A_i\big|(\beta_j,\beta_{j+1})=(0,0)\bigr\}\\
           \leq & \sum_{i}\mathbb{P}\bigl\{W_i>t_p\big|(\beta_j,\beta_{j+1})=(0,0)\bigr\}
           = L_p p^{-u}.
    \end{split}
    \eeq
    \item When $(\beta_j,\beta_{j+1})=(0,\tau_p)$,
        \vspace{0.2cm}
    \beq \label{EC-2-1}
    \mathbb{P}\bigl\{W_j>t_p\big|(\beta_j,\beta_{j+1})=(0,\tau_p)\bigr\}\geq L_p p^{-(\sqrt{u}-\rho\sqrt{r})^2-(\xi_\rho\sqrt{r}-\eta_\rho\sqrt{u})_+^2+(\sqrt{r}-\sqrt{u})_+^2},\eeq
    \beq \label{EC-2-2}
    \mathbb{P}\bigl\{W_j>t_p,A\big|(\beta_j,\beta_{j+1})=(0,\tau_p)\bigr\}\leq L_p p^{-u}\eeq
    for $A=A_{j+p+1\to j}, A_{j+1, j+p+1}$ and
    \beq \label{EC-2-3}
    \mathbb{P}\bigl\{W_j>t_p,A\big|(\beta_j,\beta_{j+1})=(0,\tau_p)\bigr\}\leq L_p p^{-(\sqrt{u}-\rho\sqrt{r})^2-(\xi_\rho\sqrt{r}-\eta_\rho\sqrt{u})_+^2+(\sqrt{r}-\sqrt{u})_+^2}\eeq
    for $A=A_{j,j+1},A_{j,j+p},A_{j\to j+p+1}$.
    \vspace{0.2cm}
            
    This time we choose \[p_2^T=
\left\{  
\begin{array}{ll}
(t_p,\rho t_p+(1-\rho^2)\tau_p,\rho \tau_p,\rho t_p-\rho^2\tau_p), & (1+\rho)\tau_p\leq t_p, \\
(t_p,t_p,\frac{\rho}{1+\rho}t_p,\frac{\rho}{1+\rho}t_p), & \tau_p \leq t_p < (1+\rho)\tau_p, \\
(t_p+\rho(\tau_p-t_p),\tau_p,\rho(\tau_p-t_p)+\frac{\rho}{1+\rho}t_p,\frac{\rho}{1+\rho}t_p), & t_p < \tau_p.
\end{array}
\right.  
\]
When $h=p_2$ and $ t_p\geq \tau_p$, variable $j$ is the first variable entering the four-variate Lasso path with $W_j=t_p$; when $h=p_2$ and $t_p< \tau_p$, variable $j+1$ is the first and $j$ is the second variable entering the Lasso path with $W_j=t_p$ and $W_{j+1}=\tau_p$. $h=p_2$ is on the boundary of the region that chooses variable $j$ as a signal. Since $h\sim \mathcal{N}(\mu_2,G)$ with $\mu_2=(\rho\tau_p,\tau_p,\rho\tau_p,0)^T$, by Lemma~\ref{lem:tool}, 
     \[
     \begin{split}
        \mathbb{P}\bigl\{W_j>t_p\big|(\beta_j,\beta_{j+1})=(0,\tau_p)\bigr\}&\geq L_p p^{-(p_2-\mu_2)^T G^{-1}(p_2-\mu_2)/2\log(p)}\\
        &=L_p p^{-(\sqrt{u}-\rho\sqrt{r})^2-(\xi_\rho\sqrt{r}-\eta_\rho\sqrt{u})_+^2+(\sqrt{r}-\sqrt{u})_+^2}.
     \end{split}
    \]
    
    When $A_{j+1,j+p+1}$ occurs, since by our argument on $\lambda_1$ and $\lambda_2$, $Z_{j+1}$ and $Z_{j+p+1}$ are the $\lambda$ value at which the variables enter the solution path in the bi-variate Lasso problem $L_{j+1,j+p+1}$. Therefore, $Z_{j+1}=|y^T x_{j+1}|, Z_{j+p+1}=|y^T\tilde{x}_{j+1}|$. We notice that $Z_{j+p+1}>Z_j>t_p$ and marginally $y^T\tilde{x}_{j+1}\sim \mathcal{N}(0,1)$, so 
\[
\begin{split}
   &\mathbb{P}\bigl\{W_j>t_p,A_{j+1,j+p+1}\big|(\beta_j,\beta_{j+1})=(0,\tau_p)\bigr\}\\
   &\leq \mathbb{P}\bigl\{|y^T\tilde{x}_{j+1}|>t_p\big|(\beta_j,\beta_{j+1})=(0,\tau_p)\bigr\}= L_p p^{-u}.
\end{split}
\]
Above inequality also holds for $A_{j+p+1\to j}$ since if variable $j+p+1$ is the first entering the Lasso path, then we must have $|y^T\tilde{x}_{j+1}|=Z_{j+p+1}>Z_j>t_p$.

When any one of $A_{j,j+1}, A_{j,j+p}, A_{j\to j+p+1}$ occurs, it implies in the bi-variate Lasso problem $L_{j,j+1}$, the largest $\lambda$ such that variable 1 enters the model for the first time is equal to $W_j$, thus larger than $t_p$. In other words, if variable $j$ is a false positive using Knockoff for variable selection, then it is also a false positive when using bi-variate Lasso $L_{j,j+1}$. This means $\mathbb{P}\bigl\{W_j>t_p,A\big|(\beta_j,\beta_{j+1})=(0,\tau_p)\bigr\}$ is upper bounded by the corresponding false positive rate of Lasso, which is $L_p p^{-(\sqrt{u}-\rho\sqrt{r})^2-(\xi_\rho\sqrt{r}-\eta_\rho\sqrt{u})_+^2+(\sqrt{r}-\sqrt{u})_+^2}$, for $A=A_{j,j+1}, A_{j,j+p}, A_{j\to j+p+1}$.

Since $A_{j+1,j+p}$ and $A_{j+p,j+p+1}$ can never occur when $W_j>0$, (\ref{EC-2-2}) and (\ref{EC-2-3}) implies 
\beq \mathbb{P}\bigl\{W_j>t_p\big|(\beta_j,\beta_{j+1})=(0,\tau_p)\bigr\}\leq L_p p^{-\min\{u,(\sqrt{u}-\rho\sqrt{r})^2+(\xi_\rho\sqrt{r}-\eta_\rho\sqrt{u})_+^2-(\sqrt{r}-\sqrt{u})_+^2\}}.\eeq
Further coupled with (\ref{EC-1-1}) and (\ref{EC-2-1}), we have 
\beq \label{EC-FP}
\mathbb{P}\bigl\{W_j>t_p,\beta_j=0 \bigr\}=L_p p^{-\min\{u,\vartheta+(\sqrt{u}-\rho\sqrt{r})^2+(\xi_\rho\sqrt{r}-\eta_\rho\sqrt{u})_+^2-(\sqrt{r}-\sqrt{u})_+^2\}}. \eeq
    \item When $(\beta_j,\beta_{j+1})=(\tau_p,0)$,
    \vspace{0.2cm}
    \beq \label{EC-3-1}
    \mathbb{P}\bigl\{W_j\leq t_p\big|(\beta_j,\beta_{j+1})=(\tau_p,0)\bigr\}\geq L_p p^{-[(\sqrt{r}-\sqrt{u})_+]^2},\eeq
    and
    \beq \label{EC-3-2}
    \mathbb{P}\bigl\{W_j\leq t_p\big|(\beta_j,\beta_{j+1})=(\tau_p,0)\bigr\}\leq L_p p^{\vartheta-f^+_{\text{Hamm}}(u,r,\vartheta)}.\eeq
    
    Let $p_3=(t_p,\rho t_p,0,\rho t_p)^T$. when $h=p_3$, variable $j$ is the first variable entering the Lasso path and $p_3$ is in the region of rejecting variable $j$ as a signal. Since $h\sim \mathcal{N}(\mu_3,G)$ with $\mu_3=(\tau_p,\rho \tau_p,0,\rho\tau_p)^T$, by Lemma~\ref{lem:tool}, 
     \[
     \begin{split}
        \mathbb{P}\bigl\{W_j\leq t_p\big|(\beta_j,\beta_{j+1})=(\tau_p,0)\bigr\}&\geq L_p p^{-(p_3-\mu_3)^T G^{-1}(p_3-\mu_3)/2\log(p)}\\
        &=L_p p^{-[(\sqrt{r}-\sqrt{u})_+]^2}.
     \end{split}
    \] 

Before we prove (\ref{EC-3-2}), we first analysis $f^+_{\text{Hamm}}(u,r,\vartheta)$. By simply calculation, we find 
the optimal value of $u$ that maximize $f^+_{\text{Hamm}}(u,r,\vartheta)$ given $r,\vartheta$ is 
\[u^*=
\left\{  
\begin{array}{ll}
\frac{1+\rho}{(\sqrt{1+\rho}+\sqrt{1-\rho})^2}r, & \vartheta\leq \frac{2\rho}{(\sqrt{1+\rho}+\sqrt{1-\rho})^2}r, \\
\frac{(r+\vartheta)^2}{4r}, & \frac{2\rho}{(\sqrt{1+\rho}+\sqrt{1-\rho})^2}r\leq \vartheta < r, \\
\vartheta, & r<\vartheta.
\end{array}
\right.  
\]
This implies $u^*\geq \frac{1+\rho}{(\sqrt{1+\rho}+\sqrt{1-\rho})^2}r$ regardless of the relationship of $\vartheta$ and $r$. Consider $r,\vartheta$ as fixed,
$f^+_{\text{Hamm}}(r,u,\vartheta)$ as a function of $u$ is monotonically non-decreasing in $[0,u^*]$ and monotonically non-increasing in $[u^*,\infty)$. $f^+_{\text{Hamm}}(r,\vartheta)=\vartheta+[(\sqrt{r}-\sqrt{u})_+-((1-\xi_\rho)\sqrt{r}-(1-\eta_\rho)\sqrt{u})_+]^2$ if and only if $u>u^*$. Since $(1-\xi_\rho)\sqrt{r}-(1-\eta_\rho)\sqrt{u^*}<0$, $(1-\xi_\rho)\sqrt{r}-(1-\eta_\rho)\sqrt{u}<0$ for all $u>u^*$, which implies $f^+_{\text{Hamm}}(r,\vartheta)=\vartheta+[(\sqrt{r}-\sqrt{u})_+]^2$ when $u>u^*$. Therefore,
\[
\begin{split}
    f^+_{\text{Hamm}}(r,u,\vartheta)=\min\{&u,\vartheta+(\sqrt{u}-|\rho|\sqrt{r})^2+((\xi_\rho \sqrt{r}-\eta_\rho \sqrt{u})_+)^2-((\sqrt{r}-\sqrt{u})_+)^2,\\
    &\vartheta+[(\sqrt{r}-\sqrt{u})_+]^2\}.
\end{split}
\]
Now, we show that (\ref{EC-3-2}) holds for $u\geq u^*$. This would implies (\ref{EC-3-2}) for all $u\geq 0$, since the false negative rate $\mathbb{P}\bigl\{W_j\leq t_p(u)\big|(\beta_j,\beta_{j+1})=(\tau_p,0)\bigr\}$ is monotone non-decreasing with $u$, so for $u< u^*$, $\mathbb{P}\bigl\{W_j\leq t_p(u)\big|(\beta_j,\beta_{j+1})=(\tau_p,0)\bigr\}\leq \mathbb{P}\bigl\{W_j\leq t_p(u^*)\big|(\beta_j,\beta_{j+1})=(\tau_p,0)\bigr\}\leq L_p p^{\vartheta-f^+_{\text{Hamm}}(r,u^*,\vartheta)}\leq L_p p^{\vartheta-f^+_{\text{Hamm}}(r,u,\vartheta)}$.

Assume $u\geq u^*$, so $u\geq \frac{1+\rho}{(\sqrt{1+\rho}+\sqrt{1-\rho})^2}r$ and 
\beq \label{ueq}
-[(\sqrt{r}-\sqrt{u})_+]^2\geq -\Big(\frac{\sqrt{1-\rho}}{\sqrt{1+\rho}+\sqrt{1-\rho}}\Big)^2r\geq -(2-\sqrt{3})(1-\rho)r\geq -\frac{1-\rho}{2}r\geq -\frac{1}{2}r.\eeq
We next prove (\ref{EC-3-2}) by showing that 
    \beq
    \mathbb{P}\bigl\{W_j\leq t_p,A \big|(\beta_j,\beta_{j+1})=(\tau_p,0)\bigr\}\leq L_p p^{-[(\sqrt{r}-\sqrt{u})_+]^2}\eeq
holds for $A=A_j,A_{j+1},A_{j+p},A_{j+p+1}$ and $u\geq u^*$. Respectively, 
\[
\begin{split}
    \mathbb{P}\bigl\{W_j\leq t_p,A_j \big|(\beta_j,\beta_{j+1})=(\tau_p,0)\bigr\}
&\leq \mathbb{P}\bigl\{|y^Tx_j|\leq t_p\big|(\beta_j,\beta_{j+1})=(\tau_p,0)\bigr\}\\
&=L_p p^{-[(\sqrt{r}-\sqrt{u})_+]^2},
\end{split}
\]
and by symmetry and (\ref{ueq}),
\[
\begin{split}
 \mathbb{P}\bigl\{W_j\leq t_p,A_{j+1} \big|(\beta_j,\beta_{j+1})=(\tau_p,0)\bigr\}=\mathbb{P}\bigl\{W_j\leq t_p,A_{j+p+1} \big|(\beta_j,\beta_{j+1})=(\tau_p,0)\bigr\}   \\
  \leq \mathbb{P}\bigl\{|y^Tx_j|\leq |y^Tx_{j+p+1}|\big|(\beta_j,\beta_{j+1})=(\tau_p,0)\bigr\} \leq L_p p^{-\frac{1-\rho}{2}r}\leq L_p p^{-[(\sqrt{r}-\sqrt{u})_+]^2},
\end{split}
\]
\[
\begin{split}
\mathbb{P}\bigl\{W_j\leq t_p,A_{j+p} \big|(\beta_j,\beta_{j+1})=(\tau_p,0)\bigr\}  
  &\leq \mathbb{P}\bigl\{|y^Tx_j|\leq |y^Tx_{j+p}|\big|(\beta_j,\beta_{j+1})=(\tau_p,0)\bigr\} \\
  &\leq L_p p^{-\frac{1}{2}r} \leq L_p p^{-[(\sqrt{r}-\sqrt{u})_+]^2}.
\end{split}
\]
(\ref{EC-3-2}) is immediate by $[(\sqrt{r}-\sqrt{u})_+]^2 \geq f^+_{\text{Hamm}}(r,u,\vartheta)-\vartheta$.
\vspace{0.2cm}

    \item When $(\beta_j,\beta_{j+1})=(\tau_p,\tau_p)$,
    \vspace{0.2cm}
    \beq \label{EC-4-1}
    \mathbb{P}\bigl\{W_j\leq t_p\big|(\beta_j,\beta_{j+1})=(\tau_p,\tau_p)\bigr\}\leq L_p p^{\vartheta-f^+_{\text{Hamm}}(u,r,\vartheta)}.\eeq
More precisely, we will prove that
    \beq \label{EC-4-1_}
    \mathbb{P}\bigl\{W_j\leq t_p\big|(\beta_j,\beta_{j+1})=(\tau_p,\tau_p)\bigr\}\leq L_p p^{-[(\sqrt{r}-\sqrt{u})_+]^2}\eeq
holds for $u\geq u^*$, thus implies (\ref{EC-4-1}). We prove (\ref{EC-4-1_}) by showing 
\beq \label{EC-4-1__}
    \mathbb{P}\bigl\{W_j\leq t_p, A\big|(\beta_j,\beta_{j+1})=(\tau_p,\tau_p)\bigr\}\leq L_p p^{-[(\sqrt{r}-\sqrt{u})_+]^2}\eeq
holds for $A=A_j,A_{j+1\to j},A_{j+1\to j+p},A_{j+1\to j+p+1},A_{j+p},A_{j+p+1}$ and $u\geq u^*$, which cover all possibilities.
Respectively,
\[ 
\begin{split}
    \mathbb{P}\bigl\{W_j\leq t_p, A_j\big|(\beta_j,\beta_{j+1})=(\tau_p,\tau_p)\bigr\}& \leq \mathbb{P}\bigl\{|y^T x_j|\leq t_p\big|(\beta_j,\beta_{j+1})=(\tau_p,\tau_p)\bigr\}\\
    &\leq L_p p^{-[((1+\rho)\sqrt{r}-\sqrt{u})_+]^2} \leq L_p p^{-[(\sqrt{r}-\sqrt{u})_+]^2},
\end{split}
\]
\[ 
\begin{split}
    \mathbb{P}\bigl\{W_j\leq t_p, A_{j+p}\big|(\beta_j,\beta_{j+1})=(\tau_p,\tau_p)\bigr\}& \leq \mathbb{P}\bigl\{|y^T x_j|\leq |y^T \tilde{x}_{j}|\big|(\beta_j,\beta_{j+1})=(\tau_p,\tau_p)\bigr\}\\
    &\leq L_p p^{-\frac{1}{2}r} \leq L_p p^{-[(\sqrt{r}-\sqrt{u})_+]^2},
\end{split}
\]
\[ 
\begin{split}
    \mathbb{P}\bigl\{W_j\leq t_p, A_{j+p+1}\big|(\beta_j,\beta_{j+1})=(\tau_p,\tau_p)\bigr\}& \leq \mathbb{P}\bigl\{|y^T x_j|\leq |y^T \tilde{x}_{j+1}|\big|(\beta_j,\beta_{j+1})=(\tau_p,\tau_p)\bigr\}\\
    &\leq L_p p^{-\frac{1}{2(1-\rho)}r} \leq L_p p^{-[(\sqrt{r}-\sqrt{u})_+]^2}.
\end{split}
\]
When $A_{j+1\to j}$ occurs, the bi-variate Lasso problem $L_{j,j+1}$ shares the same $\lambda_1$ and $\lambda_2$ with the four-variate Lasso problem. So variable $j$ is a false negative when doing variable selection using the bi-variate Lasso $L_{j,j+1}$ given $W_j\leq t_p$, which implies $\mathbb{P}\bigl\{W_j\leq t_p,A_{j+1\to j}\big|(\beta_j,\beta_{j+1})=(\tau_p,\tau_p)\bigr\}$ is upper bounded by the corresponding false negative rate of Lasso, which is $L_p p^{-(\xi_\rho\sqrt{r}-\eta_\rho \sqrt{u})_+^2}\leq L_p p^{-[(\sqrt{r}-\sqrt{u})_+]^2}$. The last inequality is equivalent to
\[(1-\sqrt{1-\rho^2})\sqrt{r}\leq \Big(1-\sqrt{\frac{1-\rho}{1+\rho}}\Big)\sqrt{u}.\]
By (\ref{ueq}), the right hand side is no smaller than $\sqrt{r}$, thus no smaller than the left hand side.

When $A_{j+1\to j+p}$ occurs, we know variable $j+p$ instead of variable $j$ is the second one entering the Lasso path. This means the $\lambda_2$ (the $\lambda$ value when the second variable entering Lasso path) of the bi-variate Lasso problem $L_{j+1,j+p}$ is larger than the $\lambda_2$ of the bi-variate Lasso problem $L_{j,j+1}$. Since we have derived the explicit expression of $\lambda_2$ in bi-variate Lasso problems, when $y^T x_{j+1}\geq 0$, we must have 
\[\max\{\frac{y^T x_j-\rho y^T x_{j+1}}{1-\rho},\frac{y^T x_j-\rho y^T x_{j+1}}{-1-\rho}\}< \max\{\frac{y^T \tilde{x}_j-\rho y^T x_{j+1}}{1-\rho},\frac{y^T \tilde{x}_j-\rho y^T x_{j+1}}{-1-\rho}\}.\]
Therefore, $A_{j+1\to j+p}$ implies one the three following events must occur:
\[y^T x_{j+1}< 0,\frac{y^T x_j-\rho y^T x_{j+1}}{1-\rho}<\frac{y^T \tilde{x}_j-\rho y^T x_{j+1}}{1-\rho}, \frac{y^T x_j-\rho y^T x_{j+1}}{1-\rho}<\frac{y^T \tilde{x}_j-\rho y^T x_{j+1}}{-1-\rho}\] 
The probability of these three events given $(\beta_j,\beta_{j+1})=(\tau_p,\tau_p)$ are $L_p p^{-(1+\rho)^2r}$, $L_p p^{-\frac{r}{2}}$ and $L_p p^{-\frac{(1+2\rho)^2(1-\rho)}{2(1+\rho)}r}$, all of which are upper bounded by $L_p p^{-[(\sqrt{r}-\sqrt{u})_+]^2}$ when $u\geq u^*$.

When $A_{j+1\to j+p+1}$ occurs, the $\lambda_2$ of the bi-variate Lasso problem $L_{j+1,j+p+1}$ is larger than the $\lambda_2$ of the bi-variate Lasso problem $L_{j,j+1}$. When $y^T x_{j+1}\geq 0$, we must have 
\[\max\{\frac{y^T x_j-\rho y^T x_{j+1}}{1-\rho},\frac{y^T x_j-\rho y^T x_{j+1}}{-1-\rho}\}< |y^T \tilde{x}_{j+1}|.\]
Therefore, $A_{j+1\to j+p+1}$ implies one the three following events must occur:
\[y^T x_{j+1}< 0,\frac{y^T x_j-\rho y^T x_{j+1}}{1-\rho}<y^T \tilde{x}_{j+1}, \frac{y^T x_j-\rho y^T x_{j+1}}{1-\rho}<-y^T \tilde{x}_{j+1}.\]
Respectively, the probability of these three events are $L_p p^{-(1+\rho)^2r}$, $L_p p^{-\frac{r}{2}}$ and $L_p p^{-\frac{(1+2\rho)^2(1-\rho)}{2(1+\rho)}r}$, all of which are upper bounded by $L_p p^{-[(\sqrt{r}-\sqrt{u})_+]^2}$ when $u\geq u^*$. From here we have verified (\ref{EC-4-1__}), thus implies (\ref{EC-4-1}).

From (\ref{EC-FP}) and (\ref{EC-3-1}), we have 
\beq 
\mathbb{P}\bigl\{W_j> t_p,\beta_j=0\bigr\}+\mathbb{P}\bigl\{W_j\leq t_p,\beta_j=\tau_p\bigr\}\geq L_p p^{-f^+_{\text{Hamm}}(r,u,\vartheta)}.
\eeq
\beq 
\begin{split}
\mathbb{P}\bigl\{W_j\leq t_p,\beta_j=\tau_p\bigr\}=&p^{-\vartheta}\times \mathbb{P}\bigl\{W_j\leq t_p\big|(\beta_j,\beta_{j+1})=(\tau_p,0)\bigr\}\\
&+p^{-2\vartheta}\times \mathbb{P}\bigl\{W_j\leq t_p\big|(\beta_j,\beta_{j+1})=(\tau_p,\tau_p)\bigr\}\\
\leq & L_p p^{-f^+_{\text{Hamm}}(r,u,\vartheta)}
\end{split}
\eeq
Since (\ref{EC-FP}) also implies $\mathbb{P}\bigl\{W_j> t_p,\beta_j=0\bigr\}\leq L_p p^{-f^+_{\text{Hamm}}(r,u,\vartheta)}$, we know 
\beq \label{EC-positive}
\mathbb{P}\bigl\{W_j> t_p,\beta_j=0\bigr\}+\mathbb{P}\bigl\{W_j\leq t_p,\beta_j=\tau_p\bigr\}= L_p p^{-f^+_{\text{Hamm}}(r,u,\vartheta)}.
\eeq
    \item When $(\beta_j,\beta_{j+1})=(-\tau_p,\tau_p)$,
    \vspace{0.2cm}
    \beq \label{EC-5-1}
    \mathbb{P}\bigl\{W_j\leq t_p\big|(\beta_j,\beta_{j+1})=(-\tau_p,\tau_p)\bigr\}\geq L_p p^{-((\xi_\rho \sqrt{r}-\eta_\rho^{-1} \sqrt{u})_+)^2},\eeq
    \beq \label{EC-5-2}
    \mathbb{P}\bigl\{W_j\leq t_p\big|(\beta_j,\beta_{j+1})=(-\tau_p,\tau_p)\bigr\}\geq L_p p^{-\frac{(1-2\rho)^2(1+\rho)}{2(1-\rho)}r},\eeq
    and
    \beq \label{EC-5-3}
    \mathbb{P}\bigl\{W_j\leq t_p\big|(\beta_j,\beta_{j+1})=(-\tau_p,\tau_p)\bigr\}\leq L_p p^{-\min\{((\xi_\rho \sqrt{r}-\eta_\rho^{-1} \sqrt{u})_+)^2,\frac{(1-2\rho)^2(1+\rho)}{2(1-\rho)}r,-\vartheta+f^+_{\text{Hamm}}(u,r,\vartheta)\}}.\eeq
 Let 
 \[p_4^T=
\left\{  
\begin{array}{ll}
(-(1-\rho)\tau_p,(1-\rho)\tau_p,\rho \tau_p,-\rho \tau_p), & (1-\rho)\tau_p\leq t_p, \\
(\rho(1-\rho)\tau_p-(1+\rho)t_p,(1-\rho)\tau_p,\rho(1-\rho)\tau_p+\frac{\rho^2}{1-\rho}t_p,-\frac{\rho}{1-\rho}t_p), & (1-\rho)\tau_p> t_p. \\
\end{array}
\right.
\]
When $h=p_4$ and $(1-\rho)\tau_p\leq t_p$, variable $j$ is the first variable entering the Lasso path with $W_j=(1-\rho)\tau_p\leq t_p$; when $h=p_4$ and $(1-\rho)\tau_p> t_p$, $j+1$ is the first and $j$ is the second variable entering the Lasso path with $W_j=t_p$. Regardless of the relationship between $\tau_p$ and $t_p$, $h=p_4$ is always in the region of rejecting $j$ as a signal. Since $h\sim \mathcal{N}(\mu_4,G)$ with $\mu_4=(-(1-\rho)\tau_p,(1-\rho)\tau_p,\rho \tau_p,-\rho \tau_p)^T$, by Lemma~\ref{lem:tool}, 
     \[
     \begin{split}
        \mathbb{P}\bigl\{W_j\leq t_p\big|(\beta_j,\beta_{j+1})=(-\tau_p,\tau_p)\bigr\}&\geq L_p p^{-(p_4-\mu_4)^T G^{-1}(p_4-\mu_4)/2\log(p)}\\
        &=L_p p^{-((\xi_\rho \sqrt{r}-\eta_\rho^{-1} \sqrt{u})_+)^2}.
     \end{split}
    \] 
Let $p_5^T=\Big(\frac{4\rho^3-2\rho^2+\rho-1}{2(1-\rho)}\tau_p,(1-\rho)\tau_p ,\frac{1+2\rho-4\rho^2}{2}\tau_p,-\frac{\rho^2}{1-\rho}\tau_p\Big)$. When $h=p_5$, variable $j+1$ is the first one entering the Lasso path with $W_{j+1}=(1-\rho)\tau_p$, if we slightly increase the value of the third coordinate of $p_5$, then it falls in the region of rejecting $j$ as a signal since variable $j+p$ is the second variable entering the Lasso path. This implies $h=p_5$ in on the boundary of the region that rejects $j$ as a signal, by Lemma~\ref{lem:tool}, 
     \[
     \begin{split}
        \mathbb{P}\bigl\{W_j\leq t_p\big|(\beta_j,\beta_{j+1})=(-\tau_p,\tau_p)\bigr\}&\geq L_p p^{-(p_5-\mu_4)^T G^{-1}(p_5-\mu_4)/2\log(p)}\\
        &=L_p p^{-\frac{(1-2\rho)^2(1+\rho)}{2(1-\rho)}r}.
     \end{split}
    \] 
Next, we show that 
    \beq \label{EC-5-4}
    \mathbb{P}\bigl\{W_j\leq t_p, A\big|(\beta_j,\beta_{j+1})=(-\tau_p,\tau_p)\bigr\}\leq L_p p^{-\min\{((\xi_\rho \sqrt{r}-\eta_\rho^{-1} \sqrt{u})_+)^2,\frac{(1-2\rho)^2(1+\rho)}{2(1-\rho)}r,-\vartheta+f^+_{\text{Hamm}}(u,r,\vartheta)\}}.\eeq
holds for $A=A_j,A_{j+1\to j},A_{j+1\to j+p},A_{j+1\to j+p+1},A_{j+p},A_{j+p+1}$, which cover all possibilities.

When $A=A_j$ or $A_{j+1 \to j}$ occurs, as previously discussed, variable $j$ is a false negative when doing variable selection using the bi-variate Lasso $L_{j,j+1}$ given $W_j\leq t_p$, which implies $\mathbb{P}\bigl\{W_j\leq t_p,A\big|(\beta_j,\beta_{j+1})=(-\tau_p,\tau_p)\bigr\}$ is upper bounded by the corresponding false negative rate of Lasso, which is $L_p p^{-((\xi_\rho \sqrt{r}-\eta_\rho^{-1} \sqrt{u})_+)^2}$. 

When $A_{j+1\to j+p}$ occurs, the $\lambda_2$ of the bi-variate Lasso problem $L_{j+1,j+p}$ is larger than the $\lambda_2$ of the bi-variate Lasso problem $L_{j,j+1}$. When $y^T x_{j+1}\geq 0$, we must have 
\[\max\{\frac{y^T x_j-\rho y^T x_{j+1}}{1-\rho},\frac{y^T x_j-\rho y^T x_{j+1}}{-1-\rho}\}< \max\{\frac{y^T \tilde{x}_j-\rho y^T x_{j+1}}{1-\rho},\frac{y^T \tilde{x}_j-\rho y^T x_{j+1}}{-1-\rho}\}.\]
Therefore, $A_{j+1\to j+p}$ implies one of the three following events must occur:
\[y^T x_{j+1}< 0,\frac{y^T x_j-\rho y^T x_{j+1}}{-1-\rho}<\frac{y^T \tilde{x}_j-\rho y^T x_{j+1}}{1-\rho}, \frac{y^T x_j-\rho y^T x_{j+1}}{-1-\rho}<\frac{y^T \tilde{x}_j-\rho y^T x_{j+1}}{-1-\rho}\] 
The probability of these three events are $L_p p^{-(1-\rho)^2r}$, $L_p p^{-\frac{(1-2\rho)^2(1+\rho)}{2(1-\rho)}r}$ and $L_p p^{-\frac{r}{2}}$, all of which are upper bounded by $L_p p^{-\min\{\frac{(1-2\rho)^2(1+\rho)}{2(1-\rho)}r,-\vartheta+f^+_{\text{Hamm}}(u,r,\vartheta)\}}$.

When $A_{j+1\to j+p+1}$ occurs, the $\lambda_2$ of the bi-variate Lasso problem $L_{j+1,j+p+1}$ is larger than the $\lambda_2$ of the bi-variate Lasso problem $L_{j,j+1}$. When $y^T x_{j+1}\geq 0$, we must have 
\[\max\{\frac{y^T x_j-\rho y^T x_{j+1}}{1-\rho},\frac{y^T x_j-\rho y^T x_{j+1}}{-1-\rho}\}< |y^T \tilde{x}_{j+1}|.\]
Therefore, $A_{j+1\to j+p+1}$ implies one of the three following events must occur:
\[y^T x_{j+1}< 0,\frac{y^T x_j-\rho y^T x_{j+1}}{-1-\rho}<y^T \tilde{x}_{j+1}, \frac{y^T x_j-\rho y^T x_{j+1}}{-1-\rho}<-y^T \tilde{x}_{j+1}.\]
The probability of these three events are $L_p p^{-(1-\rho)^2r}$,  $L_p p^{-\frac{r}{2}}$ and $L_p p^{-\frac{(1-2\rho)^2(1+\rho)}{2(1-\rho)}r}$, all of which are upper bounded by $L_p p^{-\min\{\frac{(1-2\rho)^2(1+\rho)}{2(1-\rho)}r,-\vartheta+f^+_{\text{Hamm}}(u,r,\vartheta)\}}$.

When $A_{j+p}$ occurs, then $|y^T \tilde{x}_j|>|y^T x_j|$ and $|y^T \tilde{x}_j|>|y^T x_{j+1}|$. If $y^T \tilde{x}_j>0$, we further have $(y^T \tilde{x}_j-y^T x_{j+1})+\frac{1}{2\rho+1}(y^T \tilde{x}_j+y^T x_j)> 0$; if $y^T \tilde{x}_j\leq 0$, we further have $y^T \tilde{x}_j+y^T x_{j+1}<0$. Therefore,
\[
     \begin{split}
        \mathbb{P}&\bigl\{W_j\leq t_p,A_{j+p}\big|(\beta_j,\beta_{j+1})=(-\tau_p,\tau_p)\bigr\}\\
        \leq & \mathbb{P}\bigl\{(y^T \tilde{x}_j-y^T x_{j+1})+\frac{1}{2\rho+1}(y^T \tilde{x}_j+y^T x_j)> 0\big|(\beta_j,\beta_{j+1})=(-\tau_p,\tau_p)\bigr\}\\
        &+\mathbb{P}\bigl\{y^T \tilde{x}_j+y^T x_{j+1}<0\big|(\beta_j,\beta_{j+1})=(-\tau_p,\tau_p)\bigr\}\\
        \leq& L_p p^{-\frac{2(1-2\rho)^2(1+\rho)}{3-4\rho^2}r}+L_p p^{-\frac{r}{2}}\leq L_p p^{-\min\{\frac{(1-2\rho)^2(1+\rho)}{2(1-\rho)}r,-\vartheta+f^+_{\text{Hamm}}(u,r,\vartheta)\}}.
     \end{split}
    \]
For $A=A_{j+p+1}$, (\ref{EC-5-4}) is immediate due to the symmetry between variable $j+p$ and $j+p+1$. 

Now consider the case where $\beta_j$ takes value in $\{0,-\tau_p\}$ and $\beta_{j+1}$ takes value in $\{0,\tau_p\}$, this corresponds to the $\rho<0$ case (we flipped the sign of $\rho$ and $\beta_j$ simultaneously).
By (\ref{EC-1-1}), (\ref{EC-2-1}), (\ref{EC-3-1}), (\ref{EC-5-1}) and (\ref{EC-5-2}), we know 
\beq 
\begin{split}
    \mathbb{P}\bigl\{W_j> t_p,\beta_j=0\bigr\}+\mathbb{P}\bigl\{W_j\leq t_p,\beta_j=-\tau_p\bigr\}\\
\geq L_p p^{-\min\{f^+_{\text{Hamm}}(u,r,\vartheta),2\vartheta+((\xi_\rho \sqrt{r}-\eta_\rho^{-1} \sqrt{u})_+)^2,2\vartheta+\frac{(1-2|\rho|)^2(1+|\rho|)}{2(1-|\rho|)}r\}}.
\end{split}
\eeq
Meanwhile, (\ref{EC-1-1}), (\ref{EC-2-2}), (\ref{EC-2-3}), (\ref{EC-3-2}) and (\ref{EC-5-3}) gives
\beq 
\begin{split}
    \mathbb{P}\bigl\{W_j> t_p,\beta_j=0\bigr\}+\mathbb{P}\bigl\{W_j\leq t_p,\beta_j=-\tau_p\bigr\}\\
\leq L_p p^{-\min\{f^+_{\text{Hamm}}(u,r,\vartheta),2\vartheta+((\xi_\rho \sqrt{r}-\eta_\rho^{-1} \sqrt{u})_+)^2,2\vartheta+\frac{(1-2|\rho|)^2(1+|\rho|)}{2(1-|\rho|)}r\}}.
\end{split}
\eeq
Therefore,
\beq \label{EC-negative}
\begin{split}
    \mathbb{P}\bigl\{W_j> t_p,\beta_j=0\bigr\}+\mathbb{P}\bigl\{W_j\leq t_p,\beta_j=-\tau_p\bigr\}\\
= L_p p^{-\min\{f^+_{\text{Hamm}}(u,r,\vartheta),2\vartheta+((\xi_\rho \sqrt{r}-\eta_\rho^{-1} \sqrt{u})_+)^2,2\vartheta+\frac{(1-2|\rho|)^2(1+|\rho|)}{2(1-|\rho|)}r\}}.
\end{split}
\eeq
(\ref{EC-positive}) and (\ref{EC-negative}) complete the proof for Theorem~\ref{thm:knockoff-block2}.
\end{itemize}


\section{Proof of Theorem~\ref{thm:knockoff-CI}}
The only difference of the conditional knockoff from the Equal-correlated knockoff construction is that $x_j^T \tilde{x}_j$ is changed from 0 to $\rho^2$ for $j=1,\cdots,p$. Therefore, $G=((1,\rho,\rho^2,\rho)^T,(\rho,1,\rho,\rho^2)^T$, $(\rho^2,\rho,1,\rho)^T,(\rho,\rho^2,\rho,1)^T)$ is the new gram matrix for the four-variate Lassos (\ref{Lasso-EC}). We follow the same notations and workflow from the previous proof.

\begin{itemize}
    \item When $(\beta_j,\beta_{j+1})=(0,0)$,
        \vspace{0.2cm}
    \beq \label{CI-1-1}
    \mathbb{P}\bigl\{W_j>t_p\big|(\beta_j,\beta_{j+1})=(0,0)\bigr\}=L_p p^{-u}.\eeq
    
Let $p_1=(t_p,\rho t_p,\rho^2 t_p,\rho t_p)^T$ where $t_p=\sqrt{2u\log(p)}$. When $h=p_1$, variable $j$ is the first one entering the Lasso path. Though $h=p_1$ is in the rejection region, it is also on the boundary of the region that choose variable $j$ as a signal. Since $h\sim \mathcal{N}(\mu_1,G)$ with $\mu_1=\textbf{0}$, by Lemma~\ref{lem:tool},
     \[\mathbb{P}\bigl\{W_j>t_p\big|(\beta_j,\beta_{j+1})=(0,0)\bigr\}\geq L_p p^{-(p_1-\mu_1)^T G^{-1}(p_1-\mu_1)/2\log(p)}=L_p p^{-u}.\]
     The upper bound is derived exactly the same as (\ref{EC-1-2}).
        \item When $(\beta_j,\beta_{j+1})=(0,\tau_p)$,
        \vspace{0.2cm}
    \beq \label{CI-2-1}
    \mathbb{P}\bigl\{W_j>t_p\big|(\beta_j,\beta_{j+1})=(0,\tau_p)\bigr\}= L_p p^{-(\sqrt{u}-\rho\sqrt{r})^2-(\xi_\rho\sqrt{r}-\eta_\rho\sqrt{u})_+^2+(\sqrt{r}-\sqrt{u})_+^2}.\eeq

    This time we choose \[p_2^T=
\left\{  
\begin{array}{ll}
(t_p,\rho t_p+(1-\rho^2)\tau_p,\rho^2 t_p+\rho(1-\rho^2) \tau_p,\rho t_p)^T, & (1+\rho)\tau_p\leq t_p, \\
(t_p,t_p,\rho t_p,\rho t_p)^T, & \tau_p \leq t_p < (1+\rho)\tau_p, \\
((1-\rho)t_p+\rho \tau_p,\tau_p,\rho \tau_p,\rho(1-\rho)t_p+\rho^2 \tau_p)^T, & t_p< \tau_p.
\end{array}
\right. 
\]
When $h=p_2$ and $ t_p\geq \tau_p$, variable $j$ is the first variable entering the four-variate Lasso path with $W_j=t_p$; when $h=p_2$ and $t_p< \tau_p$, variable $j+1$ is the first and $j$ is the second variable entering the Lasso path with $W_j=t_p$ and $W_{j+1}=\tau_p$. $h=p_2$ is on the boundary of the region that chooses variable $j$ as a signal. Since $h\sim \mathcal{N}(\mu_2,G)$ with $\mu_2=(\rho\tau_p,\tau_p,\rho\tau_p,\rho^2\tau_p)^T$, by Lemma~\ref{lem:tool}, 
     \beq \label{CI-2-3}
     \begin{split}
        \mathbb{P}\bigl\{W_j>t_p\big|(\beta_j,\beta_{j+1})=(0,\tau_p)\bigr\}&\geq L_p p^{-(p_2-\mu_2)^T G^{-1}(p_2-\mu_2)/2\log(p)}\\
        &=L_p p^{-(\sqrt{u}-\rho\sqrt{r})^2-(\xi_\rho\sqrt{r}-\eta_\rho\sqrt{u})_+^2+(\sqrt{r}-\sqrt{u})_+^2}.
     \end{split}
    \eeq
Next we show that 
\beq \label{CI-2-2}
    \mathbb{P}\bigl\{W_j>t_p,A\big|(\beta_j,\beta_{j+1})=(0,\tau_p)\bigr\}\leq  L_p p^{-(\sqrt{u}-\rho\sqrt{r})^2-(\xi_\rho\sqrt{r}-\eta_\rho\sqrt{u})_+^2+(\sqrt{r}-\sqrt{u})_+^2}\eeq
holds for $A=A_{j,j+1}, A_{j,j+p}, A_{j\to j+p+1},A_{j+p+1\to j},A_{j+1,j+p+1}$, which covers all possibilities.

When any one of $A_{j,j+1}, A_{j,j+p}, A_{j\to j+p+1}$ occurs, same as for EC-knockoff, it implies if variable $j$ is a false positive using Knockoff for variable selection, then it is also a false positive when using bi-variate Lasso $L_{j,j+1}$. So $\mathbb{P}\bigl\{W_j>t_p,A\big|(\beta_j,\beta_{j+1})=(0,\tau_p)\bigr\}$ is upper bounded by the corresponding false positive rate of Lasso, which is $L_p p^{-(\sqrt{u}-\rho\sqrt{r})^2-(\xi_\rho\sqrt{r}-\eta_\rho\sqrt{u})_+^2+(\sqrt{r}-\sqrt{u})_+^2}$, for $A=A_{j,j+1}, A_{j,j+p}, A_{j\to j+p+1}$.

When $A=A_{j+p+1\to j}$, $j+p+1$ is the first variable entering the model in the four-variate Lasso problem, thus it's also the first variable entering the model in the bi-variate Lasso problem $L_{j+1,j+p+1}$ and $L_{j,j+p+1}$.
Variable $j+p+1$ gets picked up as a signal in $L_{j+1,j+p+1}$ implies
\[
\begin{split}
\mathbb{P}\bigl\{W_j>t_p,A_{j+p+1\to j}\big|(\beta_j,\beta_{j+1})=(0,\tau_p)\bigr\}
&\leq L_p p^{-(\sqrt{u}-|\rho^2|\sqrt{r})^2-(\xi_{\rho^2}\sqrt{r}-\eta_{\rho^2}\sqrt{u})_+^2+(\sqrt{r}-\sqrt{u})_+^2}\\
    &\leq L_p p^{-(\sqrt{u}-|\rho|\sqrt{r})^2-(\xi_{\rho}\sqrt{r}-\eta_{\rho}\sqrt{u})_+^2+(\sqrt{r}-\sqrt{u})_+^2}
\end{split}
\]
when $u\geq (1+\rho)^2r$ or $u\leq (1+\rho^2)^2r$.

Now consider bi-variate Lasso problem $L_{j,j+p+1}$ given $(1+\rho^2)^2r<u< (1+\rho)^2r$. Variable $j, j+p+1$ both get picked up as signals with $j+p+1$ entering the model first given $W_j>t_p$. This implies $(y^T x_j,y^T \tilde{x}_{j+1})$ falls in the purple or green region of the right panel of Figure~\ref{fig:region-proof}. Marginally, $(y^T x_j,y^T \tilde{x}_{j+1})\sim \mathcal{N}((\rho\tau_p,\rho^2\tau_p)^T,[(1,\rho),(\rho,1)])$. The point in purple or green region that has the smallest ellipsoid distance to $(\rho\tau_p,\rho^2\tau_p)^T$ is $(t_p,t_p)$ when $(1+\rho^2)^2r<u< (1+\rho)^2r$, thus by Lemma~\ref{lem:tool},
\[
\begin{split}
   \mathbb{P}\bigl\{W_j>t_p,A_{j+p+1\to j}\big|(\beta_j,\beta_{j+1})=(0,\tau_p)\bigr\}&\leq L_p p^{-(\sqrt{u}-\rho\sqrt{r})^2-\frac{1-\rho}{1+\rho}u}\\
   &\leq L_p p^{-r+2\sqrt{ru}-\frac{2}{1+\rho}u}\\
   &=L_p p^{-(\sqrt{u}-|\rho|\sqrt{r})^2-(\xi_{\rho}\sqrt{r}-\eta_{\rho}\sqrt{u})_+^2+(\sqrt{r}-\sqrt{u})_+^2}
\end{split}
\]
for $u\in ((1+\rho^2)^2r,(1+\rho)^2r)$, which completes the proof of (\ref{CI-2-2}) for $A=A_{j+p+1\to j}$.

    When $A_{j+1,j+p+1}$ occurs, consider the bi-variate Lasso problem $L_{j+1,j+p+1}$. In this bi-variate Lasso problem, $\{\lambda_1,\lambda_2\}=\{Z_{j+1},Z_{j+p+1}\}$, both of which are larger than $W_j$. Thus in this bi-variate Lasso problem, both variables will be picked up as signals given $W_j>t_p$. So $(y^T x_{j+1},y^T \tilde{x}_{j+1})/\sqrt{2\log(p)}$ falls in one of the four regions in the right panel of Figure~\ref{fig:region-proof} (with $x_{j+1}^T\tilde{x}_{j+1}=\rho^2$ instead of $\rho$): the purple region, the mirror of purple region against $x=y$, the green region and the mirror of green region against $x=-y$. Since $(y^T x_{j+1},y^T \tilde{x}_{j+1})\sim \mathcal{N}((\tau_p,\rho^2\tau_p)^T,[(1,\rho^2),(\rho^2,1)])$. By Lemma~\ref{lem:tool}, we need to find the point in those regions that has the smallest ellipsoid distance to the center-$(\tau_p,\rho^2\tau_p)^T$. When $\tau_p \leq t_p$, this critical point is $(y^T x_{j+1},y^T \tilde{x}_{j+1})=(t_p,t_p)$; when $\tau_p > t_p$, this critical point is $(y^T x_{j+1},y^T \tilde{x}_{j+1})=(\tau_p,t_p+\rho(\tau_p-t_p))$. So Lemma~\ref{lem:tool} gives the probability for $\lambda_1$ and $\lambda_2$ in $L_{j+1,j+p+1}$ to be both larger than $t_p$ is 
    \[L_p p^{-(\sqrt{u}-\sqrt{r})_+^2-\frac{1-\rho^2}{1+\rho^2}u}\leq L_p p^{-(\sqrt{u}-|\rho|\sqrt{r})^2-(\xi_{\rho}\sqrt{r}-\eta_{\rho}\sqrt{u})_+^2+(\sqrt{r}-\sqrt{u})_+^2}.
    \]
   Since $A_{j+1,j+p+1}\cap \{W_j>t_p\}$ implies $\{\lambda_1>t_p\}\cap \{\lambda_2>t_p\}$ in $L_{j+1,j+p+1}$, we know
    \[
          \mathbb{P}\bigl\{W_j>t_p,A_{j+1,j+p+1}\big|(\beta_j,\beta_{j+1})=(0,\tau_p)\bigr\}
          \leq L_p p^{-(\sqrt{u}-|\rho|\sqrt{r})^2-(\xi_{\rho}\sqrt{r}-\eta_{\rho}\sqrt{u})_+^2+(\sqrt{r}-\sqrt{u})_+^2}. 
    \]
Now, we have verified (\ref{CI-2-2}). Further coupled with (\ref{CI-2-3}), we have (\ref{CI-2-1}).

    \item When $(\beta_j,\beta_{j+1})=(\tau_p,0)$,
    \vspace{0.2cm}
    \beq \label{CI-3-1}
    \mathbb{P}\bigl\{W_j\leq t_p\big|(\beta_j,\beta_{j+1})=(\tau_p,0)\bigr\}\geq L_p p^{-[(\sqrt{r}-\sqrt{u})_+]^2},\eeq
    and
    \beq \label{CI-3-2}
    \mathbb{P}\bigl\{W_j\leq t_p\big|(\beta_j,\beta_{j+1})=(\tau_p,0)\bigr\}\leq L_p p^{\vartheta-f^+_{\text{Hamm}}(u,r,\vartheta)}.\eeq
    
    Let $p_3=(t_p,\rho t_p,\rho^2 t_p,\rho t_p)^T$. when $h=p_3$, variable $j$ is the first variable entering the Lasso path and $p_3$ is in the region of rejecting variable $j$ as a signal. Since $h\sim \mathcal{N}(\mu_3,G)$ with $\mu_3=(\tau_p,\rho \tau_p,\rho^2 \tau_p,\rho\tau_p)^T$, by Lemma~\ref{lem:tool}, 
     \[
     \begin{split}
        \mathbb{P}\bigl\{W_j\leq t_p\big|(\beta_j,\beta_{j+1})=(\tau_p,0)\bigr\}&\geq L_p p^{-(p_3-\mu_3)^T G^{-1}(p_3-\mu_3)/2\log(p)}\\
        &=L_p p^{-[(\sqrt{r}-\sqrt{u})_+]^2}.
     \end{split}
    \] 

Now, we show that (\ref{CI-3-2}) holds for $u\geq u^*$, which implies (\ref{CI-3-2}) for all $u\geq 0$ as discussed in the proof of EC-knockoff. 
We prove (\ref{CI-3-2}) by showing that 
    \beq
    \mathbb{P}\bigl\{W_j\leq t_p,A \big|(\beta_j,\beta_{j+1})=(\tau_p,0)\bigr\}\leq L_p p^{-[(\sqrt{r}-\sqrt{u})_+]^2}\eeq
holds for $A=A_j,A_{j+1},A_{j+p},A_{j+p+1}$ given $u\geq u^*$. Respectively, 
\[
\begin{split}
    \mathbb{P}\bigl\{W_j\leq t_p,A_j \big|(\beta_j,\beta_{j+1})=(\tau_p,0)\bigr\}
&\leq \mathbb{P}\bigl\{|y^Tx_j|\leq t_p\big|(\beta_j,\beta_{j+1})=(\tau_p,0)\bigr\}\\
&=L_p p^{-[(\sqrt{r}-\sqrt{u})_+]^2},
\end{split}
\]
and by symmetry and (\ref{ueq}),
\[
\begin{split}
 \mathbb{P}\bigl\{W_j\leq t_p,A_{j+1} \big|(\beta_j,\beta_{j+1})=(\tau_p,0)\bigr\}=\mathbb{P}\bigl\{W_j\leq t_p,A_{j+p+1} \big|(\beta_j,\beta_{j+1})=(\tau_p,0)\bigr\}   \\
  \leq \mathbb{P}\bigl\{|y^Tx_j|\leq |y^Tx_{j+p+1}|\big|(\beta_j,\beta_{j+1})=(\tau_p,0)\bigr\} \leq L_p p^{-\frac{1-\rho}{2}r}\leq L_p p^{-[(\sqrt{r}-\sqrt{u})_+]^2},
\end{split}
\]
\[
\begin{split}
\mathbb{P}\bigl\{W_j\leq t_p,A_{j+p} \big|(\beta_j,\beta_{j+1})=(\tau_p,0)\bigr\}  
  &\leq \mathbb{P}\bigl\{|y^Tx_j|\leq |y^Tx_{j+p}|\big|(\beta_j,\beta_{j+1})=(\tau_p,0)\bigr\} \\
  &\leq L_p p^{-\frac{1-\rho^2}{2}r} \leq L_p p^{-[(\sqrt{r}-\sqrt{u})_+]^2}.
\end{split}
\]
(\ref{CI-3-2}) is immediate by $[(\sqrt{r}-\sqrt{u})_+]^2 \geq f^+_{\text{Hamm}}(r,u,\vartheta)-\vartheta$.
\vspace{0.2cm}

    \item When $(\beta_j,\beta_{j+1})=(\tau_p,\tau_p)$,
    \vspace{0.2cm}
    \beq \label{CI-4-1}
    \mathbb{P}\bigl\{W_j\leq t_p\big|(\beta_j,\beta_{j+1})=(\tau_p,\tau_p)\bigr\}\leq L_p p^{\vartheta-f^+_{\text{Hamm}}(u,r,\vartheta)}.\eeq
    We prove (\ref{CI-4-1}) by showing 
\beq \label{CI-4-1__}
    \mathbb{P}\bigl\{W_j\leq t_p, A\big|(\beta_j,\beta_{j+1})=(\tau_p,\tau_p)\bigr\}\leq L_p p^{-[(\sqrt{r}-\sqrt{u})_+]^2}\eeq
holds for $A=A_j,A_{j+1\to j},A_{j+1\to j+p},A_{j+1\to j+p+1},A_{j+p},A_{j+p+1}$ given $u\geq u^*$, which cover all possibilities.
Respectively,
\[ 
\begin{split}
    \mathbb{P}\bigl\{W_j\leq t_p, A_j\big|(\beta_j,\beta_{j+1})=(\tau_p,\tau_p)\bigr\}& \leq \mathbb{P}\bigl\{|y^T x_j|\leq t_p\big|(\beta_j,\beta_{j+1})=(\tau_p,\tau_p)\bigr\}\\
    &\leq L_p p^{-[((1+\rho)\sqrt{r}-\sqrt{u})_+]^2} \leq L_p p^{-[(\sqrt{r}-\sqrt{u})_+]^2},
\end{split}
\]
\[ 
\begin{split}
    \mathbb{P}\bigl\{W_j\leq t_p, A_{j+p}\big|(\beta_j,\beta_{j+1})=(\tau_p,\tau_p)\bigr\}& \leq \mathbb{P}\bigl\{|y^T x_j|\leq |y^T \tilde{x}_{j}|\big|(\beta_j,\beta_{j+1})=(\tau_p,\tau_p)\bigr\}\\
    &\leq L_p p^{-\frac{1-\rho^2}{2}r} \leq L_p p^{-[(\sqrt{r}-\sqrt{u})_+]^2},
\end{split}
\]
\[ 
\begin{split}
    \mathbb{P}\bigl\{W_j\leq t_p, A_{j+p+1}\big|(\beta_j,\beta_{j+1})=(\tau_p,\tau_p)\bigr\}& \leq \mathbb{P}\bigl\{|y^T x_j|\leq |y^T \tilde{x}_{j+1}|\big|(\beta_j,\beta_{j+1})=(\tau_p,\tau_p)\bigr\}\\
    &\leq L_p p^{-\frac{(1-\rho)(1+\rho)^2}{2}r} \leq L_p p^{-[(\sqrt{r}-\sqrt{u})_+]^2}.
\end{split}
\]
When $A_{j+1\to j}$ occurs, the bi-variate Lasso problem $L_{j,j+1}$ has variable $j$ is a false negative given $W_j\leq t_p$, which implies $\mathbb{P}\bigl\{W_j\leq t_p,A_{j+1\to j}\big|(\beta_j,\beta_{j+1})=(\tau_p,\tau_p)\bigr\}$ is upper bounded by the corresponding false negative rate of Lasso, which is $L_p p^{-(\xi_\rho\sqrt{r}-\eta_\rho \sqrt{u})_+^2}\leq L_p p^{-[(\sqrt{r}-\sqrt{u})_+]^2}$ for $u\geq u^*$.

When $A_{j+1\to j+p}$ occurs, we know variable $j+p$ instead of variable $j$ is the second one entering the Lasso path. This means the $\lambda_2$ (the $\lambda$ value when the second variable entering Lasso path) of the bi-variate Lasso problem $L_{j+1,j+p}$ is larger than the $\lambda_2$ of the bi-variate Lasso problem $L_{j,j+1}$. When $y^T x_{j+1}\geq 0$, we must have 
\[\max\{\frac{y^T x_j-\rho y^T x_{j+1}}{1-\rho},\frac{y^T x_j-\rho y^T x_{j+1}}{-1-\rho}\}< \max\{\frac{y^T \tilde{x}_j-\rho y^T x_{j+1}}{1-\rho},\frac{y^T \tilde{x}_j-\rho y^T x_{j+1}}{-1-\rho}\}.\]
Therefore, $A_{j+1\to j+p}$ implies one the three following events must occur:
\[y^T x_{j+1}< 0,\frac{y^T x_j-\rho y^T x_{j+1}}{1-\rho}<\frac{y^T \tilde{x}_j-\rho y^T x_{j+1}}{1-\rho}, \frac{y^T x_j-\rho y^T x_{j+1}}{1-\rho}<\frac{y^T \tilde{x}_j-\rho y^T x_{j+1}}{-1-\rho}\] 
The probability of these three events given $(\beta_j,\beta_{j+1})=(\tau_p,\tau_p)$ are $L_p p^{-(1+\rho)^2r}$, $L_p p^{-\frac{1-\rho^2}{2}r}$ and $L_p p^{-\frac{(1+\rho)^3(1-\rho)}{2(1+\rho^2)}r}$, all of which are upper bounded by $L_p p^{-[(\sqrt{r}-\sqrt{u})_+]^2}$ when $u\geq u^*$.

When $A_{j+1\to j+p+1}$ occurs, the $\lambda_2$ of the bi-variate Lasso problem $L_{j+1,j+p+1}$ is larger than the $\lambda_2$ of the bi-variate Lasso problem $L_{j,j+1}$. When $y^T x_{j+1}\geq 0$, we must have 
\[\max\{\frac{y^T x_j-\rho y^T x_{j+1}}{1-\rho},\frac{y^T x_j-\rho y^T x_{j+1}}{-1-\rho}\}< \max\{\frac{y^T \tilde{x}_{j+1}-\rho^2 y^T x_{j+1}}{1-\rho^2},\frac{y^T \tilde{x}_{j+1}-\rho^2 y^T x_{j+1}}{-1-\rho^2}\}.\]
Therefore, $A_{j+1\to j+p+1}$ implies one the three following events must occur:
\[y^T x_{j+1}< 0,\frac{y^T x_j-\rho y^T x_{j+1}}{1-\rho}<\frac{y^T \tilde{x}_{j+1}-\rho^2 y^T x_{j+1}}{1-\rho^2}, \frac{y^T x_j-\rho y^T x_{j+1}}{1-\rho}<-\frac{y^T \tilde{x}_{j+1}-\rho^2 y^T x_{j+1}}{-1-\rho^2}.\]
Respectively, the probability of these three events are $L_p p^{-(1+\rho)^2r}$, $L_p p^{-\frac{1-\rho^2}{2}r}$ and $L_p p^{-\frac{(1+\rho)^3(1-\rho)}{2(1+\rho^2)}r}$, all of which are upper bounded by $L_p p^{-[(\sqrt{r}-\sqrt{u})_+]^2}$ when $u\geq u^*$. From here we have verified (\ref{CI-4-1__}), thus implies (\ref{CI-4-1}).

From (\ref{CI-1-1}), (\ref{CI-2-1}), (\ref{CI-3-1}), (\ref{CI-3-2}) and (\ref{CI-4-1}), we have 
\beq 
\mathbb{P}\bigl\{W_j> t_p,\beta_j=0\bigr\}+\mathbb{P}\bigl\{W_j\leq t_p,\beta_j=\tau_p\bigr\}= L_p p^{-f^+_{\text{Hamm}}(r,u,\vartheta)},
\eeq
which completes the proof for positive $\rho$.

    \item When $(\beta_j,\beta_{j+1})=(-\tau_p,\tau_p)$,
    \vspace{0.2cm}
    \beq \label{CI-5-1}
    \mathbb{P}\bigl\{W_j\leq t_p\big|(\beta_j,\beta_{j+1})=(-\tau_p,\tau_p)\bigr\}\geq L_p p^{-((\xi_\rho \sqrt{r}-\eta_\rho^{-1} \sqrt{u})_+)^2}\eeq
    and
    \beq \label{CI-5-3}
    \mathbb{P}\bigl\{W_j\leq t_p\big|(\beta_j,\beta_{j+1})=(-\tau_p,\tau_p)\bigr\}\leq L_p p^{-\min\{((\xi_\rho \sqrt{r}-\eta_\rho^{-1} \sqrt{u})_+)^2,\frac{(1-\rho)^3(1+\rho)}{2(1+\rho^2)}r\}}.\eeq
 Let 
 \[p_4^T=
\left\{  
\begin{array}{ll}
(-(1-\rho)\tau_p,(1-\rho)\tau_p,\rho (1-\rho)\tau_p,-\rho(1-\rho) \tau_p), & (1-\rho)\tau_p\leq t_p, \\
(\rho(1-\rho)\tau_p-(1+\rho)t_p,(1-\rho)\tau_p,\rho(1-\rho)\tau_p,\rho^2(1-\rho)\tau_p-\rho(1+\rho)t_p), & (1-\rho)\tau_p> t_p. \\
\end{array}
\right.
\]
When $h=p_4$ and $(1-\rho)\tau_p\leq t_p$, variable $j$ is the first variable entering the Lasso path with $W_j=(1-\rho)\tau_p\leq t_p$; when $h=p_4$ and $(1-\rho)\tau_p> t_p$, $j+1$ is the first and $j$ is the second variable entering the Lasso path with $W_j=t_p$. Regardless of the relationship between $\tau_p$ and $t_p$, $h=p_4$ is always in the region of rejecting $j$ as a signal. Since $h\sim \mathcal{N}(\mu_4,G)$ with $\mu_4=(-(1-\rho)\tau_p,(1-\rho)\tau_p,\rho (1-\rho)\tau_p,-\rho(1-\rho) \tau_p)^T$, by Lemma~\ref{lem:tool}, 
     \[
     \begin{split}
        \mathbb{P}\bigl\{W_j\leq t_p\big|(\beta_j,\beta_{j+1})=(-\tau_p,\tau_p)\bigr\}&\geq L_p p^{-(p_4-\mu_4)^T G^{-1}(p_4-\mu_4)/2\log(p)}\\
        &=L_p p^{-((\xi_\rho \sqrt{r}-\eta_\rho^{-1} \sqrt{u})_+)^2}.
     \end{split}
    \] 
Next, we show that 
    \beq \label{CI-5-4}
    \mathbb{P}\bigl\{W_j\leq t_p, A\big|(\beta_j,\beta_{j+1})=(-\tau_p,\tau_p)\bigr\}\leq L_p p^{-\min\{((\xi_\rho \sqrt{r}-\eta_\rho^{-1} \sqrt{u})_+)^2,\frac{(1-\rho)^3(1+\rho)}{2(1+\rho^2)}r\}}.\eeq
holds for $A=A_j,A_{j+1\to j},A_{j+1\to j+p},A_{j+1\to j+p+1},A_{j+p},A_{j+p+1}$, which cover all possibilities.

When $A=A_j$ or $A_{j+1 \to j}$ occurs, as previously discussed, variable $j$ is a false negative in the bi-variate Lasso $L_{j,j+1}$ given $W_j\leq t_p$, which implies $\mathbb{P}\bigl\{W_j\leq t_p,A\big|(\beta_j,\beta_{j+1})=(-\tau_p,\tau_p)\bigr\}$ is upper bounded by the corresponding false negative rate of Lasso, which is $L_p p^{-((\xi_\rho \sqrt{r}-\eta_\rho^{-1} \sqrt{u})_+)^2}$. 

When $A_{j+1\to j+p}$ occurs, the $\lambda_2$ of the bi-variate Lasso problem $L_{j+1,j+p}$ is larger than the $\lambda_2$ of the bi-variate Lasso problem $L_{j,j+1}$. When $y^T x_{j+1}\geq 0$, we must have 
\[\max\{\frac{y^T x_j-\rho y^T x_{j+1}}{1-\rho},\frac{y^T x_j-\rho y^T x_{j+1}}{-1-\rho}\}< \max\{\frac{y^T \tilde{x}_j-\rho y^T x_{j+1}}{1-\rho},\frac{y^T \tilde{x}_j-\rho y^T x_{j+1}}{-1-\rho}\}.\]
Therefore, $A_{j+1\to j+p}$ implies one of the three following events must occur:
\[y^T x_{j+1}+y^T \tilde{x}_j< 0,\frac{y^T x_j-\rho y^T x_{j+1}}{-1-\rho}<\frac{y^T \tilde{x}_j-\rho y^T x_{j+1}}{1-\rho}, \frac{y^T x_j-\rho y^T x_{j+1}}{-1-\rho}<\frac{y^T \tilde{x}_j-\rho y^T x_{j+1}}{-1-\rho}\] 
The probability of these three events are $L_p p^{-\frac{(1+\rho)(1-\rho)^2}{2}r}$, $L_p p^{-\frac{(1-\rho)^3(1+\rho)}{2(1+\rho^2)}r}$ and $L_p p^{-\frac{1-\rho^2}{2}r}$, all of which are upper bounded by $L_p p^{-\frac{(1-\rho)^3(1+\rho)}{2(1+\rho^2)}r}$.

When $A_{j+1\to j+p+1}$ occurs, the $\lambda_2$ of the bi-variate Lasso problem $L_{j+1,j+p+1}$ is larger than the $\lambda_2$ of the bi-variate Lasso problem $L_{j,j+1}$. When $y^T x_{j+1}\geq 0$, we must have 
\[\max\{\frac{y^T x_j-\rho y^T x_{j+1}}{1-\rho},\frac{y^T x_j-\rho y^T x_{j+1}}{-1-\rho}\}< \max\{\frac{y^T \tilde{x}_{j+1}-\rho^2 y^T x_{j+1}}{1-\rho^2},\frac{y^T \tilde{x}_{j+1}-\rho^2 y^T x_{j+1}}{-1-\rho^2}\}.\]
Therefore, $A_{j+1\to j+p+1}$ implies one of the three following events must occur:
\[y^T x_{j+1}+y^T \tilde{x}_j< 0,\frac{y^T x_j-\rho y^T x_{j+1}}{-1-\rho}<\frac{y^T \tilde{x}_{j+1}-\rho^2 y^T x_{j+1}}{1-\rho^2}, \frac{y^T x_j-\rho y^T x_{j+1}}{-1-\rho}<\frac{y^T \tilde{x}_{j+1}-\rho^2 y^T x_{j+1}}{-1-\rho^2}.\]
The probability of these three events are $L_p p^{-\frac{(1+\rho)(1-\rho)^2}{2}r}$, $L_p p^{-\frac{1-\rho^2}{2}r}$, $L_p p^{-\frac{(1-\rho)^3(1+\rho)}{2(1+\rho^2)}r}$, all of which are upper bounded by $L_p p^{-\frac{(1-\rho)^3(1+\rho)}{2(1+\rho^2)}r}$.

When $A_{j+p}$ occurs, if $y^T \tilde{x}_j<0$, then $y^T x_{j+1}+y^T \tilde{x}_j\leq 0$, which happens with probability $L_p p^{-\frac{(1+\rho)(1-\rho)^2}{2}r}\leq L_p p^{-\frac{(1-\rho)^3(1+\rho)}{2(1+\rho^2)}r}$. If $y^T \tilde{x}_j\geq 0$, then $y^T \tilde{x}_j+\frac{1-\rho}{2}y^T x_j-\frac{1+\rho}{2}y^T x_{j+1}\geq 0$, which happens with probability $L_p p^{-\frac{2(1-\rho)^3}{3+\rho^2}r}\leq L_p p^{-\frac{(1-\rho)^3(1+\rho)}{2(1+\rho^2)}r}$. Therefore, (\ref{CI-5-4}) holds for $A_{j+p}$ and also for $A_{j+p+1}$ due to symmetry. We thus complete the proof for (\ref{CI-5-4}).

Now consider the case where $\beta_j$ takes value in $\{0,-\tau_p\}$ and $\beta_{j+1}$ takes value in $\{0,\tau_p\}$, this corresponds to the $\rho<0$ case (we flipped the sign of $\rho$ and $\beta_j$ simultaneously).
By (\ref{CI-1-1}), (\ref{CI-2-1}), (\ref{CI-3-1}) and (\ref{CI-5-1}), we know 
\beq 
\begin{split}
    \mathbb{P}\bigl\{W_j> t_p,\beta_j=0\bigr\}+\mathbb{P}\bigl\{W_j\leq t_p,\beta_j=-\tau_p\bigr\}\\
\geq L_p p^{-\min\{f^+_{\text{Hamm}}(u,r,\vartheta),2\vartheta+((\xi_\rho \sqrt{r}-\eta_\rho^{-1} \sqrt{u})_+)^2\}}.
\end{split}
\eeq
Meanwhile, (\ref{CI-1-1}), (\ref{CI-2-1}), (\ref{CI-3-2}) and (\ref{CI-5-3}) gives
\beq 
\begin{split}
    \mathbb{P}\bigl\{W_j> t_p,\beta_j=0\bigr\}+\mathbb{P}\bigl\{W_j\leq t_p,\beta_j=-\tau_p\bigr\}\\
\leq L_p p^{-\min\{f^+_{\text{Hamm}}(u,r,\vartheta),2\vartheta+((\xi_\rho \sqrt{r}-\eta_\rho^{-1} \sqrt{u})_+)^2,2\vartheta+\frac{(1-|\rho|)^3(1+|\rho|)}{2(1+|\rho|^2)}r\}}.
\end{split}
\eeq
The proof is complete once we show that 
\beq \label{reductio}
\min\{f^+_{\text{Hamm}}(u,r,\vartheta),2\vartheta+((\xi_\rho \sqrt{r}-\eta_\rho^{-1} \sqrt{u})_+)^2\}\leq 2\vartheta+\frac{(1-|\rho|)^3(1+|\rho|)}{2(1+\rho^2)}r.
\eeq 
Otherwise, there exists a tuple of $(\vartheta,r,\rho,u,r)$ such that
\begin{equation}\label{CI-un1}
    2\vartheta+\frac{(1-|\rho|)^3(1+|\rho|)}{2(1+\rho^2)}r< 2\vartheta+((\xi_\rho \sqrt{r}-\eta_\rho^{-1} \sqrt{u})_+)^2
\end{equation}
and 
\begin{equation}\label{CI-un2}
    2\vartheta+\frac{(1-|\rho|)^3(1+|\rho|)}{2(1+\rho^2)}r< \vartheta+(\sqrt{u}-|\rho|\sqrt{r})^2+((\xi_\rho \sqrt{r}-\eta_\rho \sqrt{u})_+)^2-((\sqrt{r}-\sqrt{u})_+)^2
\end{equation}
are satisfied simultaneously.

By (\ref{CI-un1}), $\xi_\rho \sqrt{r}-\eta_\rho^{-1} \sqrt{u}>0$, which implies $(1-|\rho|)\sqrt{r}>\sqrt{u}$. Therefore, the right hand side of (\ref{CI-un2}) simplifies to $\vartheta+\frac{1-|\rho|}{1+|\rho|}u$. By (\ref{CI-un2}), we know
\[\frac{(1-|\rho|)^3(1+|\rho|)}{2(1+\rho^2)}r\leq \vartheta+\frac{(1-|\rho|)^3(1+|\rho|)}{2(1+\rho^2)}r<\frac{1-|\rho|}{1+|\rho|}u.\]
Plug this into the right hand side of (\ref{CI-un1}), we have
\begin{equation}
\begin{split}
        2\vartheta+\frac{(1-|\rho|)^3(1+|\rho|)}{2(1+\rho^2)}r< 2\vartheta+((\xi_\rho \sqrt{r}-\eta_\rho^{-1} \sqrt{u})_+)^2\\
        \leq 2\vartheta + \Big(\sqrt{1-\rho^2}-\sqrt{\frac{(1-|\rho|)(1+|\rho|)^3}{2(1+\rho^2)}}\Big)^2r,
\end{split}
\end{equation}
which can only be true when $\rho^2>1$. By reduction, we proved (\ref{reductio}).

\end{itemize}

\section{Proof of Theorem~\ref{thm:OLS}}

The least-squares estimator satisfies that $\hat{\beta}\sim {\cal N}_p(\beta, G^{-1})$. It gives $\hat{\beta}_j\sim {\cal N}(\beta_j, \omega_j)$. Applying Lemma~\ref{lem:tool} to $X_p=\hat{\beta}_j$ and $S=\{x\in\mathbb{R}: x\geq \sqrt{u}\}$, we have 
\[
    \mathbb{P}(|\hat{\beta}_j|>t_p(u)|\beta_j=0) = L_pp^{-\omega_j^{-1}u},\qquad
    \mathbb{P}(|\hat{\beta}_j|\leq t_p(u)|\beta_j=\tau_p) = L_pp^{-\omega_j^{-1}(\sqrt{r}-\sqrt{u})_+^2}.
\]
It follows that
\begin{align*}
    \FP_p(u)& =\sum_{j=1}^p (1-\epsilon_p) \cdot \mathbb{P}(W_j^*>t_p(u)|\beta_j=0)=L_p\sum_{j=1}^p p^{-\omega_j^{-1}u},\cr
     \FN_p(u)& =\sum_{j=1}^p \epsilon_p\cdot \mathbb{P}(W_j^*<t_p(u)|\beta_j=\tau_p)=L_pp^{-\vartheta}\sum_{j=1}^p p^{-\omega_j^{-1}(\sqrt{r}-\sqrt{u})_+^2}.
\end{align*}
For the block-wise diagonal design \eqref{block}, $\omega_j=(1-\rho^2)^{-1}$ for all $1\leq j\leq p-1$.

\section{Proof of Theorem~\ref{thm:knockoffOLS-block}} 
By the property of least-square coefficients, \[(\hat{\beta}_1,\cdots,\hat{\beta}_p,\tilde{\beta}_1,\cdots,\tilde{\beta}_p)\sim \mathcal{N}_{2p}\big((\beta_1,\cdots,\beta_p,0,\cdots,0),(G^*)^{-1}\big).\]
Consider the joint distribution of $\hat{\beta}_j$ and $\tilde{\beta}_{j}$ which are the regression coefficient of $x_j$ and $\tilde{x}_j$, we know that $(\hat{\beta}_j,\tilde{\beta}_{j})\sim \mathcal{N}_2\big((\beta_j,0),A_j\big)$ where $A_j$ has $\omega_{1j}$ as its diagonal element and $\omega_{2j}$ as its off-diagonal elements. Then theorem~\ref{thm:knockoffOLS-block} is immediate from the following lemma: 
\begin{lem}\label{lem:ols}
If $(Z_j,\tilde{Z}_j)$ follows ${\cal N}_2\Bigl((\beta_j,0)^T,\;\;\Sigma\Bigr)$ with $\Sigma=((\sigma_1,\sigma_2),(\sigma_2,\sigma_1))$, then \begin{equation}\label{ols_lemma1}
    \mathbb{P}(|Z_j|>\sqrt{2u\log(p)}, |Z_j|\geq |\tilde{Z}_j|\big|\beta_j=0)=L_p p^{-u/\sigma_1}
\end{equation} 
and 
\begin{equation}\label{ols_lemma2}
\begin{split}
    &\mathbb{P}(|Z_j|\leq \sqrt{2u\log(p)}\text{ or }|Z_j|< |\tilde{Z}_j|\big|\beta_j=\sqrt{2r\log(p)})\\
    =&L_p p^{-\min\{(\sqrt{r}-\sqrt{u})_+^2/\sigma_1,r/(2\max\{\sigma_1+\sigma_2,\sigma_1-\sigma_2\})\}}.
\end{split}
\end{equation}
\end{lem}

Next, we prove Lemma~\ref{lem:ols}. To compute the left hand side of (\ref{ols_lemma1}), we only need to find the $t$ such that ellipsoid $(x,y)\Sigma^{-1}(x,y)^T=t^2$ is tangent with $x=\pm\sqrt{2u\log(p)}$. This is because when we increase the radius of the ellipsoid, it must intersect with $x=\pm\sqrt{2u\log(p)}$ first amongst the boundaries of the region that pick variable $j$ as a signal. When they intersect, 
\[t^2=\frac{1}{\sigma_1^2-\sigma_2^2}(\sigma_1 x^2-2\sigma_2 xy +\sigma_1 y^2)=\frac{1}{\sigma_1^2-\sigma_2^2} \Big(\sigma_1\Big(y-\frac{\sigma_2}{\sigma_1}x\Big)^2+\Big(\sigma_1-\frac{\sigma_2^2}{\sigma_1}\Big)x^2\Big)\geq \frac{2u\log(p)}{\sigma_1}.\]
When $t^2=\frac{2u\log(p)}{\sigma_1}$, the tangent points are $(\pm\sqrt{2u\log(p)},\pm\frac{\sigma_2}{\sigma_1}\sqrt{2u\log(p)})$. By Lemma~\ref{lem:tool}, we verified (\ref{ols_lemma1}).

For (\ref{ols_lemma2}), when $r<u$, the center of the bi-variate normal is in the region of rejecting variable $j$ as a signal thus the false positive rate is $L_p$. When $r>u$, we need to find the $t$ such that ellipsoid $(x-\beta_j,y)\Sigma^{-1}(x-\beta_j,y)^T=t^2$ is tangent with either $x=\pm\sqrt{2u\log(p)}$ or $y=\pm x$. When the ellipsoid intersects with $x=\pm\sqrt{2u\log(p)}$,
\[t^2=\frac{1}{\sigma_1^2-\sigma_2^2} \Big(\sigma_1\Big(y-\frac{\sigma_2}{\sigma_1}(x-\beta_j)\Big)^2+\Big(\sigma_1-\frac{\sigma_2^2}{\sigma_1}\Big)(x-\beta_j)^2\Big)\geq \frac{2(\sqrt{u}-\sqrt{r})^2\log(p)}{\sigma_1},\]
therefore, they are tangent at $(\pm \sqrt{2u\log(p)},\frac{\sigma_2}{\sigma_1}(\pm \sqrt{2u\log(p)}-\beta_j))$ when $t^2=\frac{2(\sqrt{u}-\sqrt{r})^2\log(p)}{\sigma_1}$.

Meanwhile, since the long/short shaft of the ellipsoid are paralleled with $y=\pm x$, the tangent points of ellipsoid with $y=\pm x$ must be $(\beta_j/2,\beta_j/2)$ and $(\beta_j/2,-\beta_j/2)$, which gives $t^2=\frac{r\log(p)}{\sigma_1+\sigma_2}$ and $\frac{r\log(p)}{\sigma_1-\sigma_2}$. From here we can conclude the "distance" between the center of the normal distribution and the region that reject variable $j$ as a signal is \[\min\{\frac{2(\sqrt{r}-\sqrt{u})_+^2\log(p)}{\sigma_1},\frac{r\log(p)}{\sigma_1+\sigma_2},\frac{r\log(p)}{\sigma_1-\sigma_2}\}.\]
By Lemma~\ref{lem:tool}, we know 
\[
    \mathbb{P}(|Z_j|\leq \sqrt{2u\log(p)}\big|\beta_j=\sqrt{2r\log(p)})=L_p p^{-\min\{(\sqrt{r}-\sqrt{u})_+^2/\sigma_1,r/(2\max\{\sigma_1+\sigma_2,\sigma_1-\sigma_2\})\}}.
\]

\vskip 0.2in
\bibliography{VS}

\end{document}